\documentclass[runningheads,envcountsect]{svjour3}

\usepackage{amsmath,amssymb,mathptmx}
\usepackage{cite}
\usepackage[width=114mm, left=48mm, height=18.6cm] {geometry}
\usepackage{exscale}
\usepackage{color}
\usepackage[weather,alpine,misc,geometry]{ifsym}
\usepackage{hyperref}
\hypersetup{
    colorlinks=true,
    linkcolor=blue, 
    citecolor=blue, 
    urlcolor=blue  } 
\smartqed  
\newcommand{\al}{\alpha}

\newcommand{\la}{\lambda}
\newcommand{\de}{\delta}
\newcommand{\eps}{\varepsilon}
\newcommand{\bx}{\bar x}

\newcommand{\iv}{^{-1} }

\newcommand {\R} {\mathbb R}
\newcommand {\N} {\mathbb N}

\newcommand {\B} {\mathbb B}

\newcommand {\dom} {{\rm dom}\,}
\newcommand {\epi} {{\rm epi}\,}

\newcommand {\bd} {{\rm bd}\,}

\newcommand {\sd} {\partial}

\renewcommand{\iff}{$ \Leftrightarrow\ $}
\newcommand{\folgt}{$ \Rightarrow\ $}

\def\nbh{neighbourhood}
\def\es{\emptyset}

\def\RHS{right-hand side}
\def\SVM{set-valued mapping}
\def\EVP{Ekeland variational principle}
\def\Fr{Fr\'echet}

\newcommand{\norm}[1]{\left\Vert#1\right\Vert}

\newcommand{\set}[1]{\left\{#1\right\}}

\newcommand{\blue}[1]{\textcolor{blue}{#1}}

\newcommand{\red}[1]{\textcolor{red}{#1}}

\newcommand{\ang}[1]{\left\langle #1 \right\rangle}
\newcommand{\qdtx}[1]{\quad\mbox{#1}\quad}
\newcommand{\AND}{\quad\mbox{and}\quad}
\newcounter{mycount}
\def\cnta{\setcounter{mycount}{\value{enumi}}}
\def\cntb{\setcounter{enumi}{\value{mycount}}}
\newcommand{\AK}[1]{\todo[inline]{AK {#1}}}

\newcommand{\Tr}{{\rm tr} [\Omega_1,\ldots,\Omega_n](\bx)}

\newcommand{\HB}[1]{\todo[inline,color=green!40]{HB {#1}}}

\spnewtheorem{Remark}{Remark}[section]{\bf}{\rm}
\spnewtheorem{Example}{Example}[section]{\bf}{\rm}
\spnewtheorem*{Proof}{Proof.}{\bf}{\rm}

\title{
Extremality, Stationarity and Generalized Separation of Collections of Sets
}

\author{
Hoa T. Bui
\and
Alexander Y. Kruger
}

\institute{
Hoa T. Bui
\and
Alexander Y. Kruger (\Letter\,)
\at
Centre for Informatics and Applied Optimization,
School of Science, Engineering and Information Technology,
Federation University Australia, POB 663, Ballarat, VIC 3350, Australia\\
\email{h.bui@federation.edu.au, a.kruger@federation.edu.au}
}

\date{Received: date / Accepted: date}
\begin{document}

\maketitle
\abstract{
The core arguments used in various proofs of the \emph{extremal principle} and its extensions as well as in primal and dual characterizations of \emph{approximate stationarity} and \emph{transversality} of collections of sets are exposed, analysed and refined, leading to a unifying theory, encompassing all existing approaches to obtaining `extremal' statements.
For that, we examine and clarify quantitative relationships between the parameters involved in the respective definitions and statements.
Some new characterizations of extremality properties are obtained.
}



\section{Introduction}

Models involving pairs, or more generally finite (and even infinite) collections of sets, are pretty common in various fields of mathematics, especially optimization and neighbouring areas.
For instance, the classical \emph{separation theorem} \cite{Roc70,Zal02,Pen16} for convex sets, one of the key results of nonlinear analysis, is instrumental when establishing multiplier rules or subdifferential calculus rules;
the ubiquitous \emph{feasibility problems} \cite{BauBor96,KruLukTha18} cover solving systems of (generalized) equations, while appropriate \emph{transversality} (regular intersection) conditions produce \emph{constraint qualifications} and convergence estimates for \emph{alternating projections} (\emph{von Neuman method}) \cite{DruIofLew15,BauDaoNolPha16,KruTha16,NolRon16}, or more generally \emph{cyclic projections}, which in turn provide a convenient model for convergence analysis of computational algorithms.
On the other hand, (local) \emph{optimality} in optimization problems can naturally be interpreted as `irregular' (extremal) intersection of certain sets, and sufficient transversality conditions immediately produce necessary optimality conditions.

In the framework of convex analysis, the property of two sets that one set does not meet the interior of the other set (while both are nonempty), which is the key assumption of the conventional separation theorem, provides an example of the absence of regular intersection.
Thus, the separation theorem can be interpreted as a characterization of the absence of regular intersection of convex sets.
On the other hand, the opposite property, i.e. when the intersection of one set with the interior of the other set is nonempty, is a typical qualification condition in various convex analysis statements \cite{Zal02}.

There have been many successful applications of convex analysis, and particularly the separation theorem in the nonconvex settings by considering appropriate local convex approximations of sets and functions.
The most prominent example in the optimization area is probably given by the \emph{Dubovitskii--Milyutin formalism} \cite{DubMil65}.
The Clarke tangent cone \eqref{TCC} is an important example of a local convex approximation of a set.

However, in many situations it is impossible to construct satisfactory local convex approximations of nonconvex sets.
For instance, if a set consists of two intersecting lines on the plain, it is easy to check that its Clarke tangent cone at the point of intersection is trivial (contains only the zero vector), and provides no meaningful information about the set.
The powerful tools of convex analysis generally fail outside of the comfortable convex setting.

A new tool with the potential to substitute the conventional separation theorem in the general nonconvex environment---currently known as the \emph{Extremal principle}---was suggested in 1979--1980 in \cite{KruMor79,KruMor80.2,KruMor80}.
It was observed that the separation theorem actually characterizes a kind of extremal arrangement of sets.
Building on this observation,
the key assumption of the conventional separation theorem, that one set does not meet the interior of the other set, was replaced by a more general geometric (local) \emph{extremality} property that arbitrarily small translations of the sets make their intersection empty (in a neighbourhood of a given \emph{extremal point}); cf. Definition~\ref{D2.5}(i) and (ii).
This property assumes neither convexity of the sets nor that one of the sets has nonempty interior, while still embracing many conventional (and generalized) optimality notions.
Moreover, it is applicable to any finite ($n\ge2$) collections of sets (cf. \cite{KruMor79,KruMor80,KruMor80.2}) and, if the space is Banach and the sets are closed, the function measuring the distance between points in different translated sets satisfies the assumptions of the \EVP.
Employing an appropriate sum rule to the perturbation of the distance function arising from the application of the \EVP, one can formulate dual characterizations of extremality in terms of appropriate normals to individual sets.
Such dual conditions can be interpreted as a kind of generalized separation; cf. conditions (ii) and (iii) in Theorem~\ref{T2.9}.

All existing infinite dimensional versions of the extremal principle are based on the \EVP\ and differ mainly in the kind of sum rule used in the proof and, accordingly, the kind of normals used in the statement.
In the beginning, the statements were restricted to \emph{Fr\'echet smooth} spaces (i.e. Banach spaces admitting an equivalent norm, Fr\'echet differentiable away from zero); cf. \cite{KruMor80,KruMor80.2,Kru81.2,Kru85.1}.
This allowed application of the \emph{differentiable sum rule} (Lemma~\ref{SR}(ii)); the generalized separation conditions were formulated in terms of sets of $\eps$-normals \cite{Kru81.1}.
The Asplund space \emph{approximate sum rule} for Fr\'echet subdifferentials (Lemma~\ref{SR}(iii)) established by Fabian \cite{Fab89} in the end of the 1980s allowed extending the extremal principle to Asplund spaces and replacing sets of $\eps$-normals in the generalized separation conditions with the \Fr\ normal cones.
This Asplund space version of the extremal principle appeared in 1996 in \cite{MorSha96}.
It was also shown in \cite{MorSha96} that with the \Fr\ normal cones the result cannot be in general extended beyond Asplund spaces.

Similarly, plugging the \emph{Clarke--Rockafellar subdifferential sum rule} (Lemma~\ref{SR}\allowbreak(iv)) into the proof, one can formulate a general Banach space version of the extremal principle in terms of the Clarke normal cones.
Of course, since the \Fr\ normal cone is a subset of the Clarke one and the latter is often unreasonably large, when the space is Asplund, the version in terms of \Fr\ normals gives better estimates, unless the sets are convex.
In the latter case, the estimates coincide.
Several publications have gone a little further and considered abstract subdifferentials possessing appropriate sum rules (in the respective \emph{trustworthy} \cite{Iof98} spaces) and the corresponding normal cones; see, e.g., \cite{Iof98,BorJof98,Mor06.1,KruLop12.1}.

Since its inception almost 40 years ago,
the extremal principle has indeed impacted strongly on the nonconvex analysis and optimization substituting the conventional separation theorem.
Several examples of its application to proving necessary optimality conditions (multiplier rules) can be found already in the very first publications \cite{KruMor79,KruMor80,KruMor80.2}.
Numerous applications of the extremal principle in multiple publications on optimization and analysis are exposed and commented on in the monograph \cite{Mor06.1}.
The original publications have been followed by several studies of the concept of extremality of collections of sets and the extremal principle, resulting in its additional characterizations and several extensions.

It was established in \cite{Kru03,Kru04} that, similar to the classical Lagrange multiplier rule, the conclusion of the extremal principle (the generalized separation), being a dual necessary condition of extremality, actually characterizes a property which is weaker than local extremality.
This property can be interpreted as a kind of stationarity of a collection of sets.
The explicit primal space definition of this property was introduced (cf. Definition~\ref{D2.5}(iv)), and it was shown (by refining slightly the original proof of the extremal principle) that for this property the conclusion of the extremal principle is not only necessary (in the Asplund space setting), but also sufficient, thus producing the ultimate conventional version of the extremal principle, known as the \emph{extended extremal principle}; cf. Theorem~\ref{T2.9}.
This stationarity property called \emph{approximate stationarity} \cite{Kru09} happens to be strongly connected to (in fact the negation of) another important geometric property of a collection of sets called \emph{transversality} (also known under other names; cf. Definition~\ref{D2.15} and Theorem~\ref{T2.13}), having its roots in the classical differential geometry and instrumental for the convergence analysis of alternating projections; see a discussion of the role of the latter property in \cite{KruLukTha18}, while a table illustrating the evolution of the terminology can be found in \cite[Section~2]{KruLukTha17}.
The connections of the approximate stationarity and (extended) extremal principle can be traced further to the fundamental for variational analysis \emph{metric regularity} property of \SVM s and the corresponding coderivative criterion; see \cite{Kru09}.

The powerful (extended) extremal principle, despite its recognized universality and wide applicability, has its limitations.
The second author encountered a problem while attempting to extend the extremal principle to infinite collections of sets in \cite{KruLop12.1}.
The initial idea to consider families of finite subcollections and apply the conventional extremal principle to each of them failed, as well as the belief of the authors in the unlimited universality of the conventional extremal principle.
Uniform estimates were required, holding for all finite subcollections, and the conventional extremal principle was unable to provide such estimates.
The solution was found within the conventional (extended) extremal principle, more precisely, in its proof, providing another piece of evidence for the well-known fact that important mathematical results are often deeper than their statements.

The conventional extended extremal principle asserts the equivalence of two properties of a collection of sets: one in the primal space (approximate stationarity) and another one in the dual space (generalized separability); see Theorem~\ref{T2.9}, which for completeness gives two (equivalent) versions of the latter property.
Both (in the case of Theorem~\ref{T2.9} all three) properties are formulated `for any $\eps$ there exist \ldots', and as such, the $\eps$'s and the other parameters in different properties seem completely independent.
However, this cannot be true.
To prove that one property of this kind implies another one, there must be a way, given any $\eps$ in the second property, to construct another one to be plugged into the first property to produce the needed estimates, i.e. the parameters must be related.
The relationship between the parameters in the two parts of the conventional (extended) extremal principle is exactly what was required to establish the uniform estimates needed in \cite{KruLop12.1}, and it indeed could be found in the proof of the conventional result.
This observation made the authors of \cite{KruLop12.1} carve out the core part of (the proof of) the conventional (extended) extremal principle with $\eps$ and other parameters fixed (a kind of $\eps$-extremality), and formulate it as a separate statement in \cite[Theorem~3.1]{KruLop12.1} (see Proposition~\ref{P5.16}(iii)), exposing the relationship between the parameters hidden in the proof of the conventional statement.
As a result, both the conventional statement and its extension to infinite collections of sets are corollaries of \cite[Theorem~3.1]{KruLop12.1}.

There are many equivalent formulations of the primal space approximate stationarity property (as well as other extremality and stationarity properties) and the dual space generalized separability properties, all formulated in the form `for any $\eps$ there exist \ldots'; see Sect.~\ref{S3}.
The relationships between the parameters involved in these formulations can also be established, and are important for identifying and analysing the core arguments in the conventional proofs of metric and dual characterizations of the extremality and stationarity and their extensions.
This analysis is performed in Sects.~\ref{S5} and \ref{S6}.

Another successful `surgical operation' on the proof of the conventional extremal principle was done earlier by Zheng and Ng in \cite{ZheNg05}, producing an $\eps$-separability characterization of another kind of $\eps$-extremality property with the explicit relationship between the $\eps$'s in both properties.
Unlike the conventional extremal principle and its extension in \cite[Theorem~3.1]{KruLop12.1} which assume the sets to have a common point, \cite[Lemmas~2.2 and 2.2']{ZheNg05} assume, on the contrary, that the intersection of the sets is empty.
At the same time, their proofs follow the original ideas from \cite{KruMor80,KruMor80.2}, utilizing the \EVP\ and either the Asplund space fuzzy or the general Banach space Clarke--Rockafellar subdifferential sum rule; cf. Lemma~\ref{SR}.
The statements have been further polished and analysed in a sequence of subsequent papers \cite{ZheNg06,LiNgZhe07,LiTanYuWei08,ZheNg11}.
In particular, it was proved by Guoyin Li et al in \cite[Theorem~3.1]{LiTanYuWei08} that the conclusion of \cite[Lemma 2.2]{ZheNg05} is actually equivalent to the \EVP, and as such implies (the Banach space with Clarke normal cones version of) the extremal principle.
This fact confirms the need to move from the conventional `for any $\eps$ \ldots' extremality statements to more subtle ones with $\eps$ fixed.
The most advanced version of the Zheng and Ng lemma was given in \cite[Theorems~3.1 and 3.4]{ZheNg11} (\emph{unified separation theorems}; see Theorem~\ref{ZhNg}), where, besides other improvements, an additional condition was added to the concluding part, relating the dual vectors involved in the $\eps$-se\-parability characterization with certain primal space vectors involved in the original $\eps$-extremality property.

A partial comparison of the assumptions and conclusions in the two existing $\eps$-extremality statements in \cite[Theorem~3.1]{KruLop12.1} and \cite[Theorems~3.1 and 3.4]{ZheNg11} has been done recently in \cite{BuiKru18}, where it was also noted that they are formulated using in a sense different languages and are in general incomparable.
In the current paper, we formulate in Theorem~\ref{T5.1} (and prove) a new general $\eps$-extremality statement which exposes the core arguments and the role of the parameters in the conventional proofs of dual characterizations of the extremality and stationarity and, in particular, implies \cite[Theorem~3.1]{KruLop12.1} and \cite[Theorems~3.1 and 3.4]{ZheNg11}.
In view of the above-mentioned result by Guoyin Li et al, the conclusion of Theorem~\ref{T5.1} is also equivalent to the \EVP.
We also establish in Sect.~\ref{S6} a series of other consequences of Theorem~\ref{T5.1}, covering primal and dual space conditions involved in (hopefully) all known formulations of extremality/stationarity and generalized separability properties.

The structure of the paper is as follows.
Section~\ref{S2} contains some preliminary definitions and facts used throughout the paper.
Some of the facts are new.
In Sect.~\ref{S2.1} dedicated to the \EVP, two new versions of this classical result are formulated: the \emph{Geometric Ekeland Variational Principle (GEVP)} and the \emph{Asymmetric Geometric Ekeland Variational Principle (AGEVP)}, and the equivalence of these two statements to the conventional \EVP\ is proved in Proposition~\ref{P2.1}.
In Sect.~\ref{S2.2}, we recall the definitions of the \Fr\ and Clarke normal cones and subdifferentials, several versions of subdifferential sum rules and a few other calculus facts needed in the sequel.
Section~\ref{S2.3} is dedicated to the extremality and stationarity properties, and the (extended) extremal principle, and collects in one place all the relevant definitions, primal space metric and dual space normal characterizations as well as some discussions and historical comments.
In Sect.~\ref{S2.4}, we prove several assertions containing elementary arguments which are used in proving the equivalence of various generalized separation statements and, more importantly, provide quantitative estimates for parameters involved in such statements.
Such arguments are usually hidden within proofs.
We make them explicit here for use in the subsequent sections and also outside this paper.

Section~\ref{S3} recalls several known modifications of the extremality and stationarity properties and their dual characterizations, which are going to be important for our study of the $\eps$-versions of these properties in the subsequent sections, and introduces the new \emph{approximate $\al$-stationarity} property.
In Sect.~\ref{S4}, we discuss several ways of defining `distances' between $n$ sets, quantifying the closeness of a finite number of sets, i.e. how `far apart' they are, or, at least, whether they have a common point.
Such distances play an important role when studying extremality, stationarity and regularity properties of collections of sets.
An example of a distance of this kind was considered in \cite{ZheNg11} under the name \emph{nonintersect index}.
Primal space characterizations of `$\eps$-closest' points of a finite collection of sets are established in the form common for conventional formulations of approximate stationarity, thus, building a bridge between the two seemingly different languages used in defining `extremal' properties of collections of sets.

Sections~\ref{S5} and \ref{S6} are dedicated to the comparison of, respectively, primal and dual space `fixed $\eps$' conditions involved in definitions and characterizations of extremality, stationarity and generalized separation properties.
In Sect.~\ref{S5}, we extend the Asymmetric Geometric Ekeland Variational Principle (AGEVP) from Sect.~\ref{S2} to the case of $n\ge2$ sets and establish a metric counterpart of the `fixed $\eps$' conditions involved in the definition of local extremality, thus producing the core component of the equivalent metric characterizations of the approximate $\al$-stationarity, and consequently of the approximate stationarity and transversality.
Section~\ref{S6} presents a series of `generalized separation' statements, providing dual characterizations of typical `extremal' arrangements of collections of sets, discussed in the preceding sections, and examines the relationships between them.
These statements refine core arguments, which can be found in various existing versions of the (extended) extremal principle, as well as some new extensions.
Their conclusions combine assertions in terms of the \Fr\ normal cones in Asplund spaces and in terms of the Clarke normal cones in general Banach spaces.
All the separation statements in this section as well as the conventional extremal principle and its extensions in \cite{KruLop12.1} and \cite{ZheNg11} are consequences of the general Theorem~\ref{T5.1}.
Its proof uses the extension of the AGEVP from Sect.~\ref{S5} as a replacement for the conventional \EVP\ in generalized separation statements.
The dual characterizations of the `extremal' arrangements of collections of sets given in the statements in the first part of Sect.~\ref{S6} are partially reversed in several statements at the end of the section in the setting of a general normed vector space.
This leads, in particular, to a full dual characterization of the approximate $\al$-stationa\-rity, and as a consequence of the approximate stationa\-rity and transversality.

\section{Preliminaries}\label{S2}

Our basic notation is standard; see, e.g., \cite{RocWet98,Mor06.1,Pen13,DonRoc14,Pen16,Iof17}.
Throughout the paper, $X$ is either a metric or (more often) a normed vector space.
In the latter case, we often require it to be Banach or Asplund.
The distance and the norm are denoted by $d(\cdot,\cdot)$ and $\|\cdot\|$, respectively.
We use the same symbols to denote distances and norms in all spaces (primal and dual).
When considering products of spaces, we usually assume them equipped with the maximum distance or norm.
$B_\de(x)$ and $\overline{B}_\de(x)$ denote, respectively, the open and closed balls with centre $x$ and radius $\de>0$.
Given a point $x$ and a set $A$ in $X$,
$d(x,A):=\inf_{a\in A}d(x,a)$ denotes the distance from $x$ to $A$; in particular $d(x,\es):=+\infty$.
Given two subsets $A,B\subset X$,
$d(A,B):=\inf_{a\in A}d(a,B)$ denotes the distance between $A$ and $B$.
Given a set $A$, a point $a\in A$ and a number $\de>0$, we call the set $A\cap B_\de(a)$ a \emph{localization} of the set $A$ near $a$.
If $X$ is a normed vector space,
its topological dual is denoted by $X^*$, while $\langle\cdot,\cdot\rangle$ denotes the bilinear form defining the pairing between the two spaces.
The open unit balls in $X$ and $X^*$ are denoted by $\B$ and $\B^*$, respectively.
$\mathbb{N}$ stands for the set of all positive integers.
We also use the notation $\R_\infty:=\R\cup\{+\infty\}$.
Given a function $f:X\to\R_\infty$, its domain is the set $\dom f:=\{x\in X:\, f(x)<\infty\}$.

{\bf \subsection{Ekeland Variational Principle}\label{S2.1}}

Not surprisingly, the following classical result due to Ekeland \cite{Eke74} (see also \cite{Mor06.1,DonRoc14,Pen13,Pen16,Iof17}) plays the key role in the subsequent studies.

\spnewtheorem*{evp}{Ekeland Variational Principle (EVP)} {\bf}{\it}
\begin{evp}
Suppose $X$ is a complete metric space, $f: X\to\R_\infty$ is lower semicontinuous, $\bx\in X$, and $\varepsilon>0.$
If
$$
f(\bar{x})<\inf_{X}f+\varepsilon,
$$
then, for any $\la>0$, there exists an $\hat{x}\in X$ such that
\begin{enumerate}
\item $d(\hat{x},\bx)<\lambda$;
\item $f(\hat{x})\le f(\bx)$;
\item
$f(x)+(\varepsilon/\lambda)d(x,\hat{x})>f(\hat{x})$ for all $x\in X\setminus\{\bx\}.$
\end{enumerate}
\end{evp}

In the current paper dedicated to collections of sets, we are going to use the next two geometric versions of the \EVP.
They characterise the mutual arrangement of a pair of sets in a complete metric space with respect to a pair of points, being almost (up to $\eps)$ closest points of these sets and, similarly to the conventional EVP, establish the existence of another pair of points arbitrarily close (up to an additional parameter or a pair of parameters) to the given one and minimizing a certain perturbed function.
The perturbed functions in both assertions involve the distance between localizations of the sets near these points, i.e. intersections of the sets with neighbourhoods of the points.

\spnewtheorem*{gevp}{Geometric Ekeland Variational Principle (GEVP)} {\bf}{\it}
\spnewtheorem*{agevp}{Asymmetric Geometric Ekeland Variational Principle (AGEVP)} {\bf}{\it}
\begin{gevp}
Suppose $X$ is a complete metric space, $A$ and $B$ are closed subsets of $X$, $a\in A$, $b\in B$, and $\eps>0$.
If
\sloppy
\begin{equation}\label{gevp1}
d(a,b)<d(A,B)+\eps,
\end{equation}
then, for any $\la>0$, there exist $\hat a\in A\cap B_\la(a)$ and $\hat b\in B\cap B_\la (b)$ such that
\begin{enumerate}
\item
$d(\hat a,\hat b)\le d(a,b)$;
\item
$d(A\cap B_{\xi}(\hat a),B\cap B_{\xi}(\hat b)) +\frac{\xi\eps}{\la}>d(\hat a,\hat b)$ for all $\xi>0$.
\end{enumerate}
\end{gevp}
\begin{agevp}
Suppose $X$ is a complete metric space, $A$ and $B$ are closed subsets of $X$, $a\in A$, $b\in B$, and $\eps>0$.
If condition \eqref{gevp1} is satisfied,
then, for any $\la,\rho>0$, there exist $\hat a\in A\cap B_{\la}(a)$ and $\hat b\in B\cap B_{\rho}(b)$ such that
\begin{enumerate}
\item
$d(\hat a,\hat b)\le d(a,b)$;
\item
$d(A\cap B_{\xi\la}(\hat a),B\cap  B_{\xi\rho}(\hat b)) +{\xi\eps}>d(\hat a,\hat b)$ for all $\xi>0$.
\end{enumerate}
\end{agevp}

Unlike GEVP, its asymmetric version AGEVP allows for balls with different radii to be used in the localizations of the sets.
This feature is going to play an important role in our subsequent analysis.
It is easy to see that GEVP is a particular case of AGEVP with $\rho=\la$.
We now show that the two geometric versions of the \EVP\ formulated above are both equivalent to the conventional one.

\begin{proposition}\label{P2.1}
EVP \iff GEVP \iff AGEVP.
\end{proposition}
\begin{Proof}
We first show that EVP \folgt AGEVP.
Let the assumptions of AGEVP be satisfied.
Given numbers $\la,\rho>0$, we consider the space $X\times X$ with a metric defined as follows:
\begin{gather}\label{P2.1P-1}
d_{\la,\rho}((x,y),(u,v)) :=\max\left\{\frac{1}{\la}d(x,u), \frac{1}{\rho}d(y,v)\right\}\quad (x,y,u,v\in X).
\end{gather}
Since $X$ is complete, $(A\times B,d_{\la,\rho})$ is a complete space.
Choose an $\eps'\in]0,\eps[$ such that \eqref{gevp1} is satisfied with $\eps'$ in place of $\eps$.
EVP applied to the function $d$ on $(A\times B,d_{\la,\rho})$ gives the existence of points $\hat a\in A$ and $\hat b\in B$ such that
\begin{gather}\label{P2.1P-2}
d_{\la,\rho}((\hat a,\hat b),(a,b))<1,\quad
d(\hat a,\hat b)\le d(a,b),
\\\label{P2.1P-3}
d(x,y)-d(\hat a,\hat b)+\eps' d_{\la,\rho}((x,y),(\hat a,\hat b))\ge0
\qdtx{for all}
(x,y)\in A\times B.
\end{gather}
In view of the definition \eqref{P2.1P-1}, the first inequality in \eqref{P2.1P-2} is equivalent to the following two: $d(\hat a,a)<\la$ and $d(\hat b,b)<\rho$.
Given any $\xi>0$, $x\in A\cap B_{\xi\la}(\hat a)$ and $y\in B\cap B_{\xi\rho}(\hat b)$, by \eqref{P2.1P-1} and \eqref{P2.1P-3}, we have, respectively, $d_{\la,\rho}((x,y),(\hat a,\hat b))<\xi$ and
$$
d(x,y)\ge d(\hat a,\hat b)-\eps' d_{\la,\rho}((x,y),(\hat a,\hat b))>d(\hat a,\hat b)-{\xi\eps'}.
$$
It follows that $d(A\cap B_{\xi\la}(\hat a), B\cap B_{\xi\rho}(\hat b))\ge d(\hat a,\hat b)-\xi\eps'> d(\hat a,\hat b)-\xi\eps$.
This proves AGEVP.

The implication AGEVP \folgt GEVP is straightforward.

To complete the proof, we next show that GEVP \folgt EVP.
Let the assumptions of EVP be satisfied.
Choose a positive number $\eps'<\eps$ such that $f(\bar{x})<\inf_{X}f+\varepsilon'$ and another number $\al$ such that
$0<\alpha<\la(\eps'{}\iv-\eps\iv)$.
We are going to consider the space $X\times\R$ with the metric $d:=d_X+\alpha|\cdot|$, which makes $X\times\R$ a complete metric space, two closed subsets of $X\times \R$: $A:=\{(x,y):\, x\in X,\; y\ge f(x)\}$ and $B:=X\times\{M\}$, where $M:=\inf_X f$, and two points $a:=(\bx,f(\bx))\in A$ and $b:=(\bx,M)\in B$.
We have $d(A,B)=0$ and
$$
d(a,b)=\alpha(f(\bx)-M)<\alpha\eps'=d(A,B)+\alpha\eps'.
$$
GEVP gives the existence of points $\hat a=(\hat x,\hat y)\in A\cap B_{\la}(a)$ and $\hat b=(\hat x',M)\in B\cap B_\la (b)$ satisfying
\begin{gather}\label{P2.1P-4}
d(\hat a,\hat b)\le d(a,b),
\\\label{P2.1P-5}
d(A\cap B_{\xi}(\hat a),B\cap B_{\xi}(\hat b)) +\frac{\xi \alpha\eps'}{\la}>d(\hat a,\hat b)
\qdtx{for all}
\xi>0.
\end{gather}
Observe that $d(\hat x,\bx)<\la$, and consequently, $(\hat x,M)\in B\cap B_\la(b)$.
Moreover, $\hat x'=\hat x$.
Indeed, if $\hat x\neq\hat x'$, then we can take $\xi:=d(\hat x,\hat x')$.
Then, $(\hat x,\hat y)\in A\cap B_{\xi}(\hat a)$, $(\hat x, M)\in B\cap\overline{B}_{\xi}(\hat b)$ and
\begin{align*}
d((\hat x,\hat y),(\hat x, M))+
\frac{\xi \alpha\eps'}{\la}
&<\alpha(\hat y-M)+\left(1-\frac{\eps'}{\eps}\right) d(\hat x,\hat x')
\\
&<\alpha(\hat y-M)+d(\hat x,\hat x')=d(\hat a,\hat b),
\end{align*}
which contradicts \eqref{P2.1P-5}; hence $\hat x'=\hat x$ and $\hat b=(\hat x,M)$.
Consequently, condition \eqref{P2.1P-4} reduces to $\alpha(\hat y-M)\le \alpha(f(\bx)-M)$, which implies that $f(\hat x)\le\hat y\le f(\bx)$.

It remains to prove condition (iii) in EVP.
Let $x\ne\hat x$.
If $f(x)\ge f(\hat x)$, the condition holds trivially.
Let $f(x)<f(\hat x)$, and take $\xi:={d(x,\hat x) +\alpha(\hat y-f(x))}$.
Then $(x,f(x))\in A\cap \overline{B}_{\xi}(\hat a)$, $(x,M)\in B\cap B_{\xi}(\hat b)$, and by \eqref{P2.1P-5},
\begin{align*}
\alpha(f(x)-M)+
\frac{\alpha\eps'}{\la}
(d(x,\hat x)+\alpha(\hat y-f(x)))
&=d((x,f(x)),(x,M))+
\frac{\xi \alpha\eps'}{\la}
\\
&\ge d(\hat a,\hat b)=\alpha(\hat y-M).
\end{align*}
Thus,
$$
f(x)+\frac{\eps'}{\la}(d(x,\hat x)+\alpha(\hat y-f(x)))
\ge\hat y,
$$
or equivalently,
$$
f(x)+\frac{\eps'}{\la-\alpha\eps'}d(x,\hat x)\ge\hat y.
$$
By the definition of $\al$, we have
$\la-\alpha\eps'>\la\eps'\eps\iv$.
Hence,
$$
f(x)+\frac{\eps}{\la}d(x,\hat x)> f(x)+\frac{\eps'}{\la-\alpha\eps'}d(x,\hat x)\ge\hat y\ge f(\hat x).
$$
The proof is complete.
\qed\end{Proof}

Condition (ii) in GEVP and AGEVP corresponds to condition (iii) in the conventional EVP, while condition (i) corresponds to the pair of conditions (i) and (ii).
The $\eps$-closeness condition \eqref{gevp1} is going to play an important role in our analysis.
It will be further discussed in Sects.~\ref{S4} and \ref{S5}.

{\bf \subsection{Normal Cones and Subdifferentials}\label{S2.2}}

In Sect.~\ref{S6}, we use
dual tools -- normal cones and subdifferentials, usually in the Fr\'echet or Clarke sense.
Given a subset $A$ of a normed vector space $X$ and a point $\bx\in A$, the set (cf. \cite{Kru03})
\begin{gather}\label{NC}
N_{A}^F(\bx):= \left\{x^\ast\in X^\ast:\,
\limsup_{x\to\bx,\,x\in A\setminus\{\bx\}} \frac {\langle x^\ast,x-\bx\rangle}
{\|x-\bx\|} \le 0 \right\}
\end{gather}
is the \emph{Fr\'echet normal cone} to $A$ at $\bx$.
It is a nonempty
closed convex cone, often trivial
(i.e. $N_{A}^F(\bx)=\{0\}$).
The \emph{Clarke normal cone} to $A$ at $\bx$ is defined as the set (cf. \cite{Cla83})
\begin{gather}\label{NCC}
N_{A}^C(\bx):= \left\{x^\ast\in X^\ast:\,
\ang{x^\ast,z}\le0
\qdtx{for all}
z\in T_{A}^C(\bx)\right\},
\end{gather}
where $T_{A}^C(\bx)$ is the \emph{Clarke tangent cone} to $A$ at $\bx$:
\begin{align}\notag
T_{A}^C(\bx):= \Big\{z\in X:\,
&\forall x_k\stackrel{A}{\rightarrow}\bx,\;\forall t_k\downarrow0,\;\exists z_k\to z
\qdtx{such that} 
\\\label{TCC}
&x_k+t_kz_k\in A \qdtx{for all}
k\in\N\Big\}.
\end{align}
The set \eqref{NCC} is a nonempty
weak$^*$
closed convex cone, and $N_{A}^F(\bx)\subset N_{A}^C(\bx)$.
If $A$ is a convex set, then \eqref{NC} and \eqref{NCC} reduce to the normal cone in the sense of convex analysis (cf., e.g., \cite[Proposition~1.19]{Kru03}, \cite[Proposition~2.4.4]{Cla83}):
\begin{gather*}\label{CNC}
N_{A}(\bx):= \left\{x^*\in X^*:\, \langle x^*,x-\bx \rangle \leq 0 \qdtx{for all} x\in A\right\}.
\end{gather*}
We will often use the generic notation $N$ for both Fr\'echet and Clarke normal cones, specifying wherever necessary that either $N:=N^F$ or $N:=N^C$.

Given a function $f:X\to\R_\infty$ and a point $\bx\in\dom f$, the \emph{Fr\'echet} and \emph{Clarke subdifferentials} of $f$ at $\bx$ can be defined via the respective normal cones to the epigraph $\epi f:=\{(x,\al)\in X\times\R:\, f(x)\le\al\}$ of $f$ as follows:
\begin{gather*}
\partial^F{f}(\bx):= \left\{x^\ast\in X^\ast:\,
(x^*,-1)\in N_{\epi f}^F(\bx,f(\bx))\right\},
\\
\partial^C{f}(\bx):= \left\{x^\ast\in X^\ast:\,
(x^*,-1)\in N_{\epi f}^C(\bx,f(\bx))\right\}.
\end{gather*}
These are closed convex sets, the first one often empty, and $\partial^F{f}(\bx)\subset\partial^C{f}(\bx)$.
If $f$ is convex, they both reduce to the subdifferential in the sense of convex analysis (cf., e.g., \cite[Proposition~1.2]{Kru03}, \cite[Proposition~2.2.7]{Cla83}):
\begin{gather*}
\partial{f}(\bx):= \left\{x^\ast\in X^\ast:\,
f(x)-f(\bx)-\langle{x}^\ast,x-\bx\rangle\ge 0 \qdtx{for all} x\in X \right\}.
\end{gather*}
By convention, we set $N_{A}^F(\bx)=N_{A}^C(\bx):=\es$ if $\bx\notin A$ and $\partial^F{f}(\bx)=\partial^C{f}(\bx):=\es$ if $\bx\notin\dom f$.
It is easy to check that $N_{A}^F(\bx)=\partial^Fi_A(\bx)$ and $N_{A}^C(\bx)=\partial^Ci_A(\bx)$, where $i_A$ is the \emph{indicator function} of $A$: $i_A(x)=0$ if $x\in A$ and $i_A(x)=\infty$ if $x\notin A$; cf., e.g., \cite[Proposition~1.18]{Kru03}, \cite[Proposition~2.4.12]{Cla83}.

Several kinds of \emph{subdifferential sum rules} are used in Sect.~\ref{S6} when deducing dual space results.
They are collected in the next lemma.

\begin{lemma}[Subdifferential sum rules] \label{SR}
Suppose $X$ is a normed vector space, $f_1,f_2:X\to\R_\infty$, and $\bx\in\dom f_1\cap\dom f_2$.
\begin{enumerate}
\item
{\bf Convex sum rule}. Suppose
$f_1$ and $f_2$ are convex and $f_1$ is continuous at a point in $\dom f_2$.
Then,
$$
\partial (f_1+f_2) (\bar x) = \sd f_1(\bx) +\partial f_2(\bx).
$$

\item
{\bf Differentiable sum rule}. Suppose
$f_1$ is Fr\'echet differentiable at $\bx$.
Then,
$$
\partial^F(f_1+f_2) (\bar x) = \nabla f_1(\bx) +\partial^Ff_2(\bx).
$$

\item
{\bf Fuzzy sum rule}. Suppose $X$ is Asplund,
$f_1$ is Lipschitz continuous, and
$f_2$
is lower semicontinuous in a neighbourhood of $\bar x$.
Then, for any $\varepsilon>0$, there exist $x_1,x_2\in X$ with $\|x_i-\bar x\|<\varepsilon$, $|f_i(x_i)-f_i(\bar x)|<\varepsilon$ $(i=1,2)$, such that
$$
\partial^F(f_1+f_2) (\bar x) \subset \partial^Ff_1(x_1) +\partial^Ff_2(x_2) + \varepsilon\B^\ast.
$$

\item
{\bf Clarke--Rockafellar sum rule}. Suppose
$f_1$ is Lipschitz continuous and
$f_2$
is lower semicontinuous in a neighbourhood of $\bar x$.
Then,
$$
\partial^C(f_1+f_2)(\bar x)\subset\sd^C f_1(\bx) +\partial^Cf_2(\bx).
$$
\end{enumerate}
\end{lemma}

The first sum rule in the lemma above is the conventional subdifferential sum rule of convex analysis; see, e.g., \cite[Theorem~0.3.3]{IofTik79} and
\cite[Theorem~2.8.7]{Zal02}.
Together with the second one, theses are examples of \emph{exact} sum rules.
The third sum rule is known as the \emph{fuzzy} or \emph{approximate} sum rule (Fabian \cite{Fab89}) for Fr\'echet subdifferentials in Asplund spaces; cf., e.g., \cite[Rule~2.2]{Kru03} and \cite[Theorem~2.33]{Mor06.1}.
Note that, unlike the sum rules in parts (i) and (ii) of the lemma, the subdifferentials in the \RHS\ of the inclusion are computed not at the reference point, but at some points nearby.
This explains the name.
The fourth sum rule is formulated in terms of Clarke subdifferentials.
It was established in Rockafellar \cite[Theorem~2]{Roc79}.
Similar to the previous one, it is valid generally only as inclusion.
Nevertheless, it is another example of exact sum rule.

Recall that a Banach space is \emph{Asplund} if every continuous convex function on an open convex set is Fr\'echet differentiable on a dense subset \cite{Phe93}, or equivalently, if the dual of each its separable subspace is separable.
We refer the reader to \cite{Phe93,Mor06.1,BorZhu05} for discussions about and characterizations of Asplund spaces.
All reflexive, particularly, all finite dimensional Banach spaces are Asplund.

The following facts are immediate consequences of the definition of the \Fr\ subdifferential and normal cone (cf., e.g., \cite[Propositions~1.10 and 1.29]{Kru03}).

\begin{lemma}\label{L2.3}
Suppose $X$ is a normed vector space and $f:X\to\R_\infty$.
If $\bx\in\dom f$ is a point of local minimum of $f$, then $0\in\sd^Ff(\bx)$.
\end{lemma}

\begin{lemma}\label{L2.4}
Suppose $X_1$ and $X_2$ are normed vector spaces, $\bx_1\in A_1\subset X_1$, and $\bx_2\in A_2\subset X_2$.
Then,
$$
N_{A_1\times A_2}^F(\bx_1,\bx_2)=N_{A_1}^F(\bx_1)\times N_{A_2}^F(\bx_2).
$$
\end{lemma}

{\bf \subsection{Extremality, Stationarity and (Extended) Extremal Principle}\label{S2.3}}

In this subsection, we recall and discuss the conventional definitions of extremality, local extremality, stationarity and approximate stationarity of collections of sets, and the (extended) extremal principle.
Here and in the rest of the paper, we consider $n$ sets $\Omega_1,\ldots,\Omega_n$ ($2\le n<\infty$) and write $\{\Omega_1,\ldots,\Omega_n\}$ to denote the collection of the sets as a single object.

The next definition collects several extremality and stationarity properties of collections of sets.

\begin{definition}\label{D2.5}
Suppose $\Omega_1,\ldots,\Omega_n$ are subsets of a normed vector space $X$ and $\bx\in\cap_{i=1}^n\Omega_i$.
The collection $\{\Omega_1,\ldots,\Omega_n\}$ is \begin{enumerate}
\item
\textbf{extremal} iff
for any $\eps>0$, there exist vectors
$a_i\in{X}$ $(i=1,\ldots,n)$ satisfying
\begin{gather}\label{P1}\tag{P1}
\bigcap_{i=1}^n(\Omega_i-a_i)=\emptyset\AND
\max_{1\le i\le n}\norm{a_i}<\eps;
\end{gather}
\item
\textbf{locally extremal} at $\bx$ iff there exists a number $\rho\in]0,\infty]$ such that,
for any $\eps>0$, there are vectors
$a_i\in{X}$ $(i=1,\ldots,n)$ satisfying
\begin{gather}\label{P2}\tag{P2}
\bigcap_{i=1}^n(\Omega_i-a_i)\cap{B}_\rho(\bar{x})
=\emptyset\AND
\max_{1\le i\le n}\norm{a_i}<\eps;
\end{gather}
\item
\textbf{stationary} at $\bx$ iff
for any $\eps>0$,
there exist a number $\rho\in]0,\eps[$
and vectors
$a_i\in{X}$ $(i=1,\ldots,n)$ satisfying
\begin{gather*}
\bigcap_{i=1}^n(\Omega_i-a_i)\cap{B}_\rho(\bar{x})
=\emptyset\AND
\max_{1\le i\le n}\norm{a_i}<\eps\rho;
\end{gather*}
\item
\textbf{approximately stationary} at $\bx$ iff
for any $\eps>0$,
there exist a number $\rho\in]0,\eps[$, points $\omega_i\in\Omega_i\cap B_\eps(\bx)$
and vectors
$a_i\in{X}$ $(i=1,\ldots,n)$ satisfying
\begin{gather}\label{P3}\tag{P3}
\bigcap_{i=1}^n(\Omega_i-\omega_i-a_i)\cap(\rho\B)
=\emptyset\AND
\max_{1\le i\le n}\norm{a_i}<\eps\rho.
\end{gather}
\end{enumerate}
\end{definition}

The formulas in Definition~\ref{D2.5} and several their modifications discussed in Sect.~\ref{S3} are central for our analysis and are going to be extensively referred to throughout the paper.
We will use special P tags: (P1), (P2), \ldots for such formulas as well as some other formulas in Sect.~\ref{S3} involved in primal space metric characterizations of extremality and stationarity.
D tags: (D1), (D2), \ldots are reserved for the formulas involved in the corresponding dual space characterizations of extremality and stationarity.

Condition (i) (condition (ii)) in Definition~\ref{D2.5} means that an appropriate arbitrarily small shift of the sets makes them nonintersecting (in a neighbourhood of $\bx$).
This is a very general model embracing many optimality notions.
The other two conditions correspond to more subtle properties of optimization problems, closer to stationarity.

The properties in parts (i) and (ii) of Definition~\ref{D2.5} were introduced in \cite{KruMor80} and \cite{Kru81.2}, respectively; see also \cite{Kru03,Mor06.1}.
The properties in parts (iii) and (iv) first appeared in \cite{Kru04} and \cite{Kru98}, respectively; see also \cite{Kru06,BuiKru18}.
Property (iv) was referred to in \cite{Kru98} as \emph{extremality near $\bx$}.
The name \emph{approximate stationarity} was suggested in \cite{Kru09}.

Unlike condition (ii), in conditions (iii) and (iv) the magnitudes of the `shifts' of the sets are related to that of the neighbourhood in which the sets become nonintersecting, namely $\max_{1\le i\le n} \norm{a_i}/\rho<\eps$.
Compared to (iii), in condition (iv), instead of the common point $\bx$, each set $\Omega_i$ is considered near its own point $\omega_i$.

The relationships between the properties in Definition~\ref{D2.5} are straightforward.
The equivalences in part (ii) of the proposition below were proved in \cite[Proposition~14]{Kru05}.

\begin{proposition}\label{P2.6}
Suppose $\Omega_1,\ldots,\Omega_n$ are subsets of a normed vector space $X$ and $\bx\in\cap_{i=1}^n\Omega_i$.
\begin{enumerate}
\item
For the properties in Definition~\ref{D2.5}, the following implications hold true:
{\rm (i) \folgt (ii) \folgt (iii) \folgt (iv)}.
\item
If the sets are convex, then the implications in the previous item hold as equivalences:
{\rm (i) \iff (ii) \iff (iii) \iff (iv)}.
\item
If the collection $\{\Omega_1,\ldots,\Omega_n\}$ is locally extremal at $\bx$ with $\rho=\infty$, then it is extremal.
\item
If the collection $\{\Omega_1,\ldots,\Omega_n\}$ is locally extremal at $\bx$ with some $\rho\in]0,\infty]$, then the collection of $n+1$ sets $\Omega_1,\ldots,\Omega_n$, ${B}_\rho(\bar{x})$ is extremal.
\end{enumerate}
\end{proposition}

It is easy to check that all the implications in Proposition~\ref{P2.6}(i) can be strict;
see examples in \cite{Kru05,Kru09,BuiKru18}.

Thanks to Proposition~\ref{P2.6}(i), approximate stationarity is the weakest of the four properties in Definition~\ref{D2.5}.
It happens to be an important type of mutual arrangement of a collection of sets in space.
The next two statements provide several, respectively, primal space metric and dual space (in terms of Fr\'echet normals) necessary and sufficient criteria for this property.
The characterizations formulated in these statements are obviously necessary for each of the other three properties in Definition~\ref{D2.5}.

\begin{theorem}[Approximate stationarity: metric criteria]\label{T2.8}
Suppose $\Omega_1,\ldots,\Omega_n$ are subsets of a normed vector space $X$ and $\bx\in\cap_{i=1}^n\Omega_i$.
The following conditions are equivalent:
\begin{enumerate}
\item
the collection $\{\Omega_1,\ldots,\Omega_n\}$ is approximately stationary at $\bx$;
\item
for any $\eps>0$, there exist vectors $a_i\in X$ $(i=1,\ldots,n)$ such that
\begin{equation*}
\eps d\left(\bx,\bigcap_{i=1}^n(\Omega_i-a_i)\right)> \max_{1\le i\le n} d\left(\bx,\Omega_i-a_i\right)
\qdtx{and}
\max_{1\le i\le n}\norm{a_i}<\eps;
\end{equation*}
\item
for any $\eps>0$, there exist a point $x\in B_\eps(\bx)$ and vectors $a_i\in X$ $(i=1,\ldots,n)$ such that
\begin{equation*}
\eps d\left(x,\bigcap_{i=1}^n(\Omega_i-a_i)\right) >\max_{1\le i\le n}d(x,\Omega_i-a_i)
\qdtx{and}
\max_{1\le i\le n}\norm{a_i}<\eps.
\end{equation*}
\cnta
\end{enumerate}
\end{theorem}

The equivalence of the conditions (i) and (iii) in Theorem~\ref{T2.8} is a consequence of \cite[Theorem~1]{Kru05}; cf. \cite[Theorem~1]{Kru06}.
The equivalence of the all three conditions follows from a more general statement in
Corollary~\ref{C4.9} in Sect.~\ref{S5}.
Each of the conditions (ii) and (iii) can be used as an equivalent definition of approximate stationarity.

The next theorem
can be considered as a generalization of the classical convex \emph{separation theorem} to collections of nonconvex sets.

\begin{theorem}[Extended extremal principle]\label{T2.9}
Suppose $\Omega_1,\ldots,\Omega_n$ are closed subsets of an Asplund space $X$ and
${\bx\in\cap_{i=1}^n\Omega_i}$.
With $N$ standing for the \Fr\ normal cone ($N:=N^F$), the following conditions are equivalent:
\sloppy
\begin{enumerate}
\item
the collection $\{\Omega_1,\ldots,\Omega_n\}$ is approximately stationary at $\bx$;
\item
for any $\eps>0$, there exist points
$\omega_i\in\Omega_i\cap{B}_\eps(\bx)$ and vectors $x_i^*\in X^*$ $(i=1,\ldots,n)$ such that
\begin{gather}\label{D1}\tag{D1}
\norm{\sum_{i=1}^n x^*_i}<\eps,\quad
x^*_i\in N_{\Omega_i}(\omega_i)\;(i=1,\ldots,n)\AND
\sum_{i=1}^{n}\norm{x_i^*}=1;
\end{gather}
\item
for any $\eps>0$, there exist points
$\omega_i\in\Omega_i\cap{B}_\eps(\bx)$ and vectors $x_i^*\in X^*$ $(i=1,\ldots,n)$ such that
\begin{gather}\label{D2}\tag{D2}
\sum_{i=1}^n x^*_i=0,\quad
\sum_{i=1}^nd(x^*_i,N_{\Omega_i}(\omega_i))<\eps\AND
\sum_{i=1}^{n}\norm{x_i^*}=1.
\end{gather}
\end{enumerate}
\end{theorem}

Conditions (ii) and (iii) in Theorem~\ref{T2.9} have been used interchangeably (together with several their modifications) since 1980 in the concluding part of the \emph{extremal principle} as dual necessary conditions for extremality
\cite{KruMor80} or local extremality 
\cite{Kru81.2,Kru85.1,KruMor80.2},
first in \emph{Fr\'echet smooth} spaces (i.e. Banach spaces admitting an equivalent norm Fr\'echet differentiable away from zero), and since 1996 \cite{MorSha96} in \emph{Asplund} spaces.
The proof of the extremal principle is based on the applications of the two fundamental variational analysis tools: the Ekeland Variational Principle (EVP) and the fuzzy sum rule for Fr\'echet subdifferentials (Lemma~\ref{SR}(iii)).
We refer the readers to \cite[Section~2.6]{Mor06.1} and \cite{BuiKru18} for more historical comments.

The fact that the generalized separation actually characterizes a weaker than (local) extremality property of approximate stationarity, and this characterization is necessary and sufficient was first established (in a slightly different form) in \cite{Kru98} in the setting of a Fr\'echet smooth Banach space
and extended to Asplund spaces in \cite{Kru00}.
The full proof of the extended extremal principle appeared in \cite{Kru02}, while the name \emph{Extended extremal principle} was introduced in \cite{Kru03}.
It is worth noting that the proof of the necessity (of either condition (ii) or condition (iii)) in the extended extremal principle follows that of the conventional extremal principle and only refines some estimates, while the sufficiency is almost straightforward and is valid in arbitrary normed vector spaces.

\begin{Remark}
\begin{enumerate}
\item
Conditions (ii) and (iii) in Theorem~\ref{T2.9} represent two kinds of widely used \emph{generalized (approximate) separation} of a collection of sets: they claim the existence of $n$ vectors $x_i^*\in X^*$ $(i=1,\ldots,n)$, satisfying the normalization condition $\sum_{i=1}^{n}\norm{x_i^*}=1$, and being either
normal to the respective sets at some points close (up to $\eps$) to $\bx$ with their sum being almost (up to $\eps$) 0 (condition (ii)), or
almost (up to $\eps$) normal with their sum equal exactly 0 (condition (iii)).

The equivalence of the conditions (ii) and (iii) in Theorem~\ref{T2.9} is not difficult to check directly (even in the setting of an arbitrary normed vector space and without assuming the closedness of the sets) using elementary arguments which involve small perturbations of the vectors and then scaling the perturbed vectors to ensure the normalization condition.
Such arguments have been used in many proofs and are scattered across a number of publications.
They were made explicit (in the case of two sets) in \cite[Lemma~1]{BuiKru18}.
In the next subsection, we formulate a more general statement, which will be used also in Sect.~\ref{S6}.
\item
The inequality in \eqref{D2} is only meaningful when $\eps\le1$, because otherwise it is a direct consequence of the equalities.
If condition \eqref{D2} holds with $\eps=1$, it also holds with some $\eps<1$.
Thanks to these observations, when applying Theorem~\ref{T2.9} or its extensions, one can always assume that $\eps<1$.
\item
Theorem~\ref{T2.9} (as well as the conventional extremal principle) with conditions (ii) and (iii) in their current form cannot be extended beyond Asplund spaces.
However,
replacing in these conditions \Fr\ normal cones with normal cones corresponding to other subdifferentials possessing reasonable (approximate or exact) sum rules in the respective \emph{trustworthy} \cite{Iof98} spaces, one can employ basically the same routine to show the necessity of the amended conditions in these spaces; see, e.g., \cite{Iof98,BorJof98,Mor06.1,KruLop12.1} and the discussion in Introduction.
Thanks to Lemma~\ref{SR}(iv) and (i), in general Banach spaces one can use Clarke normal cones or even conventional normal cones in the sense of convex analysis if the sets are convex.
Note that the sufficiency of the conditions (ii) and (iii) in Theorem~\ref{T2.9} is only valid for the Fr\'echet normal cones, though in general normed vector spaces.
Thus, with such extensions, we only have (i) \folgt (ii) \iff (iii), unless the sets are convex.
\end{enumerate}
\end{Remark}

{\bf \subsection{Perturbations and Scaling of Vectors} \label{S2.4}}

Below, we present several assertions containing elementary arguments which are used in proving the equivalence of the conditions (ii) and (iii) in Theorem~\ref{T2.9} and similar facts.
We start with a general statement, which will be used also in Sect.~\ref{S6}.

\begin{lemma}\label{L2.12}
Suppose $K_1,\ldots,K_n$ are cones in a normed vector space, $\eps>0$, $\rho>0$ and $\la>0$.
Suppose also that vectors $z_1,\ldots,z_{n}$ satisfy
\begin{gather}\label{L2.12-1}
\la\sum_{i=1}^nd(z_i,K_i) +\rho\norm{\sum_{i=1}^nz_i}<{\eps},\quad
\sum_{i=1}^n\norm{z_i}=1.
\end{gather}
\begin{enumerate}
\item
If $\eps+\la\le\rho$, then there exist vectors $\hat z_i$ $(i=1,\ldots,n)$ satisfying the following conditions:
\begin{gather}\label{L2.12-2}
\sum_{i=1}^n\hat z_i=0,\quad
\sum_{i=1}^nd(\hat z_i,K_i)< \frac{\eps}{\la},\quad
\sum_{i=1}^{n}\norm{\hat z_i}=1.
\end{gather}
\item
If $\eps+\rho\le\la$, then there exist vectors $\hat z_i$ $(i=1,\ldots,n)$ satisfying the following conditions:
\begin{gather}\label{L2.12-3}
\norm{\sum_{i=1}^n\hat z_i}<\frac{\eps}{\rho},\quad
\hat z_i\in K_i\; (i=1,\ldots,n),\quad
\sum_{i=1}^{n}\norm{\hat z_i}=1.
\end{gather}
\item
Moreover, if the underlying space is dual to a normed vector space $X$, and
\begin{gather}\label{L2.12-4}
\sum_{i=1}^n\ang{z_{i},x_{i}}\ge\tau\max_{1\le{i}\le{n}} \|x_{i}\|
\end{gather}
for some vectors $x_i\in X$ $(i=1,\ldots,n)$, not all zero, and a number $\tau\in]0,1]$, then the vectors $\hat z_i$ $(i=1,\ldots,n)$ in parts {\rm (i)} or {\rm (ii)} satisfy
\begin{gather}\label{L2.12-5}
\sum_{i=1}^n\ang{\hat z_{i},x_{i}} >\hat\tau\max_{1\le{i}\le{n}} \|x_{i}\|,
\end{gather}
where $\hat\tau:=\frac{\tau\rho-\eps}{\rho+\eps}$ under the assumptions in part {\rm (i)}, and $\hat\tau :=\frac{\tau\la-\eps}{\la+\eps}$ under the assumptions in part {\rm (ii)}.
\end{enumerate}
\end{lemma}

\begin{Proof}
\begin{enumerate}
\item
Let $\eps+\la\le\rho$.
Set $z:=\sum_{i=1}^nz_i$ and $v_i:=z_i-\frac{1}{n}z$ $(i=1,\ldots,n)$.
Then, $\sum_{i=1}^nv_i=0$ and, by \eqref{L2.12-1}, $\|z\|<\frac{\eps}{\rho}<1$, 
\begin{align*}
\sum_{i=1}^n\|v_i\|&\le\sum_{i=1}^n\|z_i\|+\|z\| <1+\frac{\eps}{\rho}\AND 
\\
\sum_{i=1}^n\|v_i\|&\ge\sum_{i=1}^n\|z_i\|-\|z\| =1-\|z\|>0.
\end{align*}
Set $\hat z_i:=v_i/\sum_{i=1}^n\|v_i\|$.
Then, $\sum_{i=1}^n\hat z_i=0$, $\sum_{i=1}^n\|\hat z_i\|=1$ and
\begin{align*}
\sum_{i=1}^nd(\hat z_i,K_i) &\le\frac{\sum_{i=1}^nd(z_i,K_i) +\|z\|} {\sum_{i=1}^n\|v_i\|} =\frac{\la\sum_{i=1}^nd(z_i,K_i) +\rho\|z\|-(\rho-\la)\|z\|} {\la\sum_{i=1}^n\|v_i\|}
\\
&<\frac{\eps-(\rho-\la)\|z\|} {\la(1-\|z\|)} \le\frac{\eps-\eps\|z\|} {\la(1-\|z\|)}=\frac{\eps}{\la},
\end{align*}
i.e. all conditions in \eqref{L2.12-2} are satisfied.
\item
Let $\eps+\rho\le\la$.
By \eqref{L2.12-1}, there exist vectors $v_i\in K_i$ $(i=1,\ldots,n)$ such that
\begin{gather*}
\la\sum_{i=1}^n\|z_i-v_i\| +\rho\norm{\sum_{i=1}^nz_i}<\eps.
\end{gather*}
In particular,
$
\sum_{i=1}^n\norm{z_i-v_i}<\frac{\eps}{\la}<1.
$
Hence,
\begin{align*}
\sum_{i=1}^n\norm{v_i}&\le \sum_{i=1}^n\norm{z_i}+ \sum_{i=1}^n\norm{z_i-v_i}=1 +\sum_{i=1}^n\norm{z_i-v_i},
\\
\sum_{i=1}^n\norm{v_i}&\ge \sum_{i=1}^n\norm{z_i}- \sum_{i=1}^n\norm{z_i-v_i}=1- \sum_{i=1}^n\norm{z_i-v_i}>0,
\\
\norm{\sum_{i=1}^n v_i} &\le\sum_{i=1}^n \norm{z_i-v_i} +\norm{\sum_{i=1}^nz_i} <\frac{1}{\rho}\left(\eps-\left(\la-\rho\right) \sum_{i=1}^n \norm{z_i-v_i}\right)
\\
&\le\frac{\eps}{\rho}\left(1-\sum_{i=1}^n \norm{z_i-v_i}\right).
\end{align*}
Set $\hat z_i:=v_i/\sum_{i=1}^n\|v_i\|$.
Then, $\hat z_i\in K_i$ $(i=1,\ldots,n)$, $\sum_{i=1}^n\|\hat z_i\|=1$ and
$$
\norm{\sum_{i=1}^n\hat z_i} =\frac{\norm{\sum_{i=1}^nv_i}} {\sum_{i=1}^n\norm{v_i}} <\frac{\eps(1-\sum_{i=1}^n\norm{z_i-v_i})} {\rho(1-\sum_{i=1}^n\norm{z_i-v_i})}
=\frac{\eps}{\rho},
$$
i.e. all conditions in \eqref{L2.12-3} are satisfied.
\item
Suppose that the underlying space is dual to a normed vector space $X$ and condition \eqref{L2.12-4} is satisfied
for some vectors $x_i\in X$ $(i=1,\ldots,n)$, not all zero, and a number $\tau\in]0,1]$.
Then, using the notations introduced above, we have:
\begin{align*}
\sum_{i=1}^n\ang{\hat z_{i},x_i} &\ge\frac{\sum_{i=1}^n\ang{z_i,x_i} -(\sum_{i=1}^n\|z_{i}-v_i\|) \max_{1\le{i}\le{n}}\|x_i\|} {\sum_{i=1}^n \norm{v_i}}
\\
&\ge\frac{\tau-\sum_{i=1}^n\|z_{i}-v_i\|} {\sum_{i=1}^n \norm{v_i}} \max_{1\le{i}\le{n}}\|x_i\|.
\end{align*}
Employing the estimates in part (i), we obtain
\begin{align*}
\frac{\tau-\sum_{i=1}^n\|z_{i}-v_i\|} {\sum_{i=1}^n \norm{v_i}}&=\frac{\tau-\|z\|} {\sum_{i=1}^n \norm{v_i}}
>\frac{\tau-\frac{\eps}{\rho}} {1+\frac{\eps}{\rho}} =\frac{\tau\rho-\eps}{\rho+\eps},
\end{align*}
while the estimates in part (ii) give
\begin{align*}
\frac{\tau-\sum_{i=1}^n\|z_{i}-v_i\|} {\sum_{i=1}^n \norm{v_i}}&\ge\frac{\tau-\sum_{i=1}^n\|z_{i}-v_i\|} {1+\sum_{i=1}^n\|z_{i}-v_i\|}
>\frac{\tau-\frac{\eps}{\la}} {1+\frac{\eps}{\la}}
=\frac{\tau\la-\eps}{\rho+\la}.
\end{align*}
Thus, in both cases we arrive at \eqref{L2.12-5}.
\qed\end{enumerate}
\end{Proof}

The next two corollaries present two important special cases of Lemma~\ref{L2.12}.

\begin{corollary}\label{C2.13}
Suppose $K_1,\ldots,K_n$ are cones in a normed vector space, $\eps\in]0,1[$, and vectors $z_1,\ldots,z_{n}$ satisfy
\begin{gather*}
\norm{\sum_{i=1}^nz_i}<{\eps},\quad
z_i\in K_i\; (i=1,\ldots,n),\quad \sum_{i=1}^n\norm{z_i}=1.
\end{gather*}
Then, there exist vectors $\hat z_i$ $(i=1,\ldots,n)$ satisfying the following conditions:
\begin{gather*}
\sum_{i=1}^n\hat z_i=0,\quad
\sum_{i=1}^nd(\hat z_i,K_i)< \frac{\eps}{1-\eps},\quad
\sum_{i=1}^{n}\norm{\hat z_i}=1.
\end{gather*}
Moreover, if the underlying space is dual to a normed vector space $X$, and condition \eqref{L2.12-4} is satisfied
for some vectors $x_i\in X$ $(i=1,\ldots,n)$, not all zero, and a number $\tau\in]0,1]$, then the vectors $\hat z_i$ $(i=1,\ldots,n)$ satisfy condition \eqref{L2.12-5} with
$\hat\tau :=\frac{\tau-\eps}{1+\eps}$.
\end{corollary}

\begin{Proof}
Apply Lemma~\ref{L2.12}(i) with $\rho=1$ and $\la=1-\eps$.
\qed\end{Proof}

\begin{corollary}\label{C2.14}
Suppose $K_1,\ldots,K_n$ are cones in a normed vector space, $\eps\in]0,1[$, and vectors $z_1,\ldots,z_{n}$ satisfy
\begin{gather*}
\sum_{i=1}^nz_i=0,\quad
\sum_{i=1}^nd(z_i,K_i)<{\eps},\quad \sum_{i=1}^n\norm{z_i}=1.
\end{gather*}
Then, there exist vectors $\hat z_i$ $(i=1,\ldots,n)$ satisfying the following conditions:
\begin{gather}\label{C2.14-2}
\norm{\sum_{i=1}^n\hat z_i}<\frac{\eps}{1-\eps},\quad
\hat z_i\in K_i\; (i=1,\ldots,n),\quad
\sum_{i=1}^{n}\norm{\hat z_i}=1.
\end{gather}
Moreover, if the underlying space is dual to a normed vector space $X$, and condition \eqref{L2.12-4} is satisfied
for some vectors $x_i\in X$ $(i=1,\ldots,n)$, not all zero, and a number $\tau\in]0,1]$, then the vectors $\hat z_i$ $(i=1,\ldots,n)$ satisfy condition \eqref{L2.12-5} with
$\hat\tau :=\frac{\tau-\eps}{1+\eps}$.
\end{corollary}

\begin{Proof}
Apply Lemma~\ref{L2.12}(ii) with $\rho=1-\eps$ and $\la=1$.
\qed\end{Proof}

As an immediate consequence of Corollaries~\ref{C2.13} and \ref{C2.14}, we obtain the following important assertion.

\begin{corollary}\label{C2.20}
Suppose $\Omega_1,\ldots,\Omega_n$ are subsets of a normed vector space $X$ and
$\bx\in\cap_{i=1}^n\Omega_i$.
Conditions {\rm (ii)} and {\rm (iii)} in Theorem~\ref{T2.9} are equivalent.
\end{corollary}

When $n=2$, the main estimate in Corollary~\ref{C2.13} can be improved.

\begin{proposition}
Suppose $K_1$ and $K_2$ are cones in a normed vector space, $\eps>0$, and vectors $z_1$ and $z_2$ satisfy
\begin{gather*}
\norm{z_1+z_2}<{\eps},\quad
z_1\in K_1,\; z_2\in K_2,\quad \norm{z_1}+\norm{z_2}=1.
\end{gather*}
Then, there exist vectors $\hat z_1$ and $\hat z_2$ satisfying the following conditions:
\begin{gather*}
\hat z_1+\hat z_2=0,\quad
\hat z_1\in K_1,\quad d(\hat z_2,K_2)<\eps,\quad
\norm{\hat z_1}+\norm{\hat z_2}=1.
\end{gather*}
\end{proposition}

\begin{Proof}
Without loss of generality, we can assume that $\norm{z_2}\le\frac{1}{2}\le\norm{z_1}$.
Set $\hat z_1:=\frac{z_1}{2\norm{z_1}}$ and $\hat z_2:=-\hat z_1$.
Then $\hat z_1+\hat z_2=0$, $\norm{\hat z_1}=\norm{\hat z_2} =\frac{1}{2}$, $\hat z_1\in K_1$ and
$$
d(\hat z_2,K_2)\le\norm{\hat z_2 -\frac{z_2}{2\norm{z_1}}} =\frac{\norm{z_1+z_2}}{2\norm{z_1}} <\frac{\eps}{2\norm{z_1}}\le\eps.
$$
This completes the proof.
\qed\end{Proof}

Next, we formulate an asymmetric modification of Corollary~\ref{C2.14} which will be used in Sect.~\ref{S6}.

\begin{proposition}\label{P2.15}
Suppose $K_1,\ldots,K_n$ are cones in a normed vector space, $\eps\in]0,1[$ and vectors $z_1,\ldots,z_{n}$ satisfy conditions
\begin{gather*}
\sum_{i=1}^nz_i=0,\quad
\sum_{i=1}^nd(z_i,K_i)<{\eps},\quad \sum_{i=1}^{n-1}\norm{z_i}=1.
\end{gather*}
Then, there exist vectors $\hat z_i$ $(i=1,\ldots,n)$ satisfying the following conditions:
\begin{gather}\label{P2.15-2}
\sum_{i=1}^{n}\hat z_i=0,\quad
\hat z_i\in K_i\; (i=1,\ldots,n-1),\quad
d\left(\hat z_n,K_n\right)<\frac{\eps}{1-\eps}\AND
\sum_{i=1}^{n-1}\norm{\hat z_i}=1.
\end{gather}
\end{proposition}

\begin{Proof}
Take $y_i \in K_i$ ($i=1,\ldots,n$) such that
$
\sum_{i=1}^{n}\norm{z_i- y_i}<\eps.
$
Then
$\norm{\sum_{i=1}^{n}y_i}<\eps$ and $\sum_{i=1}^{n-1}\norm{y_i}>1-\eps.$
It follow that $d(-\sum_{i=1}^{n-1}y_i,K_n)\le d(-\sum_{i=1}^{n-1}y_i,y_n)<\eps$.
Hence, vectors $\hat z_i:=y_i/\sum_{i=1}^{n-1} \norm{y_i}$ ($i=1,\ldots,n-1$) and $\hat z_n:=-\sum_{i=1}^{n-1}\hat z_i$.
satisfy all the conditions in \eqref{P2.15-2}.
\qed\end{Proof}

\section{Modifications of the Extremality and Stationarity Properties}\label{S3}

There exist several modifications of the properties discussed in Sect.~\ref{S2.3}, scattered in the literature.
Below, we briefly discuss some of them which are going to be important for our subsequent study.

The properties in Definition~\ref{D2.5} involve translations of all the sets.
It is easy to see that in all the properties it is sufficient to consider translations of all but one sets leaving the remaining set unchanged.
This simple observation leads to asymmetric conditions in the next proposition which can be useful, especially in the case $n=2$.

\begin{proposition}\label{P2.7}
Suppose $\Omega_1,\ldots,\Omega_n$ are subsets of a normed vector space $X$ and $\bx\in\cap_{i=1}^n\Omega_i$.
The collection $\{\Omega_1,\ldots,\Omega_n\}$ is \begin{enumerate}
\item
extremal if and only if,
for any $\eps>0$, there exist vectors
$a_i\in{X}$ $(i=1,\ldots,n-1)$ satisfying
\begin{gather}\label{P4}\tag{P4}
\bigcap_{i=1}^{n-1}(\Omega_i-a_i)\cap\Omega_n =\emptyset\AND
\max_{1\le i\le n-1}\norm{a_i}<\eps;
\end{gather}
\item
locally extremal at $\bx$ if and only if there exists a number $\rho\in]0,\infty]$ such that,
for any $\eps>0$, there exist vectors
$a_i\in{X}$ $(i=1,\ldots,n-1)$ satisfying
\begin{gather}\label{P5}\tag{P5}
\bigcap_{i=1}^{n-1}(\Omega_i-a_i) \cap\Omega_n\cap{B}_\rho(\bar{x})=\emptyset\AND
\max_{1\le i\le n-1}\norm{a_i}<\eps;
\end{gather}
\item
stationary at $\bx$ if and only if,
for any $\eps>0$,
there exist a number $\rho\in]0,\eps[$
and vectors
$a_i\in{X}$ $(i=1,\ldots,n-1)$ satisfying
\begin{gather*}
\bigcap_{i=1}^{n-1}(\Omega_i-a_i)\cap\Omega_n \cap{B}_\rho(\bar{x})
=\emptyset\AND
\max_{1\le i\le n-1}\norm{a_i}<\eps\rho;
\end{gather*}
\item
approximately stationary at $\bx$ if and only if,
for any $\eps>0$,
there exist a number $\rho\in]0,\eps[$, points $\omega_i\in\Omega_i\cap B_\eps(\bx)$ $(i=1,\ldots,n)$
and vectors
$a_i\in{X}$ $(i=1,\ldots,n-1)$ satisfying
\begin{gather}\label{P6}\tag{P6}
\bigcap_{i=1}^{n-1}(\Omega_i-\omega_i-a_i) \cap(\Omega_n-\omega_n) \cap(\rho\B)
=\emptyset\AND
\max_{1\le i\le n-1}\norm{a_i}<\eps\rho.
\end{gather}
\end{enumerate}
\end{proposition}

From now on, we mostly focus on characterizations of approximate stationarity.

Observe that, by scaling the vectors, the normalization condition $\sum_{i=1}^{n}\norm{x_i^*}=1$ in the dual generalized separation properties \eqref{D1} and \eqref{D2} can be dropped if the inequalities there are amended to
$$\sum_{i=1}^nd(x^*_i,N_{\Omega_i}(\omega_i)) <\eps\sum_{i=1}^{n}\norm{x_i^*}\AND\norm{\sum_{i=1}^n x^*_i} <\eps\sum_{i=1}^{n}\norm{x_i^*},$$
respectively.
Note that each of the amended inequalities still implies that $\sum_{i=1}^{n}\norm{x_i^*}>0$.
Moreover, the normalization condition $\sum_{i=1}^{n}\norm{x_i^*}=1$ involving $n$ vectors can be replaced by the similar asymmetric condition involving $n-1$ vectors: $\sum_{i=1}^{n-1}\norm{x_i^*}=1$.
This observation can be especially useful in the case $n=2$.

Observe further that parameter $\eps$ in Definition~\ref{D2.5}(iv) of approximate stationarity and its reformulation in Proposition~\ref{P2.7}(iv) as well as the metric and dual characterizations in Theorems~\ref{T2.8}(ii) and (iii) and \ref{T2.9}(ii) and (iii) plays multiple roles.
To get a deeper insight into the approximate stationarity property, it makes sense to split the parameter $\eps$ into two components.
From now on, we will use the letter $\al$ to denote the component controlling the size of the shifts of the sets.
This parameter is going to be crucial for quantifying the corresponding transversality property, playing the role of the rate/modulus of the property.

Based on the above observations, we now formulate a list of primal and dual equivalent characterizations of approximate stationarity, complementing Definition~\ref{D2.5}(iv), Proposition~\ref{P2.7}(iv), Theorem~\ref{T2.8}(ii) and (iii) and Theorem~\ref{T2.9}(ii) and (iii), which will be used in the sequel.

\begin{proposition}\label{P2.11}
Suppose $\Omega_1,\ldots,\Omega_n$ are subsets of a normed vector space $X$ and $\bx\in\cap_{i=1}^n\Omega_i$.
The following conditions are equivalent:
\begin{enumerate}
\item
the collection $\{\Omega_1,\ldots,\Omega_n\}$ is approximately stationary at $\bx$;
\item
for any $\eps>0$ and $\al>0$,
there exist a number $\rho\in]0,\eps[$, points $\omega_i\in\Omega_i\cap B_\eps(\bx)$
and vectors
$a_i\in{X}$ $(i=1,\ldots,n)$ satisfying
\begin{gather}\tag{P7}\label{P7}
\bigcap_{i=1}^n(\Omega_i-\omega_i-a_i)\cap(\rho\B)
=\emptyset\AND
\max_{1\le i\le n}\norm{a_i}<\al\rho;
\end{gather}
\item
for any $\eps>0$ and $\al>0$,
there exist a number $\rho\in]0,\eps[$, points $\omega_i\in\Omega_i\cap B_\eps(\bx)$ $(i=1,\ldots,n)$
and vectors
$a_i\in{X}$ $(i=1,\ldots,n-1)$ satisfying
\begin{gather}\tag{P8}\label{P8}
\bigcap_{i=1}^{n-1}(\Omega_i-\omega_i-a_i) \cap(\Omega_n-\omega_n) \cap(\rho\B)
=\emptyset\AND
\max_{1\le i\le n-1}\norm{a_i}<\al\rho;
\end{gather}
\item
for any $\eps>0$ and $\al>0$, there exist vectors $a_i\in X$ $(i=1,\ldots,n)$ such that
\begin{gather}\tag{P9}\label{P9}
\al d\left(\bx,\bigcap_{i=1}^n(\Omega_i-a_i)\right)> \max_{1\le i\le n} d\left(\bx,\Omega_i-a_i\right)
\qdtx{and}
\max_{1\le i\le n}\norm{a_i}<\eps;
\end{gather}
\item
for any $\eps>0$ and $\al>0$, there exist a point $x\in B_\eps(\bx)$ and vectors $a_i\in X$ $(i=1,\ldots,n)$ such that
\begin{gather}\tag{P10}\label{P10}
\al d\left(x,\bigcap_{i=1}^n(\Omega_i-a_i)\right) >\max_{1\le i\le n}d(x,\Omega_i-a_i)
\qdtx{and}
\max_{1\le i\le n}\norm{a_i}<\eps.
\end{gather}
\cnta
\end{enumerate}
With $N$ standing for either Clarke ($N:=N^C$) or \Fr\ ($N:=N^F$) normal cone, the following conditions are equivalent to conditions {\rm (ii)} and {\rm (iii)} in Theorem~\ref{T2.9}:
\begin{enumerate}
\cntb
\item
for any $\eps>0$ and $\al>0$, there exist points
$\omega_i\in\Omega_i\cap{B}_\eps(\bx)$ and vectors $x_i^*\in X^*$ $(i=1,\ldots,n)$ such that
\begin{gather}\tag{D3}\label{D3}
\norm{\sum_{i=1}^n x^*_i}<\al,\quad
x^*_i\in N_{\Omega_i}(\omega_i)\;(i=1,\ldots,n)\AND
\sum_{i=1}^{n}\norm{x_i^*}=1;
\end{gather}
\item
for any $\eps>0$ and $\al>0$, there exist points
$\omega_i\in\Omega_i\cap{B}_\eps(\bx)$ and vectors $x_i^*\in X^*$ $(i=1,\ldots,n)$ such that
\begin{gather}\tag{D4}\label{D4}
\sum_{i=1}^n x^*_i=0,\quad
\sum_{i=1}^nd(x^*_i,N_{\Omega_i}(\omega_i))<\al\AND
\sum_{i=1}^{n}\norm{x_i^*}=1;
\end{gather}
\item
for any $\eps>0$ and $\al>0$, there exist points
$\omega_i\in\Omega_i\cap{B}_\eps(\bx)$ and vectors $x_i^*\in X^*$ $(i=1,\ldots,n)$ such that
\begin{gather}\tag{D5}\label{D5}
\norm{\sum_{i=1}^n x^*_i}<\al,\quad
x^*_i\in N_{\Omega_i}(\omega_i)\;(i=1,\ldots,n)\AND
\sum_{i=1}^{n-1}\norm{x_i^*}=1;
\end{gather}
\item
for any $\eps>0$ and $\al>0$, there exist points
$\omega_i\in\Omega_i\cap{B}_\eps(\bx)$ and vectors $x_i^*\in X^*$ $(i=1,\ldots,n)$ such that
\begin{gather}\tag{D6}\label{D6}
\sum_{i=1}^n x^*_i=0,\quad
\sum_{i=1}^nd(x^*_i,N_{\Omega_i}(\omega_i))<\al\AND
\sum_{i=1}^{n-1}\norm{x_i^*}=1.
\end{gather}
\end{enumerate}
If $X$ is Asplund and $N=N^F$, then all conditions {\rm (i)--(ix)} are equivalent.
\end{proposition}

\begin{Remark}
The maximum in each of the conditions in Definition~\ref{D2.5}, Theorem~\ref{T2.8}, Proposition~\ref{P2.7} and the first (primal space) part of Proposition~\ref{P2.11} can be replaced by the sum.
Moreover, any norm on $\R^n$ (or $\R^{n-1}$ in the case of Proposition~\ref{P2.7} and Proposition~\ref{P2.11}(iii)) can be used instead.
The sum of the norms in the normalization conditions in \eqref{D1} and \eqref{D2} and the second (dual space) part of Proposition~\ref{P2.11} stands for the corresponding dual norm and can be replaced by the maximum,
or any other norm on $\R^n$ or $\R^{n-1}$.
\end{Remark}

The `$\al$-version' of the approximate stationarity based on its equivalent representation in Proposition~\ref{P2.11}(i) is going to be used in the subsequent study.

\begin{definition}[Approximate $\al$-stationarity] \label{D2.17}
Suppose $\Omega_1,\ldots,\Omega_n$ are subsets of a normed vector space $X$, $\bx\in\cap_{i=1}^n\Omega_i$, and $\al>0$.
The collection $\{\Omega_1,\ldots,\Omega_n\}$ is approximately $\al$-stationary at $\bx$ iff
for any $\eps>0$,
there exist a number $\rho\in]0,\eps[$, points $\omega_i\in\Omega_i\cap B_\eps(\bx)$
and vectors
$a_i\in{X}$ $(i=1,\ldots,n)$ satisfying conditions \eqref{P7}.
\end{definition}

In view of the above discussion, the next assertion is straightforward.

\begin{proposition}\label{P2.18}
Suppose $\Omega_1,\ldots,\Omega_n$ are subsets of a normed vector space $X$ and $\bx\in\cap_{i=1}^n\Omega_i$.
The collection $\{\Omega_1,\ldots,\Omega_n\}$ is approximately stationary at $\bx$ if and only if it is approximately $\al$-stationary at $\bx$ for all $\al>0$.
\end{proposition}

In our deeper analysis of the core arguments in proofs of metric and dual characterizations of the extremality and stationarity in Sects.~\ref{S5} and \ref{S6}, we will study properties which correspond to fixing, besides $\al$, also other parameters involved in the primal and dual properties discussed above.

Theorems~\ref{T2.8} and \ref{T2.9} as well as the equivalent characterizations in Proposition~\ref{P2.11} can be `reversed' into statements providing primal and dual space  criteria for the absence of the approximate stationarity, which turns out to be an important \emph{regularity/transversality} property of collections of sets,
which plays an important role in constraint qualifications, qualification conditions in subdifferential/coderivative calculus and convergence analysis of computational algorithms \cite{Kru05,Kru06,Kru09,KruTha13,KruTha15,KruLukTha18, LewLukMal09, BauLukPhaWan13.2,HesLuk13,DruIofLew15, Iof17}.

\begin{definition}[Transversality]\label{D2.15}
Suppose $\Omega_1,\ldots,\Omega_n$ are subsets of a normed vector space $X$ and $\bx\in\cap_{i=1}^n\Omega_i$.
The collection $\{\Omega_1,\ldots,\Omega_n\}$ is
transversal at $\bx$ iff there exist numbers $\al>0$ and $\eps>0$ such that
\begin{gather}\label{D2.15-1}
\bigcap_{i=1}^n(\Omega_i-\omega_i-a_i)\cap(\rho\B)
\ne\emptyset
\end{gather}
for all numbers $\rho\in]0,\eps[$, points $\omega_i\in\Omega_i\cap B_\eps(\bx)$
and vectors
$a_i\in{X}$ $(i=1,\ldots,n)$ satisfying $\max_{1\le i\le n}\norm{a_i}<\al\rho$.
\sloppy
\end{definition}

The property in the above definition was first considered explicitly in \cite{Kru05}.
It has been renamed multiple times (\emph{regularity, strong regularity, property} (UR)$_S$, \emph{uniform regularity}, and finally \emph{transversality}).
A table illustrating the evolution of the terminology can be found in \cite[Section~2]{KruLukTha17}.
The transversality property of collections of sets happens to be a direct counterpart of the \emph{metric regularity} property of \SVM s; cf. \cite{Kru05,Kru06,Kru09}.

The number $\al$ involved in the definition provides a quantitative characterization of the property.
The supremum of all such numbers (with the convention that the supremum of the empty subset of $\R_+$ equals 0), denoted $\Tr$, is the \emph{modulus of transversality} of the collection $\Omega_1,\ldots,\Omega_n$ at $\bx$.
Thus, the case $\Tr=0$ corresponds to approximate stationarity.
\sloppy

In view of Proposition~\ref{P2.11},
the next two statements are direct consequences of Theorems~\ref{T2.8} and \ref{T2.9}, respectively.

\begin{theorem}[Transversality: metric criteria]\label{T2.13}
Suppose $\Omega_1,\ldots,\Omega_n$ are subsets of a normed vector space $X$ and ${\bx\in\cap_{i=1}^n\Omega_i}$.
The following conditions are equivalent:
\begin{enumerate}
\item
the collection $\{\Omega_1,\ldots,\Omega_n\}$ is transversal at $\bx$;
\item
there exist numbers $\al>0$ and $\eps>0$ such that
\begin{equation*}
\al d\left(\bx,\bigcap_{i=1}^n(\Omega_i-a_i)\right) \le\max_{1\le i\le n} d\left(\bx,\Omega_i-a_i\right)
\end{equation*}
for all vectors $a_i\in X$ $(i=1,\ldots,n)$ satisfying $\max_{1\le i\le n}\norm{a_i}<\eps$;
\item
there exist numbers $\al>0$ and $\eps>0$ such that
\begin{equation*}
\al d\left(x,\bigcap_{i=1}^n(\Omega_i-a_i)\right) \le\max_{1\le i\le n}d(x,\Omega_i-a_i)
\end{equation*}
for all points $x\in B_\eps(\bx)$ and vectors $a_i\in X$ $(i=1,\ldots,n)$ satisfying $\max_{1\le i\le n}\norm{a_i}<\eps$.
\end{enumerate}
\end{theorem}

The equivalence of conditions (i) and (iii) in Theorem~\ref{T2.13} recaptures \cite[Theorem~1]{Kru05}; cf. \cite[Theorem~1]{Kru06} and \cite[Definition~2]{KruLukTha17}.
Adding condition (ii) to the list of equivalent conditions sheds additional light on the transversality property.

\begin{theorem}[Transversality: dual criteria]\label{T2.14}
Suppose $\Omega_1,\ldots,\Omega_n$ are closed subsets of an Asplund space $X$ and
${\bx\in\cap_{i=1}^n\Omega_i}$.
With $N$ standing for the \Fr\ normal cone ($N:=N^F$), the following conditions are equivalent:
\begin{enumerate}
\item
the collection $\{\Omega_1,\ldots,\Omega_n\}$ is transversal at $\bx$;
\item
there exist numbers $\al>0$ and $\eps>0$ such that
$\norm{\sum_{i=1}^n x^*_i}>\al$
for all points
$\omega_i\in\Omega_i\cap{B}_\eps(\bx)$ and vectors $x_i^*\in N_{\Omega_i}(\omega_i)$ $(i=1,\ldots,n)$  satisfying $\sum_{i=1}^{n}\norm{x_i^*}=1$;
\sloppy
\item
there exist numbers $\al>0$ and $\eps>0$ such that
$\sum_{i=1}^nd(x^*_i,N_{\Omega_i}(\omega_i))>\al$
for all points
$\omega_i\in\Omega_i\cap{B}_\eps(\bx)$ and vectors $x_i^*\in X^*$ $(i=1,\ldots,n)$ satisfying $\sum_{i=1}^n x^*_i=0$ and $\sum_{i=1}^{n}\norm{x_i^*}=1$.
\end{enumerate}
\end{theorem}

In view of Theorem~\ref{T2.14},
transversality of a collection of sets is equivalent to the absence of the generalized separation.

\begin{Remark}\label{R3.9}
\begin{enumerate}
\item
In view of Corollaries~\ref{C2.13} and \ref{C2.14},
conditions {\rm (ii)} and {\rm (iii)} in Theorem~\ref{T2.14} are equivalent.
\item
The supremums of all numbers $\al$ in parts (ii) and (iii) of Theorem~\ref{T2.13} and part (ii) of Theorem~\ref{T2.14} equal $\Tr$; see \cite[Theorem~1]{Kru05}, \cite[Theorem~4(vi)]{Kru09} and Corollary~\ref{C4.9}.
The supremum of all numbers $\al$ in part (iii) of Theorem~\ref{T2.14} can be different from $\Tr$, but its relationship with $\Tr$ can be easily established using elementary arguments discussed in Sect. \ref{S2.4}.
\end{enumerate}
\end{Remark}

\if{
\AK{4/07/18.
Can the relation between the two be estimated?}
}\fi

\begin{example}
In the space $\R^2$ equipped with the maximum norm (hence, the dual norm is the sum norm),
consider the two perpendicular lines: $\Omega_1:=\{(t,0):\, t\in \R\}$ and $\Omega_2:=\{(0,t):\, t\in \R\}$.
Then, we have $\bx:=(0,0)\in\Omega_1\cap\Omega_2$,
$N_{\Omega_1}^F(x)=\{(0,t):\, t\in \R\}$ for any $x\in\Omega_1$ and $N_{\Omega_1}^F(x)=\{(t,0):\, t\in \R\}$ for any $x\in\Omega_2$.

If $x_1^*=(0,t_1)$ and $x_2^*=(t_2,0)$ are normal vectors to $\Omega_1$ and $\Omega_2$, respectively, then $\norm{x_1^*}+\norm{x_2^*}=\norm{x_1^*+x_2^*} =|t_1|+|t_2|$.
Hence, the supremum of all $\al$ in part (ii) of Theorem~\ref{T2.14} is $1$ (and is equal to tr$[\Omega_1,\Omega_2]$).

If $x_1^*,x_2^*\in\left(\R^2\right)^*$, $x_1^*+x_2^*=0$ and $\norm{x_1^*}+\norm{x_2^*}=1$, then $x_1^*=-x_2^*=(t_1,t_2)$ for some numbers $t_1$ and $t_2$ satisfying $|t_1|+|t_2|=\frac{1}{2}$, and the distances from $x_1^*$ and $x_2^*$ to the corresponding normal cones equal $|t_1|$ and $|t_2|$, respectively.
Hence, the supremum of all $\al$ in part (iii) of Theorem~\ref{T2.14} is $\frac{1}{2}$.
\end{example}

\begin{Remark}
Since the approximate stationarity and transversality properties are complementary to each other,
it would be natural to refer to the negation of the approximate $\al$-stationarity property, i.e. the existence of a number $\eps>0$ such that condition \eqref{D2.15-1} holds
for all numbers $\rho\in]0,\eps[$, points $\omega_i\in\Omega_i\cap B_\eps(\bx)$
and vectors
$a_i\in{X}$ $(i=1,\ldots,n)$ satisfying $\max_{1\le i\le n}\norm{a_i}<\al\rho$, as \emph{$\al$-transversality} at $\bx$.
We do not consider this property in the current paper.
\sloppy
\end{Remark}

\section{Distances Between $n$ Sets}\label{S4}

When studying mutual arrangement of collections of sets in space, particularly their extremality, stationarity and regularity properties, we need to be able to estimate the `distance' between the sets, i.e. how `far apart' they are, or, at least, whether they have a common point.
In the case of two sets, the conventional distance \begin{equation}\label{d0}
d(\Omega_1,\Omega_2):= \inf_{\omega_1\in\Omega_1,\;\omega_2\in\Omega_2} d(\omega_1,\omega_2)
\end{equation}
does the job.
The general case of $n\ge2$ sets is not that straightforward.
In this section, we discuss several candidates for the role of `distance'.

We start with discussing distances between $n\ge2$ points.

\subsection{Distances Between $n$ Points}

Our aim in this subsection is to consider ways of defining \emph{$n$-point distances} estimating quantitatively the overall `closeness' of a collection of $n\ge2$ points in a metric space $(X,d)$.

Given $n\ge2$ points $\omega_1,\ldots,\omega_n$ in $X$,
one can use one of the following two quantities:
\begin{gather}\label{d1}
d_1(\omega_1,\ldots,\omega_n):=\max_{1\le i\le n-1} d(\omega_i,\omega_n),
\\\label{d2}
d_2(\omega_1,\ldots,\omega_n):=\inf_{x\in X}\max_{1\le i\le n} d(\omega_i,x).
\end{gather}
When $n=2$, \eqref{d1} reduces to the conventional distance $d(\omega_1,\omega_2)$.
However, when $n>2$, this distance is not symmetric: the last point in the list plays a special role, and the quantity itself depends on the choice of the last point.
In contrast to \eqref{d1}, in definition \eqref{d2} all points play the same role, and it involves minimization over an additional parameter $x$.
Formulas \eqref{d1} and \eqref{d2} produce in general different numbers, whatever the choice of the last point in \eqref{d1} is.

In the setting of a normed vector space,
the following symmetric distance
can be of interest:
\begin{gather}\label{d3}
d_3(\omega_1,\ldots,\omega_n) :=\max_{1\le i\le n} \norm{\omega_i-\frac{1}{n}\sum_{j=1}^n\omega_j}.
\end{gather}

\begin{example}\label{E3.1}
Consider three points in $\R$: $0$, $1$ and $5$.
By \eqref{d1}, we have $d_1(0,5,1)=\max\{|0-1|,|5-1|\}=4$, while
$d_1(0,1,5)=d_1(1,5,0)=5$.
The infimum in definition \eqref{d2} is achieved at $x=2.5$ and equals $d_2(0,1,5)=\max\{|0-2.5|,|1-2.5|,|5-2.5|\}=2.5$.
The average $\frac{1}{3}(0+1+5)$  equals 2, and formula \eqref{d3} gives $d_3(0,1,5)=\max\{|0-2|,|1-2|,|5-2|\}=3$.
Thus, for these three points all three definitions \eqref{d1}, \eqref{d2} and \eqref{d3} give different numbers.
\end{example}

Note the obvious connection between the distances $d_1$ and $d_2$:
\begin{equation*}
d_2(\omega_1,\ldots,\omega_n):=\inf_{x\in X} d_1(\omega_1,\ldots,\omega_n,x).
\end{equation*}
Observe also that $d_1(\omega_1,\ldots,\omega_n)$ is actually the usual (maximum) distance between the points $(\omega_1,\ldots,\omega_{n-1})$ and $(\omega_n,\ldots,\omega_n)$ in $X^{n-1}$, while $d_2(\omega_1,\ldots,\omega_n)$ represents the distance from the points $(\omega_1,\ldots,\omega_{n})\in X^n$ to the `diagonal' subspace $\{(\omega,\ldots,\omega)\in X^n:\, \omega\in X\}$.
If the points $\omega_1,\ldots,\omega_n$ do not all coincide, then the quantity computed in accordance with formula \eqref{d2} will remain strictly positive if the infimum there is taken not over the whole space, but over any its subset.
This can be a way of defining `localized distances'.

Some properties of the quantities \eqref{d1} and \eqref{d2} are collected in the next proposition.

\begin{proposition}\label{P3.3}
Suppose $\omega_1,\ldots,\omega_n$ $(n\ge2)$ are points in a metric space $X$.
\begin{enumerate}
\item
$d_2(\omega_1,\ldots,\omega_n)\le d_1(\omega_1,\ldots,\omega_n)\le2 d_2(\omega_1,\ldots,\omega_n)$.
\item
If $X$ is a normed vector space with the distance induced by the norm and $n=2$, then
$d_1(\omega_1,\omega_2) =2d_2(\omega_1,\omega_2)$, i.e. the last inequality in {\rm(i)} holds as equality.
\cnta
\end{enumerate}
Suppose additionally that $X$ is a normed vector space.
\begin{enumerate}
\cntb
\item
$d_2(\omega_1,\ldots,\omega_n)\le d_3(\omega_1,\ldots,\omega_n)\le2 d_2(\omega_1,\ldots,\omega_n)$.
\item
If $n=2$, then
$d_2(\omega_1,\omega_2) =d_3(\omega_1,\omega_2)$, i.e. the first inequality in {\rm(iii)} holds as equality.
\end{enumerate}
\end{proposition}

\begin{Proof}
\begin{enumerate}
\item
The first inequality follows immediately from the definitions:
\begin{align*}
d_2(\omega_1,\ldots,\omega_n) &=\inf_{x\in X}\max_{1\le i\le n} d(\omega_i,x)
\le \max_{1\le i\le n} d(\omega_i,\omega_n)
\\
&=\max_{1\le i\le n-1} d(\omega_i,\omega_n) =d_1(\omega_1,\ldots,\omega_n).
\end{align*}
To prove the second inequality, first fix an $x\in X$.
\begin{align*}
d_1(\omega_1,\ldots,\omega_n) &=\max_{1\le i\le n-1} d(\omega_i,\omega_n)
\\
&\le\max_{1\le i\le n-1} (d(\omega_i,x)+d(\omega_n,x))
\le2\max_{1\le i\le n} d(\omega_i,x).
\end{align*}
Taking the infimum over $x\in X$ in the \RHS\ of the above inequality, we arrive at the second inequality in (i).
\item
Let $X$ be a normed vector space with the distance induced by the norm and $n=2$. Then
\begin{align*}
d_1(\omega_1,\omega_2)=\norm{\omega_1-\omega_2} &=2\max\left\{\norm{\frac{\omega_1+\omega_2}{2}-\omega_1}, \norm{\frac{\omega_1+\omega_2}{2}-\omega_2}\right\}
\\
&\ge2\inf_{x\in X} \max\{\norm{x-\omega_1},\norm{x-\omega_2}\}=2 d_2(\omega_1,\omega_2).
\end{align*}
Combining this with the second inequality in (i) proves (ii).
\item
Suppose that $X$ is a normed vector space.
The first inequality follows immediately from the definitions.
To prove the second inequality, first fix an $x\in X$.
Then
\begin{align*}
d_3(\omega_1,\ldots,\omega_n) &=\max_{1\le i\le n} \norm{\omega_i-\frac{1}{n}\sum_{j=1}^n\omega_j}
\le\max_{1\le i\le n}\norm{\omega_i-x} +\norm{\frac{1}{n}\sum_{j=1}^n\omega_j-x}
\\
&\le\max_{1\le i\le n} \norm{\omega_i-x} +\frac{1}{n}\sum_{j=1}^n\norm{\omega_j-x} \le2\max_{1\le i\le n} \norm{\omega_i-x}.
\end{align*}
Taking the infimum over $x\in X$ in the \RHS\ of the second inequality, we arrive at the last inequality in (iii).
\item
Let $n=2$.
Using definitions \eqref{d3} and \eqref{d1}, and the equality in (ii), we obtain:
\begin{align*}
d_3(\omega_1,\omega_2) &=\max\left\{\norm{\omega_1-\frac{\omega_1+\omega_2}{2}}, \norm{\omega_2-\frac{\omega_1+\omega_2}{2}}\right\}
\\
&=\frac{1}{2} \norm{\omega_1-\omega_2}=\frac{1}{2} d_1(\omega_1,\omega_2)=d_2(\omega_1,\omega_2).
\end{align*}
\qed\end{enumerate}
\end{Proof}

It follows from part (ii) of Proposition~\ref{P3.3} that quantity \eqref{d2} does not reduce to the conventional distance when $n=2$: in the setting of a normed vector space it equals $\frac{1}{2}d(\omega_1,\omega_2)$.
In part (ii) of Proposition~\ref{P3.3}, the infimum in formula \eqref{d2} is computed explicitly.
Unfortunately, when $n>2$, this seems impossible in general even in the setting of a normed vector space; see the discussion in \cite[Section~6]{BuiLinRos}.

\begin{Remark}
In the setting of a normed vector space, if $n=2$ and $\omega_1\ne\omega_2$, then, in view of Proposition~\ref{P3.3}(ii) and (iv), the first inequality in Proposition~\ref{P3.3}(i) and the second inequality in Proposition~\ref{P3.3}(iii) are strict.
If $n>2$, then all the inequalities in Proposition~\ref{P3.3}(i) and (iii) can be strict.
This fact is illustrated by Example~\ref{E3.1}.
\end{Remark}

\subsection{Distances Between $n$ Sets}

Now, we employ the distances between collections of points discussed in the previous subsection to quantify `closeness' of collections of sets.
Given $n\ge2$ subsets $\Omega_1,\ldots,\Omega_n$ of a metric space $X$ and an $n$-point distance $d$, the distance between $\Omega_1,\ldots,\Omega_n$ is defined in the usual way:
\begin{equation}\label{d}
d(\Omega_1,\ldots,\Omega_n) :=\inf_{\omega_1\in\Omega_1,\ldots,\omega_n\in\Omega_n} d(\omega_1,\ldots,\omega_n).
\end{equation}

Applying construction \eqref{d} to the $n$-point distances \eqref{d1}, \eqref{d2} and \eqref{d3}, we obtain the following definitions of particular distances between $n$ sets, respectively:
\begin{align}\label{dd1}
d_1(\Omega_1,\ldots,\Omega_n) &:=\inf_{\omega_1\in\Omega_1,\ldots,\omega_n\in\Omega_n} \max_{1\le i\le n-1} d(\omega_i,\omega_n),
\\\label{dd2}
d_2(\Omega_1,\ldots,\Omega_n) &:=\inf_{\omega_1\in\Omega_1,\ldots,\omega_n\in\Omega_n, x\in X}\max_{1\le i\le n} d(\omega_i,x),
\\\label{dd3}
d_3(\Omega_1,\ldots,\Omega_n) &:=\inf_{\omega_1\in\Omega_1,\ldots,\omega_n\in\Omega_n} \max_{1\le i\le n} \norm{\omega_i-\frac{1}{n}\sum_{j=1}^n\omega_j}.
\end{align}
When $n=2$, definition \eqref{dd1} reduces to the conventional distance \eqref{d0}.
Definition \eqref{dd3} is meaningful in the setting of a normed vector space only.

Proposition \ref{P3.3} leads to similar relations for the distances between $n$ sets.

\begin{proposition}\label{P3.9}
Suppose $\Omega_1,\ldots,\Omega_n$ are subsets of a metric space $X$.
\begin{enumerate}
\item
$d_2(\Omega_1,\ldots,\Omega_n)\le d_1(\Omega_1,\ldots,\Omega_n)\le2 d_2(\Omega_1,\ldots,\Omega_n)$.
\item
If $n=2$, then
$d_1(\Omega_1,\Omega_2)=2d_2(\Omega_1,\Omega_2)$.
\item
$d_1(\Omega_1,\ldots,\Omega_n)>0
\quad\Longleftrightarrow\quad
d_2(\Omega_1,\ldots,\Omega_n)>0$.
\cnta
\end{enumerate}
Suppose additionally that $X$ is a normed vector space.
\begin{enumerate}
\cntb
\item
$d_2(\Omega_1,\ldots,\Omega_n)\le d_3(\Omega_1,\ldots,\Omega_n)\le2 d_2(\Omega_1,\ldots,\Omega_n)$.
\item
If $n=2$, then
$d_2(\Omega_1,\Omega_2) =d_3(\Omega_1,\Omega_2)$.
\item
$d_1(\Omega_1,\ldots,\Omega_n)>0
\quad\Longleftrightarrow\quad
d_2(\Omega_1,\ldots,\Omega_n)>0
\quad\Longleftrightarrow\quad
d_3(\Omega_1,\ldots,\Omega_n)>0$.
\end{enumerate}
\end{proposition}

\begin{Proof}
Conditions (i), (ii) and (iv), (v) are consequences of the corresponding conditions in Proposition \ref{P3.3}.
Condition (iii) is a consequence of condition (i).
Condition (vi) is a consequence of conditions (iii) and (iv).
\qed\end{Proof}

Another distance can be of interest and is going to be used in the sequel.
Given a subset $\Omega_{n+1}$ of $X$, define
\begin{align*}
d_{\Omega_{n+1}}(\Omega_1,\ldots,\Omega_n) &:=\inf_{\omega_1\in\Omega_1,\ldots,\omega_n\in\Omega_n, x\in\Omega_{n+1}}\max_{1\le i\le n} d(\omega_i,x).
\end{align*}
Now observe that $d_{\Omega_{n+1}}$ as well as the symmetric distance $d_2$ \eqref{dd2} are particular cases of the asymmetric distance $d_1$ \eqref{dd1} applied to $n+1$ sets:
\begin{align}\notag
&d_{\Omega_{n+1}}(\Omega_1,\ldots,\Omega_n)=
d_1(\Omega_1,\ldots,\Omega_{n},\Omega_{n+1})
\qdtx{and}\label{3.11}
\\
&d_2(\Omega_1,\ldots,\Omega_n)=
d_1(\Omega_1,\ldots,\Omega_n,X).
\end{align}
This observation makes the asymmetric distance $d_1$ a rather general quantitative measure of closeness of collections of sets.
Apart from the most straightforward case $\Omega_{n+1}:=X$ used in the above example, another useful particular case is given by $\Omega_{n+1}:=B_\rho(\bx)$ where $\bx\in X$ is a fixed point (related to the sets $\Omega_1,\ldots,\Omega_n)$ and $\rho>0$.
This allows one to examine closeness of sets in a neighbourhood of the given point.

\begin{Remark}\label{R3.10}
\begin{enumerate}
\item
As it was observed earlier, the distance $d_1$ \eqref{dd1} depends in general on the order of the sets.
However, thanks to Proposition~\ref{P3.9}(iii), if it is strictly positive for some permutation of the sets, it remains strictly positive for any other permutation.
\item
The maximum operation in all the above definitions of distances between $n$ points and $n$ sets corresponds to the maximum norm in either $\R^{n-1}$ or $\R^{n}$.
It can be replaced in all these definitions and subsequent statements by any other finite dimensional norm producing different but in a sense equivalent `distances'.
The $p$-norm version of the quantity \eqref{dd1} was considered in \cite{ZheNg11} under the name \emph{($p$-weighted) nonintersect index}.
In the current paper, for the sake of simplicity of presentation only the maximum norm is considered.
\end{enumerate}
\end{Remark}

The next proposition and its corollary characterize a set of points in the given collection of sets, which are almost closest (up to $\eps)$ points of these sets with respect to the chosen $n$-point distance.
(When $n=2$, such points play a key role in the geometric versions of the Ekeland variational principle considered in Sect.~\ref{S2}.)
It also introduces a two-step procedure, which is going to be used in the sequel.
Given a collection of sets $\Omega_1,\ldots,\Omega_n$ with empty intersection and a collection of points $\omega_i\in\Omega_i$ $(i=1,\ldots,n)$, we
\begin{enumerate}
\item[1)]
consider another collection of sets $\Omega_1-\omega_1,\ldots,\Omega_n-\omega_n$, whose intersection is obviously nonempty, and
\item[2)]
construct `small' (up to $\eps)$ translation vectors $a_1,\ldots,a_n$ such that the translated sets $\Omega_1-\omega_1-a_1,\ldots,\Omega_n-\omega_n-a_n$ have empty intersection again.
\end{enumerate}
Thus, the proposition translates $\eps$-closeness of a collection of points, which served as the key assumption when proving unified separation theorems in \cite{ZheNg11}, into the language of $\eps$-translations of the sets employed in Definition~\ref{D2.5} of extre\-mality/stationarity properties, the corresponding Proposition~\ref{P2.7} and parts (ii) and (iii) of Proposition~\ref{P2.11}.

For simplicity and in view of the observed above universality of the distance $d_1$ \eqref{dd1}, we first consider this distance.

\begin{proposition}\label{P3.11}
Suppose $\Omega_1,\ldots,\Omega_n$ are subsets of a normed vector space $X$, $\cap_{i=1}^n\Omega_i=\emptyset$, $\omega_i\in\Omega_i$ $(i=1,\ldots,n)$, and $\eps>0$.
If
\begin{equation}\label{P3.11-1}
d_1(\omega_1,\ldots,\omega_n)< d_1(\Omega_1,\ldots,\Omega_n)+\eps,
\end{equation}
then $M:=d_1(\omega_1,\ldots,\omega_n)>0$ and condition \eqref{P8} is satisfied,
where
$a_i:=\frac{\eps'}{M}(\omega_n-\omega_i)$ $(i=1,\ldots,n-1)$
and
$\eps'$ is either any number in $[0,\eps[$ satisfying $M-d_1(\Omega_1,\ldots,\Omega_n)<\eps'\le M$ if $d_1(\Omega_1,\ldots,\Omega_n)>0$, or $\eps'=M$ if $d_1(\Omega_1,\ldots,\Omega_n)=0$.
\sloppy
\end{proposition}

\begin{Proof}
In view of $\cap_{i=1}^n\Omega_i=\emptyset$, we have $M>0$.
Let condition \eqref{P3.11-1} be satisfied.

Suppose first that $d_1(\Omega_1,\ldots,\Omega_n)>0$.
Choose a positive number $\eps'<\eps$ such that $M-d_1(\Omega_1,\ldots,\Omega_n)<\eps'\le M$.
Then, $\max_{1\le i\le n-1}\norm{a_i}=\eps'<\eps$.
Suppose that the first condition in \eqref{P8} does not hold.
Then, there exists a point ${x\in\cap_{i=1}^{n-1}(\Omega_i-\omega_i-a_i) \cap(\Omega_n-\omega_n)}$, and consequently, $\hat\omega_i:=\omega_i+a_i+x\in\Omega_i$ $(i=1,\ldots,n-1)$ and $\hat\omega_n:=\omega_n+x\in\Omega_n$.
Thus,
$$\norm{\hat\omega_i-\hat\omega_n} =\norm{\omega_i-\omega_n -\frac{\eps'}{M}(\omega_i-\omega_n)} =\left(1-\frac{\eps'}{M}\right)\norm{\omega_i-\omega_n} \quad (i=1,\ldots,n-1),$$
and consequently,
$$d_1(\Omega_1,\ldots,\Omega_n)\le d_1(\hat\omega_1,\ldots,\hat\omega_n) =\left(1-\frac{\eps'}{M}\right) d_1(\omega_1,\ldots,\omega_n)=M-\eps'.$$
This contradicts the choice of $\eps'$.
Hence, conditions \eqref{P8} hold true.

In the case $d_1(\Omega_1,\ldots,\Omega_n)=0$,
set $a_i:=\omega_n-\omega_i$ $(i=1,\ldots,n-1)$.
Then, by \eqref{P3.11-1}, $\max_{1\le i\le n-1} \norm{a_i}=d_1(\omega_1,\ldots,\omega_n)<\eps$ and
\begin{equation*}
\bigcap_{i=1}^{n-1}(\Omega_i-\omega_i-a_i) \cap(\Omega_n-\omega_n) =\bigcap_{i=1}^n(\Omega_i-\omega_n)
=\bigcap_{i=1}^n\Omega_i-\omega_n
=\emptyset,
\end{equation*}
i.e. conditions \eqref{P8} hold true.
\qed\end{Proof}

Applying Proposition~\ref{P3.11} to the collection of $n+1$ sets $\Omega_1,\ldots,\Omega_n,X$, we arrive at the following statement in terms of the distance $d_2$ \eqref{dd2}.
\begin{corollary}
Suppose $\Omega_1,\ldots,\Omega_n$ are subsets of a normed vector space $X$, $\cap_{i=1}^n\Omega_i=\emptyset$, $\omega_i\in\Omega_i$ $(i=1,\ldots,n)$, $x\in X$ and $\eps>0$.
If
\begin{equation}\label{3.14}
\max_{1\le i\le n}d(\omega_i,x)< d_2(\Omega_1,\ldots,\Omega_n)+\eps,
\end{equation}
then $M:=\max_{1\le i\le n}\norm{\omega_i-x}>0$ and conditions \eqref{P1} hold true with $\Omega'_i:=\Omega_i-\omega_i$ in place of $\Omega_i$ and
$a_i:=\frac{\eps'}{M}(x-\omega_i)$ $(i=1,\ldots,n)$, where
$\eps'$ is either any number in $]0,\eps[$ satisfying $M-d_2(\Omega_1,\ldots,\Omega_n)<\eps'\le M$ if $d_2(\Omega_1,\ldots,\Omega_n)>0$, or $\eps'=M$ if $d_2(\Omega_1,\ldots,\Omega_n)=0$.
\end{corollary}

\section{Localizations and Translations of the Sets}\label{S5}

In this section, we continue studying mutual arrangement of $n$ sets in space started in Sect.~\ref{S4}.

The next lemma extends the Asymmetric Geometric Ekeland Variational Principle (AGEVP) from Sect.~\ref{S2} to the case of $n\ge2$ sets, and as such it is also equivalent to the Ekeland Variational Principle.
Just like in AGEVP, balls of two different radii are used in the concluding part of the lemma, one of them employed in the localizations of the first $n-1$ sets and the other one for the remaining single set.
The latter set is going to play a special role in the subsequent analysis.

The lemma translates $\eps$-closeness of a given collection of points into $\xi\eps$-closeness (with an additional parameter $\xi)$ of another collection of points with respect to a collection of certain localizations of the sets (depending on $\xi)$.

\begin{lemma}\label{L4.1}
Suppose $\Omega_1,\ldots,\Omega_n$ are closed subsets of a complete metric space $X$,
$\omega_{i}\in\Omega_i$ $(i=1,\ldots,n)$, and $\eps>0$.
If condition \eqref{P3.11-1} is satisfied,
then, for all numbers $\la,\rho>0$, there exist $\hat\omega_{i}\in\Omega_i\cap B_{\la}(\omega_{i})$ $(i=1,\ldots,n-1)$ and $\hat\omega_{n}\in\Omega_n\cap B_{\rho}(\omega_{n})$ such that
\begin{enumerate}
\item
$d_1(\hat\omega_1,\ldots,\hat\omega_n)\le d_1(\omega_1,\ldots,\omega_n)$;
\item
$d_1(\Omega_1\cap B_{\xi\la}(\hat\omega_1),\ldots, \Omega_{n-1}\cap B_{\xi\la}(\hat\omega_{n-1}), \Omega_n\cap B_{\xi\rho}(\hat\omega_{n})) +{\xi\eps}>d_1(\hat\omega_1,\ldots,\hat\omega_n)$ for all $\xi>0$.
\end{enumerate}
\end{lemma}

\begin{Proof}
The product space $X^{n-1}$ considered with the maximum metric is complete.
Set $A:=\Omega_1\times\ldots\times\Omega_{n-1}$, and $B:=\{(x,\ldots,x):\, x\in\Omega_n\}\subset X^{n-1}$, $a:=(\omega_1,\ldots,\omega_{n-1})\in A$, $b:=(\omega_n,\ldots,\omega_n)\in B$.
Then $d_1(\Omega_1,\ldots,\Omega_n)=d(A,B)$ and $d_1(\omega_1,\ldots,\omega_n)=d(a,b)$.

Applying AGEVP, we find points $\hat a =(\hat\omega_1,\ldots,\hat\omega_{n-1})
\allowbreak
\in A\cap B_{\la}(a)$ and $\hat b =(\hat\omega_n,\ldots,\hat\omega_n)\in B\cap B_{\rho}(b)$ satisfying conditions (i) and (ii) in AGEVP.
Recalling that the maximum metric is used in $X^{n-1}$, it follows that $\hat\omega_{i}\in\Omega_i\cap B_{\la}(\omega_{i})$ $(i=1,\ldots,n-1)$, $\hat\omega_{n}\in B\cap B_{\rho}(\omega_{n})$ and conditions (i) and (ii) above are satisfied.
\qed\end{Proof}

In view of the equalities in \eqref{3.11}, the following corollary involving the symmetric distance $d_2$ is immediate.

\begin{corollary}
Suppose $\Omega_1,\ldots,\Omega_n$ are closed subsets of a complete metric space $X$,
$\omega_{i}\in\Omega_i$ $(i=1,\ldots,n)$, $x\in X$ and $\eps>0$.
If condition \eqref{3.14} is satisfied,
then, for all $\la,\rho>0$, there exist $\hat\omega_{i}\in\Omega_i\cap B_{\la}(\omega_{i})$ $(i=1,\ldots,n)$ and $\hat x\in B_{\rho}(x)$ such that
\begin{enumerate}
\item
$\max_{1\le i\le n} d(\hat\omega_i,\hat x)\le \max_{1\le i\le n} d(\omega_i,x)$;
\item
$d_{B_{\xi\rho}(\hat x)}(\Omega_1\cap B_{\xi\la}(\hat\omega_1),\ldots, \Omega_n\cap B_{\xi\la}(\hat\omega_{n})) +{\xi\eps}>\max_{1\le i\le n} d(\hat\omega_i,\hat x)$ for all $\xi>0$.
\end{enumerate}
\end{corollary}

The next proposition is a consequence of Lemma~\ref{L4.1} in the Banach space setting.
It characterizes $\eps$-closest points of a collection of sets with empty intersection and involves localizations of the sets and small (up to $\xi\eps)$ translations of the first $n-1$ localizations.
The proposition transforms the nonlocal non-intersection condition $\cap_{i=1}^n\Omega_i=\emptyset$ into a non-intersection condition of translated localizations of the sets.
It also exposes the special role played by the last set in the list.

\begin{proposition}\label{P4.2}
Suppose $\Omega_1,\ldots,\Omega_n$ are closed subsets of a Banach space $X$,
$\omega_{i}\in\Omega_i$ $(i=1,\ldots,n)$, and $\eps>0$.
If $\cap_{i=1}^n\Omega_i=\emptyset$ and condition \eqref{P3.11-1} is satisfied,
then, for any $\la,\rho>0$, there are points $\hat\omega_{i}\in\Omega_i\cap B_{\la}(\omega_i)$ $(i=1,\ldots,n-1)$, $\hat\omega_n\in\Omega_n\cap B_{\rho}(\omega_n)$ such that, for any $\xi>0$, there exist vectors $a_i\in X$ $(i=1,\ldots,n-1)$ satisfying
\begin{equation}\label{P4.2-1}
\bigcap_{i=1}^{n-1} \Big(\big((\Omega_i-\hat\omega_{i}) \cap(\xi\la)\B\big)-a_i\Big) \cap (\Omega_n-\hat\omega_n)\cap(\xi\rho)\B
=\emptyset
\qdtx{and}
\max_{1\le i\le n-1}\norm{a_i}<\xi\eps.
\end{equation}
\end{proposition}

\begin{Proof}
Applying Lemma~\ref{L4.1}, we find points $\hat\omega_{i}\in\Omega_i\cap B_{\la}(\omega_{i})$ $(i=1,\ldots,n-1)$ and $\hat\omega_{n}\in\Omega_n\cap B_{\rho}(\omega_{n})$
satisfying condition (ii) in that lemma.
Given any number $\xi>0$, we can now apply Proposition~\ref{P3.11} with the function $d_1$ and sets $\Omega_1\cap B_{\xi\la}(\hat\omega_1)$,\ldots, $\Omega_{n-1}\cap B_{\xi\la}(\hat\omega_{n-1})$, $\Omega_n\cap B_{\xi\rho}(\hat\omega_{n})$, points $\hat\omega_{i}$ $(i=1,\ldots,n)$ and number $\xi\eps$ in place of sets $\Omega_1$,\ldots, $\Omega_n$, points $\omega_{i}$ $(i=1,\ldots,n)$ and number $\eps$, respectively, to find vectors $a_i\in X$ $(i=1,\ldots,n-1)$ such that condition \eqref{P4.2-1} is satisfied.
\qed\end{Proof}

\begin{Remark}
\begin{enumerate}
\item
Under the conditions of Proposition~\ref{P4.2}, $0\in\cap_{i=1}^{n}(\Omega_i-\omega_{i})$, $0\in\cap_{i=1}^{n}(\Omega_i-\hat\omega_{i})$, and the expressions involved in the first condition in \eqref{P4.2-1} correspond to localizations of the sets $\Omega_i-\hat\omega_{i}$ $(i=1,\ldots,n)$ near 0, or equivalently, localizations of the original sets $\Omega_i$ near $\hat\omega_i$ $(i=1,\ldots,n)$.
\item
There are certain similarities between Propositions~\ref{P3.11} and \ref{P4.2} in terms of both assumptions and conclusions.
There are also important differences.
The assumptions of Propositions~\ref{P3.11} are weaker: the space is not assumed to be complete and the sets are not assumed to be closed.
The concluding non-intersection condition in \eqref{P8} (Propositions~\ref{P3.11}) is formulated for the given sets and in terms of the given collection of points, while Proposition~\ref{P4.2} establishes existence of another collection of points, and the corresponding condition in \eqref{P4.2-1} is formulated for localizations of the sets near these points.
The translations of the sets are constructed in Propositions~\ref{P3.11} explicitly and are entirely determined by the given collection of points.
At the same time, in Proposition~\ref{P4.2} the size of the translations and localizations of the sets as well as the distance of the new points from the given ones are controlled by additional parameters.
These parameters, which appear in Proposition~\ref{P4.2}, represent the major feature of this statement compared to Propositions~\ref{P3.11}.
They provide an additional degree (degrees) of freedom for the applications of the result.
\end{enumerate}
\end{Remark}

Applying Proposition~\ref{P4.2} to the collection of $n+1$ sets $\Omega_1,\ldots,\Omega_n,X$, we arrive at the following statement.

\begin{corollary}
Suppose $\Omega_1,\ldots,\Omega_n$ are closed subsets of a Banach space $X$,
$\omega_{i}\in\Omega_i$ $(i=1,\ldots,n)$, $x\in X$ and $\eps>0$.
If $\cap_{i=1}^n\Omega_i=\emptyset$ and condition \eqref{3.14} is satisfied,
then, for any numbers $\la,\rho>0$, there are points $\hat\omega_{i}\in\Omega_i\cap B_{\la}(\omega_i)$ $(i=1,\ldots,n)$ such that, for any $\xi>0$, there exist vectors $a_i\in X$ $(i=1,\ldots,n)$ satisfying
\begin{equation}\label{C4.5-1}
\bigcap_{i=1}^{n} \Big(\big((\Omega_i-\hat\omega_{i}) \cap(\xi\la)\B\big)-a_i\Big) \cap(\xi\rho)\B
=\emptyset
\qdtx{and}
\max_{1\le i\le n}\norm{a_i}<\xi\eps.
\end{equation}
\end{corollary}

Observe that conditions \eqref{P4.2-1} and \eqref{C4.5-1} are exactly conditions \eqref{P5} and \eqref{P2}, respectively, applied to the localizations of the sets $\Omega_1-\hat\omega_1,\ldots,\Omega_n-\hat\omega_n$ with $\bx:=0$, and $\xi\eps$ and $\xi\rho$ in place of $\eps$ and $\rho$, respectively.

Proposition~\ref{P4.2} assumes that the sets have empty intersection.
Next, we demonstrate that it can be also applied to collections of sets having a common point.
Specifically, we consider the special case \eqref{P2} in Definition~\ref{D2.5}(ii) of local extremality, where the last set in the list (of $n+1$ sets) is a ball centred at this point.
Since the radius of the ball is allowed to be infinite, the global setting is covered too.

\begin{proposition}\label{P4.3}
Suppose $\Omega_1,\ldots,\Omega_n$ are closed subsets of a Banach space $X$, $\bx\in\cap_{i=1}^n\Omega_i$, $\eps>0$ and $\rho\in]0,\infty]$.
If conditions \eqref{P2} are satisfied
for some vectors $a_i\in{X}$ $(i=1,\ldots,n)$, then, for any $\la>0$, there are points $\omega_i\in\Omega_i\cap B_{\la}(\bx)$ $(i=1,\ldots,n)$ and a number $\delta\in]0,1[$ such that, for all $\xi\in ]0,\delta[$,
\begin{equation}\label{P4.3-2}
\bigcap_{i=1}^{n} \Big(\big((\Omega_i-\omega_{i}) \cap(\xi\la)\B\big)-a_i'\Big) \cap(\xi\rho)\B
=\emptyset
\qdtx{and}
\max_{1\le i\le n}\norm{a_i'}<\xi\eps
\end{equation}
for some vectors $a_i'\in{X}$ $(i=1,\ldots,n)$.

Moreover, if $\la\ge\rho+\eps$, then
\begin{equation}\label{P4.3-3}
\bigcap_{i=1}^n(\Omega_i-\omega_i-a'_i)\cap (\xi\rho)\B=\emptyset;
\end{equation}
if $\la+\eps\le\rho$, then
\begin{equation}\label{P4.3-4}
\bigcap_{i=1}^{n} \Big(\big((\Omega_i-\omega_{i}) \cap(\xi\la)\B\big)-a_i'\Big)=\emptyset.
\end{equation}
\end{proposition}

\begin{Proof}
Let $\rho'\in]0,\rho[$.
By \eqref{P2}, $\cap_{i=1}^{n+1}\Omega_i'=\emptyset$, where $\Omega_i':=\Omega_i-a_i$ $(i=1,\ldots,n)$ and $\Omega_{n+1}':=\overline{B}_{\rho'}(\bx)$, and $$d_1(\bx-a_1,\ldots,\bx-a_n,\bx)=\max_{1\le i\le n} \norm{a_i}<\eps.$$
Applying Proposition~\ref{P4.2}, we can find points $\omega_i\in\Omega_i\cap B_\la(\bx)$ $(i=1,\ldots,n)$ and $x\in X$ with $\norm{x}<\rho'$ such that, for any $\xi>0$, there exist vectors $a_i'\in X$ $(i=1,\ldots,n)$ satisfying
\begin{equation*}
\bigcap_{i=1}^{n} \Big(\big((\Omega_i-\omega_{i}) \cap(\xi\la)\B\big)-a_i'\Big)\cap\overline{B}_{\rho'}(x)\cap(\xi\rho)\B=\emptyset
\qdtx{and}
\max_{1\le i\le n}\norm{a_i'}<\xi\eps.
\end{equation*}
Set $\de:=(\rho'-\norm{x})/\rho$.
Then, $\delta\in]0,1[$ and, for any $\xi\in]0,\delta[$, we have $(\xi\rho)\B\subset\overline{B}_{\rho'}(x)$, and consequently, \eqref{P4.3-2} holds true.
If $\la\ge\rho+\eps$, then $(\xi\rho)\B\subset \left(\xi\la\right)\B-a_i'$ for all $i=1,\ldots,n$, and consequently, \eqref{P4.3-2} implies \eqref{P4.3-3}.
Similarly, if $\la+\eps\le\rho$, then $\left(\xi\la\right)\B-a_i'\subset(\xi\rho)\B$ for all $i=1,\ldots,n$, and consequently, \eqref{P4.3-2} implies \eqref{P4.3-4}.
\qed\end{Proof}

\if{
\begin{corollary}\label{C4.6}
	Suppose $\Omega_1,\ldots,\Omega_n$ are closed subsets of a Banach space $X$, $\bx\in\cap_{i=1}^n\Omega_i$, $\eps>0$ and $\rho\in]0,\infty]$.
	If conditions \eqref{P2} are satisfied
	for some vectors $a_i\in{X}$ $(i=1,\ldots,n)$, then, there are points $\omega_i\in\Omega_i\cap B_{\rho+\eps}(\bx)$ $(i=1,\ldots,n)$
	such that the sets $\{\Omega_1-\omega_1,\ldots,\Omega_n-\omega_n\}$ is ($\eps/\rho$)--approximately stationary at $0$.
\end{corollary}

\begin{Proof}
	Apply Proposition~\ref{P4.3} with $\la =\rho+\eps$, we find $\omega_i\in\Omega_i\cap B_{\rho+\eps}(\bx)$ $(i=1,\ldots,n)$ that \eqref{P4.3-2} (or \eqref{P4.3-3}) holds. And by Definition~\ref{D2.17}, we have the assertion.
	\qed
\end{Proof}
}\fi

\begin{Remark}
Conditions \eqref{P2} in Proposition~\ref{P4.3} are formulated for a fixed point $\bx\in\cap_{i=1}^n\Omega_i$ and fixed $\eps>0$ and $\rho\in]0,\infty]$.
It presumes a certain balance between the values of $\eps$ and $\rho$: the larger the value of $\rho$ is, the larger value of $\eps$ is needed to ensure the existence of  vectors $a'_i\in X$ $(i=1,\ldots,n)$ satisfying \eqref{P2}.
In contrast, condition \eqref{P4.3-2} involves two additional parameters: an arbitrary $\la>0$ and a sufficiently small $\xi\in]0,1[$.
Fixed $\eps$ and $\rho$ are replaced by $\xi\eps$ and $\xi\rho$, respectively, preserving their ratio, while the sets are replaced by their localizations controlled by $\xi\la$.
This advancement comes at a price: instead of a single common fixed point $\bx$, we now have to deal with a collection of individual points $\omega_i\in\Omega_i$ $(i=1,\ldots,n)$, whose distance from $\bx$ is controlled by $\la$.
Now the balance between this distance and the size of the localizations of the sets becomes important: choosing a smaller $\la$ ensures that the individual points $\omega_i$ are closer to $\bx$ while at the same time reducing the size of the localizations of the sets and, thus, weakening condition \eqref{P4.3-2}.
\end{Remark}

Given a collection of sets $\Omega_1,\ldots,\Omega_n$ $(n\ge2)$, a point $\bx\in\cap_{i=1}^n\Omega_i$ and a $\rho\in]0,\infty]$, define (cf. \cite{Kru04,Kru05}):
\begin{gather*}
\theta_\rho[\Omega_1,\ldots,\Omega_n](\bx)= \sup\set{r>0:\,
\bigcap_{i=1}^n(\Omega_i-a_i)\cap
B_\rho(\bx)\ne\emptyset
\qdtx{for all} a_i\in r\B}.
\end{gather*}
This nonnegative quantity tells us how far the sets can be pushed apart until their intersection becomes empty with respect to the fixed $\rho$-\nbh\ of $\bx$.

\begin{corollary}
Suppose $\Omega_1,\ldots,\Omega_n$ are closed subsets of a Banach space $X$ and $\bx\in\cap_{i=1}^n\Omega_i$.
If $\rho>0$ and $\eps>\theta_\rho[\Omega_1,\ldots,\Omega_n](\bx)$, then
\begin{equation*}
\inf_{\omega_i\in\Omega_i\cap B_{\rho+\eps}(\bx)\atop (i=1,\ldots,n)}\limsup_{\al\downarrow0} \al\iv \theta_\al[\Omega_1-\omega_1,\ldots,\Omega_n-\omega_n](0) \le \rho\iv\theta_\rho[\Omega_1,\ldots,\Omega_n](\bx).
\end{equation*}
Moreover, if $\bx\in\bd\cap_{i=1}^n\Omega_i$, then
\begin{equation*}
\liminf_{\omega_i\to\bx,\,\omega_i\in\Omega_i\atop (i=1,\ldots,n)}\limsup_{\al\downarrow0} \al\iv \theta_\al[\Omega_1-\omega_1,\ldots,\Omega_n-\omega_n](0) \le \liminf_{\rho\downarrow0} \rho\iv\theta_\rho[\Omega_1,\ldots,\Omega_n](\bx).
\end{equation*}
\end{corollary}

\begin{Proof}
The first assertion is a direct consequence of Proposition~\ref{P4.3}.
The second assertion is a consequence of the first one since $\bx\in\bd\cap_{i=1}^n\Omega_i$ implies that $\theta_\rho[\Omega_1,\ldots,\Omega_n](\bx)\to0$ as $\rho\downarrow0$; cf. \cite[Proposition~3]{Kru05}.
\qed\end{Proof}

The next proposition presents a metric counterpart of the conditions \eqref{P2}.
It contains the key ingredients of the metric criteria of approximate stationarity and transversality in Theorems~\ref{T2.8} and \ref{T2.13}, respectively.

\begin{proposition}\label{P4.6}
Suppose $\Omega_1,\ldots,\Omega_n$ are subsets of a normed vector space $X$, $\bx\in\cap_{i=1}^n\Omega_i$, $a_i\in{X}$ $(i=1,\ldots,n)$, $\eps>0$, $\rho>0$ and $\al:=\frac{\eps}{\rho}$.
Then
\begin{enumerate}
\item
conditions \eqref{P2} imply conditions \eqref{P9};
\item
if conditions \eqref{P9} are satisfied, then there exist a number $\rho'\in]0,\rho[$, points $\omega_i\in\Omega_i\cap B_{\eps}(\bx)$ and vectors $a_i'\in X$ $(i=1,\ldots,n)$ such that conditions \eqref{P7} are satisfied with $\rho'$ and $a_i'$ in place of $\rho$ and $a_i$.
\end{enumerate}
\end{proposition}

\begin{Proof}
\begin{enumerate}
\item
If conditions \eqref{P2} are satisfied,
then $\max_{1\le i\le n}\norm{a_i}<\eps$ and
$$
\al d\left(\bx,\bigcap_{i=1}^n(\Omega_i-a_i)\right) >\al\rho= \eps>\max_{1\le i\le n}\norm{a_i}\ge \max_{1\le i \le n} d\left(\bx,\Omega_i-a_i\right).
$$
\item
Let conditions \eqref{P9} be satisfied.
Then,
\begin{equation*}
\max_{1\le i\le n} d\left(\bx,\Omega_i-a_i\right)\le\max_{1\le i\le n} \norm{a_i}<\eps=\al\rho,
\end{equation*}
and there exists a number $\rho'\in]0,\rho[$
such that
\begin{equation*}
\al d\left(\bx,\bigcap_{i=1}^n(\Omega_i-a_i)\right) >\al\rho'> \max_{1\le i\le n} d\left(\bx,\Omega_i-a_i\right).
\end{equation*}
It follows from the first inequality above that
$\cap_{i=1}^n(\Omega_i-a_i)\cap B_{\rho'}(\bx)=\emptyset$,
while due to the second inequality, there exist $\omega_i\in\Omega_i$ $(i=1,\ldots,n)$ such that
$
\max_{1\le i\le n}\norm{a_i'}<\al\rho',
$
where $a_i':=a_i+\bx-\omega_i$ $(i=1,\ldots,n)$.
\sloppy
\qed\end{enumerate}
\end{Proof}

\begin{corollary}\label{C4.9}
Suppose $\Omega_1,\ldots,\Omega_n$ are subsets of a normed vector space $X$, $\bx\in\cap_{i=1}^n\Omega_i$ and $\al>0$.
The following conditions are equivalent:
\begin{enumerate}
\item
the collection $\{\Omega_1,\ldots,\Omega_n\}$ is approximately $\al$-stationary at $\bx$;
\item
for any $\eps>0$, there exist vectors $a_i\in X$ $(i=1,\ldots,n)$ such that conditions \eqref{P9} are satisfied;
\item
for any $\eps>0$, there exist vectors $x\in B_\eps(\bx)$ and $a_i\in X$ $(i=1,\ldots,n)$ such that conditions \eqref{P10} are satisfied. \end{enumerate}
\end{corollary}

\begin{Proof}
\underline{(i) \folgt (ii)}.
Let condition (i) be satisfied and a number $\eps>0$ be given.
Set $\eps':=\eps/(\al+1)$.
Then, conditions \eqref{P7} hold with some $\rho\in]0,\eps'[$, $\omega_i\in\Omega_i\cap B_{\eps'}(\bx)$, and $a_i\in X$ $(i=1,\ldots,n)$.
By Proposition~\ref{P4.6}(i), applied to the sets $\Omega_i-\omega_i$ $(i=1,\ldots,n)$ having the common point 0, we have
$$
\al d\left(0,\bigcap_{i=1}^n(\Omega_i-\omega_i-a_i)\right) >\max_{1\le i\le n}d(0,\Omega_i-\omega_i-a_i).
$$
It follows that the first inequality in \eqref{P9} is satisfied with
$a_i':=a_i+\omega_i-\bx$ in place of $a_i$ $(i=1,\ldots,n)$, and
$\max_{1\le i\le n}\norm{a'_i}<\al\rho+\eps'<(\al+1)\eps'=\eps$.
Hence, condition (ii) is satisfied.

\underline{(ii) \folgt (i)}.
Let condition (ii) be satisfied and a number $\eps>0$ be given.
Set $\eps':=\min\{\al,1\}\eps$, and find vectors $a_i\in X$ $(i=1,\ldots,n)$ such that conditions \eqref{P9} are satisfied with $\eps'$ in place of $\eps$.
By Proposition~\ref{P4.6}(ii), there are $\rho\in]0,\eps'/\al[$, $\omega_i\in\Omega_i\cap B_{\eps'}(\bx)$, $a_i'\in X$ $(i=1,\ldots,n)$ such that conditions \eqref{P7} hold with $a_i'$ in place of $a_i$.
Since $\eps'\le\eps$ and $\eps'/\al\le\eps$, condition (i) is satisfied.

\underline{(ii) \folgt (iii)} is obvious.

\underline{(iii) \folgt (ii)}.
Let condition (iii) be satisfied and a number $\eps>0$ be given.
Set $\eps':=\eps/2$.
Then, conditions \eqref{P10} hold with some $x\in B_{\eps'}(\bx)$, $a_i\in X$ $(i=1,\ldots,n)$ and $\eps'$ in place of $\eps$.
Set $a_i':=a_i+x-\bx$ $(i=1,\ldots,n)$.
Then conditions \eqref{P9} hold true with $a_i'$ in place of $a_i$.
Hence, condition (ii) is satisfied.
\qed\end{Proof}

In view of Proposition~\ref{P2.18},
Corollary~\ref{C4.9} immediately yields Theorems~\ref{T2.8} and~\ref{T2.13}.

\section{Dual Characterizations}\label{S6}

This section presents a series of `generalized separation' statements, providing dual characterizations of certain typical `extremal' arrangements of collections of sets, discussed in the preceding sections, and traces the relationships between them.
These statements contain core arguments, which can be found in various existing versions of the (extended) extremal principle, as well as some new extensions.

The definition of the approximate stationarity (Definition~\ref{D2.5}(iv)) and its conventional dual characterizations in the extended extremal principle (Theorem~\ref{T2.9}(ii) and (iii)) are all formulated `for any $\eps>0$ there exist \ldots'.
At the same time, the proofs of the extremal principle and its extensions establish connections between the values of $\eps$ and other parameters involved in the assumptions and the conclusions.
These connections, usually hidden in the proofs, are of importance for more subtle extremality statements.
We expose them in the statements below.

All the separation statements in this section are consequences of the next general theorem, providing dual characterizations of the slightly weakened version of the asymmetric extremality property contained in Proposition~\ref{P2.7}(i).
It involves the $d_1$ distance between $n$ sets defined in \eqref{dd1} and actually combines two statements: for closed sets in general Banach spaces and specifically in Asplund spaces.
So far it has been common to formulate (and prove!) such statements separately; cf. \cite{Mor06.1,ZheNg11}.

\begin{theorem}\label{T5.1}
Suppose $\Omega_1,\ldots,\Omega_n$ are closed subsets of a Banach space $X$,
$\bx\in\cap_{i=1}^n\Omega_i$, $a_i\in{X}$ $(i=1,\ldots,n-1)$, and $\eps>0$.
Suppose also that
\begin{gather}\label{T5.1-1}
\bigcap\limits_{i=1}^{n-1}(\Omega_i-a_i)\cap\Omega_{n}
=\emptyset,
\\\label{T5.1-2}
\max_{1\le i\le n-1}\|a_i\|< d_1(\Omega_1-a_1,\ldots,\Omega_{n-1}-a_{n-1},\Omega_{n}) +\eps
\end{gather}
(or simply $\max_{1\le i\le n-1}\|a_i\|<\eps$).
Then,
\begin{enumerate}
\item
for any $\la>0$ and $\rho>0$,
there exist points
$\omega_i\in\Omega_i\cap{B}_\la(\bx)$ $(i=1,\ldots,n-1)$, $\omega_n\in\Omega_n\cap B_\rho(\bx)$, and vectors
$x_i^*\in X^*$ $(i=1,\ldots,n)$ such that
\begin{gather}\label{T5.1-3}
\sum_{i=1}^{n}x_i^*=0,\quad \sum_{i=1}^{n-1}\norm{x_i^*}=1,
\\\label{T5.1-4}
{\la}{\sum_{i=1}^{n-1} d\left(x_i^*, {N}_{\Omega_i}(\omega_i)\right) + {\rho}d\left(x_n^*,{N}_{\Omega_n}(\omega_n)\right) <{\eps}},
\\\label{T5.1-6}
\sum_{i=1}^{n-1}\ang{x_{i}^*,\omega_n+a_i-\omega_{i}} =\max_{1\le{i}\le{n}-1} \|\omega_n+a_i-\omega_{i}\|,
\end{gather}
where $N$ in \eqref{T5.1-4} stands for the Clarke normal cone ($N:=N^C$);
\item
if $X$ is Asplund,
then, for any $\la>0$, $\rho>0$ and $\tau\in]0,1[$,
there exist points
$\omega_i\in\Omega_i\cap{B}_\la(\bx)$ $(i=1,\ldots,n-1)$, $\omega_n\in\Omega_n\cap B_\rho(\bx)$, and vectors
$x_i^*\in X^*$ $(i=1,\ldots,n)$ satisfying conditions \eqref{T5.1-3}, \eqref{T5.1-4}, where $N$ in \eqref{T5.1-4} stands for the \Fr\ normal cone ($N:=N^F$), and
\begin{gather}\label{T5.1-5}
\sum_{i=1}^{n-1}\ang{x_{i}^*,\omega_n+a_i-\omega_{i}} >\tau\max_{1\le{i}\le{n}-1} \|\omega_n+a_i-\omega_{i}\|.
\end{gather}
\end{enumerate}
\end{theorem}

\begin{Proof}
Let $\la>0$ and $\rho>0$.
In view of \eqref{T5.1-2}, we can choose positive numbers $\eps_1$ and $\eps_2$ satisfying
\begin{equation}\label{T5.7P1}
\max_{1\le i\le n-1}\|a_i\| -d_1(\Omega_1-a_1,\ldots,\Omega_{n-1}-a_{n-1},\Omega_{n}) <\eps_1<\eps_2<\eps
\end{equation}
Note that $\max_{1\le i\le n-1}\|a_i\| =d_1(\bx-a_1,\ldots,\bx-a_{n-1},\bx)$.
We can apply Lemma~\ref{L4.1} to find points $\hat\omega_i\in\Omega_i\cap B_\la(\bx)$ ${(i=1,\ldots,n-1)}$ and $\hat\omega_n\in\Omega_n\cap B_\rho(\bx)$ such that
\begin{align}\notag
\max_{1\le i\le n-1} \norm{\hat\omega_i-a_i-\hat\omega_n} 
<d_{1}(&\Omega_1\cap B_{\alpha\la}(\hat\omega_1)-a_1,\ldots,
\Omega_{n-1}\cap B_{\alpha\la}(\hat\omega_{n-1})-a_{n-1}, \\\label{T5.1P1}
&\Omega_n\cap B_{\alpha\rho}(\hat\omega_n)) +{\alpha\eps_1}
\end{align}
for all $\alpha>0$.
Consider the three functions $f_1,f_2,f_3:X^{n}\to\R_+\cup\{+\infty\}$:
\begin{align}\label{T5.1P2}
f_1(u_1,\ldots,u_{n}):=&\max_{1\le{i}\le{n}-1} \norm{u_i-a_i-u_n},
\\\label{T5.1P2.2}
f_2(u_1,\ldots,u_{n}):=
&{\eps_2} \max\left\{{\la\iv\max_{1\le{i}\le{n}-1} \|u_i-\hat\omega_{i}\|, \rho\iv\norm{u_n-\hat\omega_{n}}}\right\},
\\\label{T5.1P2.3}
f_3(u_1,\ldots,u_{n}):=&
\begin{cases}
0&\mbox{if } u_i\in\Omega_i\; (i=1,\ldots,n),
\\
\infty&\mbox{otherwise}.
\end{cases}
\end{align}
Observe that, in view of \eqref{T5.1-1}, $f_1(\hat\omega_1,\ldots,\hat\omega_n) =\max_{1\le{i}\le{n}-1} \|\hat\omega_{i}-a_i-\hat\omega_n\|>0$.
Moreover,
\begin{gather*}
f_1(u_1,\ldots,u_n)-f_1(\hat\omega_1,\ldots,\hat\omega_n) +f_2(u_1,\ldots,u_n)\ge0
\qdtx{for all}u_i\in\Omega_i\;(i=1,\ldots,n).
\end{gather*}
Indeed, assume that there are $u_i\in\Omega_i$ $(i=1,\ldots,n)$
such that the inequality does not hold.
Then, $(u_1,\ldots,u_n)\neq (\hat\omega_{i},\ldots,\hat\omega_n)$, and consequently,
$f_2(u_1,\ldots,u_n)>0$.
Set $\alpha :=f_2(u_1,\ldots,u_n)/\eps_1$.
We have $\norm{u_i-\hat\omega_{i}} \le\frac{\eps_1}{\eps_2}\alpha\la<\alpha\la$ $(i=1,\ldots,n-1)$, $\norm{u_n-\hat\omega_{n}} \le\frac{\eps_1}{\eps_2}\alpha\rho<\alpha\rho$, and
$$
\max_{1\le i\le n-1}\norm{u_i-a_i-u_n}-\max_{1\le i\le n-1}\norm{\hat\omega_i-a_i-\hat\omega_n}+{\alpha\eps_1}<0,
$$
which contradicts \eqref{T5.1P1}.
Thus,
$(\hat\omega_1,\ldots,\hat\omega_n)$ is a point of minimum of the sum $f_1+f_2+f_3$, and consequently (Lemma~\ref{L2.3}), $0\in\sd(f_1+f_2+f_3)(\hat\omega_1,\ldots,\hat\omega_n)$.

Functions $f_1$ and $f_2$ are convex and Lipschitz continuous.
It is easy to check that the
subdifferentials of $f_1$, $f_2$ and $f_3$ possess the following
properties:

1) The function $f_1$ \eqref{T5.1P2} is a composition function:
$f_1(u_1,\ldots,u_n)= g\big(A(u_1,\ldots,u_n)-(a_1,\ldots,a_{n-1})\big)$, where $A$ is the linear operator from $X^{n}$ to $X^{n-1}$: $A(u_1,\ldots,u_n):=\left(u_1-u_n,\ldots, u_{n-1}-u_n\right)$, and $g$ is the maximum norm on $X^{n-1}$: ${g(u_1,\ldots,u_{n-1}):= \max_{1\le i\le n-1}\|u_i\|}$.
The corresponding dual norm has the form $(v_1^*,\ldots,v_{n-1}^*)\mapsto \sum_{i=1}^{n-1}\|v_i^*\|$.
It is easy to check that the adjoint operator $A^*:(X^*)^{n-1}\to(X^*)^{n}$ is of the form $A^*(v_1^*,\ldots,v_{n-1}^*) =\left(v_1^*,\ldots, v_{n-1}^*,-\sum_{j=1}^{n-1}v_j^*\right)$.
When $f_1(u_1,\ldots,u_n)>0$, the subdifferential $\sd{g}\left(u_1-a_1-u_n,\ldots, u_{n-1}-a_{n-1}-u_n\right)$ is the set of $(v_{11}^*,\ldots,v_{1,n-1}^*)\in (X^*)^{n-1}$ satisfying
(see, e.g., \cite[Corollary~2.4.16]{Zal02})
\sloppy
\begin{gather}\label{T5.1P3}
\sum_{i=1}^{n-1}\|v_{1i}^*\|=1
\quad\mbox{and}\quad
\sum_{i=1}^{n-1}\ang{v_{1i}^*,u_{i}-a_i-u_n} =\max_{1\le{i}\le{n-1}} \norm{u_{i}-a_i-u_n}.
\end{gather}
Thus, in view of the convex chain rule (see, e.g., \cite[Theorem~2.8.3]{Zal02}),
if
${f}_1(u_1,\ldots,u_n)>0$,
then the subdifferential $\sd{f}_1(u_1,\ldots,u_n)$ is the set of all vectors
$\left(v_{11}^*,\ldots, v_{1,n-1}^*,
\allowbreak
-\sum_{j=1}^{n-1}v_{1j}^*\right)\in (X^*)^{n},$
where vectors $v_{1i}^*\in X^*$ $(i=1,\ldots,n-1)$ satisfy \eqref{T5.1P3}.

2) The function $f_2$ \eqref{T5.1P2.2} is a positive multiple of the norm on $X^{n}$ (translated by $(\hat\omega_1,\ldots,\hat\omega_n)$).
Its subgradients $(v_{21}^*,\ldots,v_{2n}^*)$ at any point satisfy
\begin{gather}\label{T5.1P4} \la\sum_{i=1}^{n-1}\norm{v_{2i}^*}+\rho\norm{v_{2n}^*} \le\eps_2.
\end{gather}

3) The function $f_3$ \eqref{T5.1P2.3} is the indicator function of the set $\Omega_1\times\ldots\times\Omega_n$.
Its subdifferential has a simple representation:
$\sd{f}_3(u_1,\ldots,u_n) =\prod_{i=1}^nN_{\Omega_i}(u_i)$ for all $u_i\in\Omega_i$ $(i=1,\ldots,n)$ (Lemma~\ref{L2.4}).

From this point, the proof splits into two cases.

(i)
We apply the Clarke--Rockafellar subdifferential sum rule (Lemma~\ref{SR}(iv))
to find elements of the three subdifferentials:  $(v_{11}^*,\ldots,v_{1,n-1}^*)\in \sd{g}(\hat\omega_1-a_1-\hat\omega_n,\ldots, \hat\omega_{n-1}-a_{n}-\hat\omega_n)$, $(v_{21}^*,\ldots,v_{2n}^*)\in \sd{f}_2(\hat\omega_{1},\ldots,\hat\omega_{n})$ and ${(v_{31}^*,\ldots,v_{3n}^*)\in \sd{f}_3(\hat\omega_{1},\ldots,\hat\omega_{n})}$ such that
\begin{gather*}
v_{1i}^*+v_{2i}^*+v_{3i}^*=0\;(i=1,\ldots,n-1)
\qdtx{and}
-\left(\sum_{j=1}^{n-1}v_{1j}^*\right) +v_{2n}^*+v_{3n}^*=0.
\end{gather*}
Then, $v_{3i}^*\in{N}_{\Omega_i}(\hat\omega_i)$ $(i=1,\ldots,n)$ and conditions \eqref{T5.1P3} and \eqref{T5.1P4} are satisfied.
The conclusion of part (i) of the theorem holds true with $\hat\omega_i$ in place of $\omega_i$ $(i=1,\ldots,n)$,
$x_i^*:=-v_{1i}^*$ $(i=1,\ldots,{n}-1)$ and
{$x_n^*:=\sum_{j=1}^{n-1}v_{1j}^*$}.

(ii)
Let $X$ be Asplund and
$\tau\in]0,1[$.
We can apply the fuzzy sum rule (Lemma~\ref{SR}(iii)) to the sum $(f_1+f_2)+f_3$ followed by the conventional convex sum rule (Lemma~\ref{SR}(i)) applied to $f_1+f_2$.
Choose a $\xi>0$ satisfying the following conditions:
\begin{gather}\label{T5.1P5}
{\xi<{\eps-\eps_2}},\quad
\xi<\la-\max_{1\le{i}\le{n-1}} \|\hat\omega_i-\bx\|,\quad
\xi<\rho-\|\hat\omega_n-\bx\|,
\\\label{T5.1P6}
(10-2\tau)\xi<(1-\tau) \max_{1\le{i}\le{n}-1}\|\hat\omega_{i}-a_i-\hat\omega_n\|.
\end{gather}
Since $\tau\in]0,1[$, the last inequality implies in particular that
\begin{gather}\label{T5.1P62}
2\xi <\max_{1\le{i}\le{n-1}} \|\hat\omega_{i}-a_i-\hat\omega_n\|.
\end{gather}
Applying the fuzzy sum rule, we find two points $(x_{1},\ldots,x_{n})$, $(\omega_1,\ldots,\omega_n)\in X^{n}$ such that $\omega_i\in\Omega_i$ $(i=1,\ldots,n)$ and
\begin{gather}\label{T5.1P7}
\max_{1\le i\le n}\|x_i-\hat\omega_i\|<\xi,\quad \max_{1\le i\le n}\|\omega_i-\hat\omega_i\|<\xi,
\end{gather}
and elements of the three subdifferentials: $(v_{11}^*,\ldots,v_{1,n-1}^*)\in \sd{g}(x_1-a_1-x_n,\ldots,x_{n-1}-a_{n}-x_n)$, $(v_{21}^*,\ldots,v_{2n}^*)\in \sd{f}_2(x_{1},\ldots,x_{n})$ and $(v_{31}^*,\ldots,v_{3n}^*)\in \sd{f}_3(\omega_{1},\ldots,\omega_{n})$ such that
\begin{gather*}
\sum_{i=1}^{n-1} \norm{v_{1i}^*+v_{2i}^*+v_{3i}^*}+
\norm{-\left(\sum_{j=1}^{n-1}v_{1j}^*\right) +v_{2n}^*+v_{3n}^*} <\frac{\xi}{\max\{\la,\rho\}}.
\end{gather*}
The last condition implies
\begin{gather}\label{T5.1P8}
\la\sum_{i=1}^{n-1} \norm{v_{1i}^*+v_{2i}^*+v_{3i}^*}+ \rho
\norm{-\left(\sum_{j=1}^{n-1}v_{1j}^*\right) +v_{2n}^*+v_{3n}^*} <\xi.
\end{gather}
In view of
\eqref{T5.1P5}, \eqref{T5.1P6}, \eqref{T5.1P62} and \eqref{T5.1P7},
we have the following estimates:
\begin{align*}
&\|\omega_i-\bx\|\le\|\hat\omega_i-\bx\|+\|\omega_i-\hat\omega_i\| <\la
{\quad (i=1,\ldots,n-1)},
\\
&\norm{\omega_n-\bx} <\norm{\hat{\omega}_n-\bx}+\|\omega_n-\hat\omega_n\|<\rho,
\end{align*}
(hence, $\omega_i\in\Omega_i\cap{B}_\la(\bx)$
{$(i=1,\ldots,n-1)$ and $\omega_n\in\Omega_n\cap B_\rho(\bx)$})
\begin{align}\label{T5.1P9}
\norm{\omega_i-x_i}&\le\norm{\omega_i-\hat\omega_i} +\norm{x_i-\hat\omega_i}<2\xi
{\quad (i=1,\ldots,n)},
\\\notag
\max_{1\le{i}\le{n}-1}\|x_{i}-a_i-x_n\| &\ge\max_{1\le{i}\le{n}-1} \left(\|\hat\omega_{i}-a_i-\hat\omega_n\| -\|x_i-\hat\omega_i\|-\norm{x_n-\hat\omega_n}\right)
\\\notag
&>\max_{1\le{i}\le{n}-1}\|\hat\omega_{i}-a_i-\hat\omega_n\| -2\xi>0,
\\\notag
\max_{1\le{i}\le{n}-1} \|\omega_{i}-a_i-\omega_n\|&\ge\max_{1\le{i}\le{n}-1} \left(\|\hat\omega_{i}-a_i-\hat\omega_n\| -\|\omega_i-\hat\omega_i\|-\norm{\omega_n-\hat\omega_n}\right)
\\\notag
&>\max_{1\le{i}\le{n}-1} \|\hat\omega_{i}-a_i-\hat\omega_n\|-2\xi>0,
\\\label{T5.1P10}
(1-\tau)\max_{1\le{i}\le{n}-1} \|\omega_{i}-a_i-\omega_n\| &>(1-\tau)\left(\max_{1\le{i}\le{n}-1} \|\hat\omega_{i}-a_i-\hat\omega_n\|-2\xi\right)
\\\notag
&>(10-2\tau)\xi-2(1-\tau)\xi=8\xi,
\end{align}
conditions \eqref{T5.1P3} and \eqref{T5.1P4} are satisfied and $v_{3i}^*\in{N}_{\Omega_i}(\omega_i)$ $(i=1,\ldots,n)$.
Denote $x_i^*:=-v_{1i}^*$ $(i=1,\ldots,{n}-1)$ and
{$x_n^*:=\sum_{j=1}^{n-1}v_{1j}^*$}.
Then, $\sum_{i=1}^{n}x_i^*=0$ and $\sum_{i=1}^{n}\norm{x_i^*}=1$.
Using \eqref{T5.1P3} and \eqref{T5.1P9}, we obtain the following inequalities:
\begin{align}\notag
\sum_{i=1}^{n-1}\ang{x_{i}^*,\omega_n+a_i-\omega_{i}}
&\ge\sum_{i=1}^{n-1}\ang{v_{1i}^*,x_{i}-a_i-x_n} -\sum_{i=1}^{n-1}\norm{v_{1i}^*}(\norm{\omega_i-x_{i}} +\norm{\omega_n-x_n})
\\\notag
&>\max_{1\le{i}\le{n}-1} \|x_{i}-a_i-x_n\|-4\xi
\\\notag
&\ge\max_{1\le{i}\le{n}-1}\left(\|\omega_{i}-a_i-\omega_n\| -\norm{\omega_i-x_{i}} -\norm{\omega_n-x_n}\right)-4\xi
\\\label{T5.1P11}
&>\max_{1\le{i}\le{n}-1} \|\omega_{i}-a_i-\omega_n\|-8\xi.
\end{align}
Adding \eqref{T5.1P10} and \eqref{T5.1P11}, we arrive at \eqref{T5.1-5}.
Making use of \eqref{T5.1P4}, \eqref{T5.1P8} and \eqref{T5.1P5}, we obtain the following estimates:
\begin{gather*}
\la\sum_{i=1}^{n-1}\norm{x_i^*-v_{3i}^*} +\rho\norm{x_n^*-v^*_{3n}} <\la\sum_{i=1}^{n-1}\norm{v_{2i}^*} +\rho\norm{v^*_{2n}}+\xi <{\eps_2}+\xi <{\eps}.
\end{gather*}
The last inequality yields \eqref{T5.1-4}.
\qed\end{Proof}

\begin{Remark}\label{R5.2}
\begin{enumerate}
\item
Conditions \eqref{T5.1-1} and \eqref{T5.1-2} are implied by conditions \eqref{P4}.
Hence, in view of Proposition~\ref{P2.7}(i), Theorem~\ref{T5.1} provides dual necessary characterizations of extremality.
\item
Conditions \eqref{T5.1-6} and \eqref{T5.1-5} relate the dual vectors $x_i^*$ and the primal space vectors $\omega_n+a_i-\omega_{i}$ $(i=1,\ldots,n-1)$.
Such conditions, though not common in the conventional formulations of the extremal/generalized separation statements, seem to provide important additional characterizations of the properties.
Conditions of this kind first appeared explicitly in the generalized separation theorems in \cite{ZheNg11}, where the authors also provided
motivations for
employing such conditions.
As one can see from the proof above, conditions \eqref{T5.1-6} and \eqref{T5.1-5} originate in computing the convex subdifferential of the norm in $X^{n-1}$ at a nonzero point; see \eqref{T5.1P3}.
Subdifferentiating a norm (in either $X^{n-1}$ or $X^{n}$) at a nonzero point is a necessary step in the proofs of all existing versions of the extremal principle and its extensions, starting with the very first one in \cite[Theorem~6.1]{KruMor80.2}, with conditions like \eqref{T5.1P3} hidden in the proofs.
In several statements in the rest of this section, following \cite{ZheNg11}, we make such conditions exposed.
\item
Lemma~\ref{L4.1} substitutes in the proof of Theorem~\ref{T5.1} the conventional \EVP.
\item
If condition \eqref{T5.1-2} in Theorem~\ref{T5.1} is replaced by a stronger one:
\begin{gather*}
\max_{1\le i\le n-1}\|a_i\|= d_1(\Omega_1-a_1,\ldots,\Omega_{n-1}-a_{n-1},\Omega_{n}),
\end{gather*}
which means that the infimum of $\max_{1\le i\le n-1} \|\omega_i-a_i-\omega_n\|$ over $\omega_i\in\Omega_i$ $(i=1,\ldots,n)$ is attained at ${\omega_1=\ldots=\omega_n=\bx}$, then the application of Lemma~\ref{L4.1} in the proof can be dropped, leading to an improvement in part (i): conditions \eqref{T5.1-4} and \eqref{T5.1-6} can be replaced by $x_i^*\in{N}_{\Omega_i}^C(\bx)$ $(i=1,\ldots,n)$ and
$\sum_{i=1}^{n-1}\ang{x_{i}^*,a_i} =\max_{1\le{i}\le{n}-1} \|a_i\|$, respectively.
A similar fact was observed in \cite[Theorem~3.1$^\prime$]{ZheNg11}.
\end{enumerate}
\end{Remark}

The assumption $\cap_{i=1}^n\Omega_i\ne\es$ in Theorem~\ref{T5.1} is not restrictive.
The common point
$\bx\in\cap_{i=1}^n\Omega_i$ of the collection of sets can be replaced by a collection of individual points $\omega_i\in\Omega_i$ $(i=1,\ldots,n)$.
The next statement provides dual characterizations of the slightly weakened version of the asymmetric extremality property contained in Proposition~\ref{P3.11}.

\begin{corollary}\label{C5.3}
Suppose $\Omega_1,\ldots,\Omega_n$ are closed subsets of a Banach space $X$,
$\omega_i\in\Omega_i$ $(i=1,\ldots,n)$, $a_i\in{X}$ ${(i=1,\ldots,n-1)}$, and $\eps>0$.
Suppose also that
\begin{gather}\label{C5.3-1}
\bigcap\limits_{i=1}^{n-1}(\Omega_i-\omega_i-a_i) \cap(\Omega_{n}-\omega_n)
=\emptyset,
\\\label{C5.3-2}
\max_{1\le i\le n-1}\|a_i\| <d_1(\Omega_1-\omega_1-a_1,\ldots, \Omega_{n-1}-\omega_{n-1}-a_{n-1},\Omega_{n}-\omega_{n}) +\eps
\end{gather}
(or simply $\max_{1\le i\le n-1}\|a_i\|<\eps$).
Then,
\begin{enumerate}
\item
for any $\la>0$ and $\rho>0$,
there exist points
$\omega'_i\in\Omega_i\cap{B}_\la(\omega_i)$ $(i=1,\ldots,n-1)$, $\omega'_n\in\Omega_n\cap B_\rho(\omega_n)$, and vectors
$x_i^*\in X^*$ $(i=1,\ldots,n)$ satisfying conditions \eqref{T5.1-3} and \eqref{T5.1-4} with $N:=N^C$ and $\omega'_i$ in place of $\omega_i$ $(i=1,\ldots,n)$, and
\begin{gather*}
\sum_{i=1}^{n-1} \ang{x_{i}^*,\omega_n'-\omega_n+a_i-\omega_{i}'+\omega_{i}} =\max_{1\le{i}\le{n}-1} \|\omega_n'-\omega_n+a_i-\omega_{i}'+\omega_{i}\|;
\end{gather*}
\item
if $X$ is Asplund,
then, for any $\la>0$, $\rho>0$ and $\tau\in]0,1[$,
there exist points
$\omega'_i\in\Omega_i\cap{B}_\la(\omega_i)$ $(i=1,\ldots,n-1)$, $\omega'_n\in\Omega_n\cap B_\rho(\omega_n)$, and vectors
$x_i^*\in X^*$ $(i=1,\ldots,n)$ satisfying conditions \eqref{T5.1-3} and \eqref{T5.1-4} with $N:=N^F$ and $\omega'_i$ in place of $\omega_i$ $(i=1,\ldots,n)$,
and
\begin{gather*}
\sum_{i=1}^{n-1} \ang{x_{i}^*,\omega_n'-\omega_n+a_i-\omega_{i}'+\omega_{i}} >\tau\max_{1\le{i}\le{n}-1} \|\omega_n'-\omega_n+a_i-\omega_{i}'+\omega_{i}\|.
\end{gather*}
\end{enumerate}
\end{corollary}

\begin{Proof}
The sets $\Omega_i':=\Omega_i-\omega_i$ $(i=1,\ldots,n)$ satisfy $0\in\cap_{i=1}^n\Omega_i'$ and conditions \eqref{T5.1-1} and \eqref{T5.1-2}.
The conclusion follows from Theorem~\ref{T5.1} after noticing that ${N}_{\Omega_i'}(\omega_i'-\omega_i) ={N}_{\Omega_i}(\omega_i')$ $(i=1,\ldots,n)$.
\qed\end{Proof}

Theorem~\ref{T5.1} is a particular case of Corollary~\ref{C5.3} with $\omega_i=\bx$ $(i=1,\ldots,n)$.

The following theorem is an immediate consequence of Theorem~\ref{T5.1}.
It combines two \emph{unified separation theorems} due to Zheng and Ng \cite[Theorems~3.1 and 3.4]{ZheNg11}.

\begin{theorem}\label{ZhNg}
Suppose $\Omega_1,\ldots,\Omega_n$ are closed subsets of a Banach space $X$, $\cap_{i=1}^n \Omega_i=\emptyset$, $\omega_i\in\Omega_i$ $(i=1,\ldots,n)$,
$\eps>0$ and
condition \eqref{P3.11-1} is satisfied.
Then,
\begin{enumerate}
\item
for any $\la>0$, there exist points
$\omega'_i\in\Omega_i\cap B_\la(\omega_i)$ and vectors $x^*_i\in X^*$ $(i=1,\ldots,n)$ satisfying conditions \eqref{D6} with $N$ standing for the Clarke normal cone ($N:=N^C$), $\al:=\frac{\eps}{\la}$ and $\omega'_i$ in place of $\omega_i$ $(i=1,\ldots,n)$, and
\begin{gather}\label{ZN_2}
\sum_{i=1}^{n-1}\langle x^*_i, \omega'_n-\omega'_i\rangle =\max_{1\le i \le n-1} \norm{\omega'_i-\omega'_n};
\end{gather}
\item
if $X$ is Asplund,
then, for any numbers $\la>0$ and $\tau\in]0,1[$, there exist points
$\omega'_i\in\Omega_i\cap B_\la(\omega_i)$ and vectors $x^*_i\in X^*$ $(i=1,\ldots,n)$ satisfying conditions \eqref{D6} with $N$ standing for the \Fr\ normal cone ($N:=N^F$), $\al:=\frac{\eps}{\la}$ and $\omega'_i$ in place of $\omega_i$ $(i=1,\ldots,n)$, and
\begin{gather}\label{ZN_3}
\sum_{i=1}^{n-1}\langle x^*_i, \omega'_n-\omega'_i\rangle >\tau\max_{1\le i \le n-1}\norm{\omega'_i-\omega'_n}.
\end{gather}
\end{enumerate}
\end{theorem}

\begin{Proof}
Observe that the sets $\Omega_i':=\Omega_i-\omega_i$ $(i=1,\ldots,n)$ and vectors $a_i:=\omega_n-\omega_i$ $(i=1,\ldots,n-1)$ satisfy $0\in\bigcap_{i=1}^n\Omega_i'$ and
\begin{align*}
\bigcap\limits_{i=1}^{n-1}(\Omega_i'-a_i)\cap\Omega_{n}'
&=\bigcap\limits_{i=1}^{n}(\Omega_i-\omega_n)
=\bigcap\limits_{i=1}^{n}\Omega_i-\omega_n
=\emptyset,
\\
\max_{1\le i\le n-1}\|a_i\| &=d_1(\omega_1,\ldots,\omega_n)< d_1(\Omega_1,\ldots,\Omega_n)+\eps
\\
&= d_1(\Omega_1'-a_1,\ldots,\Omega_{n-1}'-a_{n-1},\Omega_n') +\eps.
\end{align*}
Applying Theorem~\ref{T5.1} with $\rho=\la$, we arrive at the conclusions.
\qed\end{Proof}

\if{
\blue{The following corollaries are local versions and also direct consequences of Theorem~\ref{ZhNg}.}

\blue{\begin{corollary}\label{C5.7.1}
		Suppose $\Omega_1,\ldots,\Omega_n$ are closed subsets of a Asplund space $X$,
	$\bx\in\cap_{i=1}^n\Omega_i$, $a_i\in{X}$ $(i=1,\ldots,n)$, $\rho>0$ and $\eps>0$.
		Suppose also that \eqref{T5.7-3} holds, then,
for any $\la\in ]0,\rho]$ and $\tau\in]0,1[$,
there exist points
	$\omega_i\in\Omega_i\cap{B}_\la(\bx)$ $(i=1,\ldots,n)$, and vectors
	$x_i^*\in X^*$ $(i=1,\ldots,n)$ such that \eqref{T5.1-3}
\begin{gather}\label{C5.7.1-1}
\sum_{i=1}^n d(x_i^*,N_{\Omega_i}(\omega_i))< \frac{\eps}{\la}, \AND \\
\label{C5.7.1-2}\sum_{i=1}^{n-1}\ang{x_{i}^*,\omega_n+a_i-\omega_{i}-a_n} >\tau\max_{1\le{i}\le{n}-1} \|\omega_n+a_i-\omega_{i}-a_n\|.
\end{gather}	
\end{corollary}}

\AK{4/07/18.
I am unable to make sense of the above corollary.}
\blue{
\begin{Proof}
		Let \eqref{T5.7-3} hold for some $a_i\in X$ ($i=1,\ldots,n$), $\rho>0$ and $\eps>0$.
		Take some positive numbers $\eps_1<\eps$, $\rho_1<\rho$, $\la_1<\la$ such that $\la\eps_1\eps^{-1}<\la_1\le \rho_1$, and
$$
\max_{1\le i \le n} \norm{a_i} < d_{B_\rho(\bx)}(\Omega_1-a_1,\ldots,\Omega_n-a_n)+\eps< d_1\left((\Omega_1-a_1)\cap \overline{B}_{\rho_1}(\bx),\ldots,(\Omega_n-a_n)\cap \overline{B}_{\rho_1}(\bx)\right)+\eps_1.
$$
		Observe that the closed sets $(\Omega_i-a_i)\cap \overline{B}_{\rho_1}(\bx)$ ($i=1,\ldots,n$) have empty intersection and
	condition~\ref{P3.11-1} is satisfied with $\omega_i-a_i$ in place of $\omega_i$ and $(\Omega_i-a_i)\cap \overline{B}_{\rho_i}$ in place of $\Omega_i$ ($i=1,\ldots,n$).
		Applying Theorem~\ref{ZhNg}, we find points $\omega_i\in \Omega_i\cap \overline{B}_{\rho_i}(\bx)$ and vectors $x_i^*\in X^*$ ($i=1,\ldots,n$) satisfying \eqref{T5.1-3}, \eqref{C5.7.1-2}, and
	\begin{gather*}
	\sum_{i=1}^nd\left(x_i^*,N_{(\Omega_i-a_i)\cap \overline{B}_{\rho_1}(\bx)}(\omega_i-a_i)\right)<\frac{\eps_1}{\la_1}\le \frac{\eps}{\la}.
	\end{gather*}
		Observe further that $N_{(\Omega_i-a_i)\cap \overline{B}_{\rho_1}(\bx)}(\omega_i-a_i)=N_{\Omega_i}(\omega_i)$ for all $i=1,\ldots,n$, we have \eqref{C5.7.1-1} holds.
\end{Proof}
}
\red{\begin{corollary}\label{C5.8.1}
		Suppose $\Omega_1,\ldots,\Omega_n$ are closed subsets of an Asplund space $X$,
		$\bx\in\cap_{i=1}^n\Omega_i$ and $\al\in ]0,1[$.
		Suppose also that condition {\rm (i)} in Corollary~\ref{C4.9} is satisfied.
		Then, for any $\eps>0$ and $\tau \in ]0,1[$,
		there exist points $\omega_i\in \Omega_i\cap B_{\eps}(\bx)$,  and vectors $x^*_i\in X^*$ $(i=1,\ldots,n)$
		such that
		conditions
		\begin{gather}
		\label{C5.8.1-1}
		\sum_{i=1}^n x^*_i=0,\quad
		\sum_{i=1}^{n}\norm{x_i^*}=1,
		\quad
		\sum_{i=1}^nd(x^*_i,N_{\Omega_i}(\omega_i))<\al,
		\end{gather}
		and \eqref{C5.7.1-2} are satisfied.
\end{corollary}}

\begin{Proof}
	Let a number $\eps>0$ be given.
	Choose a number $\xi=\frac{\eps}{2}$.
	By condition {\rm (i)} in Corollary~\ref{C4.9}, there exist a number $\rho\in]0,\xi[$, points $\omega_i\in \Omega_i\cap B_\xi(\bx)$ and vectors $a_i\in X$ $(i=1,\ldots,n)$ such that condition \eqref{P7} is satisfied.
	We have $\xi+\rho<\eps$, and the conclusion follows from Corollary~\ref{C5.7.1}.
	\qed\end{Proof}
\blue{To have such a local version follows by Theorem~\ref{ZhNg}, alike Corollary~\ref{C5.7.1}, the restriction of $\la$ can not avoid.}
}\fi

\begin{Remark}
\begin{enumerate}
\item
In \cite{ZheNg11}, instead of the $d_1$ distance in condition \eqref{P3.11-1}, a slightly more general $p$-weighted nonintersect index was used with the corresponding $q$-weighted sums replacing the usual ones in \eqref{D6}, \eqref{ZN_2} and \eqref{ZN_3}.
This corresponds to considering $l_p$ norms on product spaces and the corresponding $l_q$ dual norms.
In the current paper, for simplicity only the maximum norm on product spaces is considered together with the corresponding sum norm in the dual space; cf. Remark~\ref{R3.10}(ii).
\item
Theorem~\ref{ZhNg} is a consequence of Theorem~\ref{T5.1}, which in turn is a consequence of the \EVP.
Thanks to \cite[Theorem~3.1]{LiTanYuWei08}, part (i) of Theorem~\ref{ZhNg} is equivalent to the \EVP.
Hence, the conclusion of Theorem~\ref{T5.1} is also equivalent to the \EVP\ (and to completeness of the space $X$).
\end{enumerate}
\end{Remark}

The next theorem is a `symmetric' version of Theorem~\ref{T5.1}.

\begin{theorem}\label{T5.7}
Suppose $\Omega_1,\ldots,\Omega_n$ are closed subsets of a Banach space $X$,
$\bx\in\cap_{i=1}^n\Omega_i$, $a_i\in{X}$ $(i=1,\ldots,n)$, $\rho>0$ and $\eps>0$.
Suppose also that
\begin{gather}\label{T5.7-5}
\bigcap_{i=1}^n(\Omega_i-a_i)\cap{B}_\rho(\bar{x})
=\emptyset,
\\\label{T5.7-3}
\max_{1\le i\le n}\|a_i\|<d_{{B}_\rho(\bar{x})} (\Omega_1-a_1,\ldots,\Omega_{n}-a_{n})+\eps
\end{gather}
(or simply $\max_{1\le i\le n}\|a_i\|<\eps$).
Then,
\begin{enumerate}
\item
for any $\la>0$,
there exist points
$\omega_i\in\Omega_i\cap{B}_\la(\bx)$ $(i=1,\ldots,n)$ and $x\in B_\rho(\bx)$, and vectors
$x_i^*\in{X}^*$ $(i=1,\ldots,n)$ such that
\begin{gather}\label{T5.7-1}
\la\sum_{i=1}^nd(x_i^*,{N}_{\Omega_i}(\omega_i)) +\rho\norm{\sum_{i=1}^nx_i^*}<{\eps},\quad
\sum_{i=1}^n\norm{x_i^*}=1,
\\\label{T5.7-4}
\sum_{i=1}^{n}\ang{x_{i}^*,x+a_i-\omega_{i}} =\max_{1\le{i}\le{n}} \|x+a_i-\omega_{i}\|,
\end{gather}
where $N$ in \eqref{T5.7-1} stands for the Clarke normal cone ($N:=N^C$);
\item
if $X$ is Asplund,
then, for any $\la>0$ and $\tau\in]0,1[$,
there exist points
$\omega_i\in\Omega_i\cap{B}_\la(\bx)$ $(i=1,\ldots,n)$ and $x\in B_\rho(\bx)$, and vectors
$x_i^*\in X^*$ $(i=1,\ldots,n)$ satisfying conditions \eqref{T5.7-1}, where $N$ stands for the \Fr\ normal cone ($N:=N^F$), and
\begin{gather}\label{T5.7-0}
\sum_{i=1}^n\ang{x_{i}^*,x+a_i-\omega_{i}} >\tau\max_{1\le{i}\le{n}} \|x+a_i-\omega_{i}\|.
\end{gather}
\end{enumerate}
\end{theorem}

\begin{Proof}
Choose an $\eps'\in]0,\eps[$ such that condition \eqref{T5.7-3} holds true with $\eps'$ in place of $\eps$, and then choose a $\rho'\in]0,\rho[$ such that $\rho-\rho'<\eps-\eps'$.
It is sufficient to apply Theorem~\ref{T5.1} to the collection of $n+1$ closed sets $\Omega_1,\ldots,\Omega_n$ and $\Omega_{n+1}:=\overline{B}_{\rho'}(\bx)$ with $\eps'$ and $\rho'$ in place of $\eps$ and $\rho$, respectively.
Notice that $\Omega_{n+1}\cap B_{\rho'}(\bx)=B_{\rho'}(\bx)$ and $N_{\Omega_{n+1}(\bx)}(x)=0$ for any $x\in B_{\rho'}(x)$.
One only needs to check the inequality in \eqref{T5.7-1}, which is straightforward:
\begin{gather*}
\la \sum_{i=1}^nd(x_i^*,{N}_{\Omega_i}(\omega_i)) +\rho\norm{\sum_{i=1}^nx_i^*} <{\eps}'+(\rho-\rho')\norm{\sum_{i=1}^nx_i^*} \le{\eps}'+(\rho-\rho')<{\eps}.
\end{gather*}
The proof is complete.
\qed\end{Proof}

\begin{Remark}
\begin{enumerate}
\item
Conditions \eqref{T5.7-5} and \eqref{T5.7-3} are implied by conditions \eqref{P2}.
Hence, Theorem~\ref{T5.7} provides dual necessary characterizations of the local extremality.
\item
The inequality in \eqref{T5.7-1} combines two constraints on the vectors
$x_i^*\in X^*$ $(i=1,\ldots,n)$: they must be close to the respective normal cones and their sum must be small.
This inequality obviously implies two separate inequalities:
\begin{gather*}
\sum_{i=1}^nd(x_i^*,{N}_{\Omega_i}(\omega_i)) <\frac{\eps}{\la}
\qdtx{and}
\norm{\sum_{i=1}^nx_i^*}<\frac{\eps}{\rho},
\end{gather*}
while the converse implication is not true in general.
Both these constraints are involved in each of the two generalized separation conditions in the (extended) extremal principle discussed in Sect.~\ref{S2} (parts (ii) and (iii) of Theorem~\ref{T2.9}).
However, each of the generalized separation conditions in Theorem~\ref{T2.9} requires actually a stronger version of one of the constraints: either the vectors must belong to the respective normal cones (part (ii) of Theorem~\ref{T2.9}) or their sum must be exactly zero (part (iii) of Theorem~\ref{T2.9}).
Fortunately, as the next two corollaries show, the required stronger versions of (one of) the constraints are consequences of the combined condition \eqref{T5.7-1}
and Lemma~\ref{L2.12}.
\end{enumerate}
\end{Remark}

\begin{corollary}\label{C5.11+}
Suppose $\Omega_1,\ldots,\Omega_n$ are closed subsets of a Banach space $X$,
$\bx\in\cap_{i=1}^n\Omega_i$, $a_i\in{X}$ $(i=1,\ldots,n)$, $\eps>0$, $\la>0$, $\eps+\la\le\rho$, and conditions \eqref{T5.7-5} and \eqref{T5.7-3} are satisfied (the latter condition can be replaced by the simpler and stronger one: $\max_{1\le i\le n}\|a_i\|<\eps$).
Then,
\begin{enumerate}
\item
there exist points
$\omega_i\in\Omega_i\cap{B}_\la(\bx)$ $(i=1,\ldots,n)$ and $x\in B_\rho(\bx)$, and vectors
$x_i^*\in X^*$ $(i=1,\ldots,n)$, satisfying conditions \eqref{D4} with $N$ standing for the Clarke normal cone ($N:=N^C$) and $\al:=\frac{\eps}{\la}$, and
condition \eqref{T5.7-0} with $\tau:=\frac{\rho-\eps}{\rho+\eps}$;
\item
if $X$ is Asplund,
then, for any $\tau\in]0,\frac{\rho-\eps}{\rho+\eps}[$, there exist points
$\omega_i\in\Omega_i\cap{B}_\la(\bx)$ $(i=1,\ldots,n)$ and $x\in B_\rho(\bx)$, and vectors
$x_i^*\in X^*$ $(i=1,\ldots,n)$ satisfying condition \eqref{D4} with $N$ standing for the \Fr\ normal cone ($N:=N^F$) and $\al:=\frac{\eps}{\la}$, and condition \eqref{T5.7-0}.
\end{enumerate}
\end{corollary}

\begin{Proof}
Set $\tau':=\tau(1+\frac{\eps}{\rho})+\frac{\eps}{\rho}$, and
observe that $\tau'=1$ in part (i), and $\tau'\in]0,1[$ in part (ii).
Applying Theorem~\ref{T5.1}, we find points
$\omega_i\in\Omega_i\cap{B}_\la(\bx)$ $(i=1,\ldots,n)$ and $x\in B_\rho(\bx)$, and vectors
$x_i^*\in X^*$ $(i=1,\ldots,n)$ such that condition \eqref{T5.7-1} with the respective normal cone is satisfied and
\begin{gather}\label{C5.11+P2}
\sum_{i=1}^n\ang{x_{i}^*,x+a_i-\omega_{i}} \ge\tau'\max_{1\le{i}\le{n}} \|x+a_i-\omega_{i}\|.
\end{gather}
The assertion follows from Lemma~\ref{L2.12}(i) and (iii) after noticing that
$\frac{\tau'\rho-\eps}{\rho+\eps} =\tau$.
\qed\end{Proof}

\begin{corollary}\label{C5.12+}
Suppose $\Omega_1,\ldots,\Omega_n$ are closed subsets of a Banach space $X$,
$\bx\in\cap_{i=1}^n\Omega_i$, $a_i\in{X}$ $(i=1,\ldots,n)$, $\eps>0$, $\rho>0$, $\eps+\rho\le\la$, and conditions \eqref{T5.7-5} and \eqref{T5.7-3} are satisfied (the latter condition can be replaced by the simpler and stronger one: $\max_{1\le i\le n}\|a_i\|<\eps$).
Then,
\begin{enumerate}
\item
there exist points
$\omega_i\in\Omega_i\cap{B}_\la(\bx)$ $(i=1,\ldots,n)$ and $x\in B_\rho(\bx)$, and vectors
$x_i^*\in X^*$ $(i=1,\ldots,n)$, satisfying conditions \eqref{D3} with $N$ standing for the Clarke normal cone ($N:=N^C$) and $\al:=\frac{\eps}{\rho}$, and
condition \eqref{T5.7-0} with $\tau:=\frac{\la-\eps}{\la+\eps}$;
\item
if $X$ is Asplund,
then, for any $\tau\in]0,\frac{\la-\eps}{\la+\eps}[$, there exist points
$\omega_i\in\Omega_i\cap{B}_\la(\bx)$ $(i=1,\ldots,n)$ and $x\in B_\rho(\bx)$, and vectors
$x_i^*\in X^*$ $(i=1,\ldots,n)$ satisfying conditions \eqref{D3} with $N$ standing for the \Fr\ normal cone ($N:=N^F$) and $\al:=\frac{\eps}{\rho}$, and condition \eqref{T5.7-0}.
\end{enumerate}
\end{corollary}

\begin{Proof}
Set $\tau':=\tau(1+\frac{\eps}{\la})+\frac{\eps}{\la}$, and
observe that $\tau'=1$ in part (i), and $\tau'\in]0,1[$ in part (ii).
Applying Theorem~\ref{T5.1}, we find points
$\omega_i\in\Omega_i\cap{B}_\la(\bx)$ $(i=1,\ldots,n)$ and $x\in B_\rho(\bx)$, and vectors
$x_i^*\in X^*$ $(i=1,\ldots,n)$ such that conditions \eqref{T5.7-1} (with the respective normal cone) and \eqref{C5.11+P2} are satisfied.
The assertion follows from Lemma~\ref{L2.12}(ii) and (iii) after noticing that
$\frac{\tau'\la-\eps}{\la+\eps} =\tau$.
\qed\end{Proof}

The versions of Theorem~\ref{T5.7} and Corollaries~\ref{C5.11+} and \ref{C5.12+} with a common point
$\bx\in\cap_{i=1}^n\Omega_i$ of a collection of sets replaced by a collection of individual points $\omega_i\in\Omega_i$ $(i=1,\ldots,n)$, presented in the next three corollaries, follow immediately.


\begin{corollary}\label{C5.8}
Suppose $\Omega_1,\ldots,\Omega_n$ are closed subsets of a Banach space $X$,
$\omega_i\in\Omega_i$, $a_i\in{X}$ $(i=1,\ldots,n)$, $\rho>0$, and $\eps>0$.
Suppose also that
\begin{gather}\label{C5.8-1}
\bigcap_{i=1}^n(\Omega_i-\omega_i-a_i)\cap(\rho\B)
=\emptyset,
\\\label{C5.8-2}
\max_{1\le i\le n}\|a_i\|<d_{\rho\B} (\Omega_1-\omega_1-a_1,\ldots,\Omega_{n}-\omega_n-a_{n}) +\eps
\end{gather}
(or simply $\max_{1\le i\le n}\|a_i\|<\eps$).
Then,
\begin{enumerate}
\item
for any $\la>0$,
there exist points
$\omega'_i\in\Omega_i\cap{B}_\la(\omega_i)$ $(i=1,\ldots,n)$ and $x\in\rho\B$, and vectors
$x_i^*\in{X}^*$ $(i=1,\ldots,n)$ satisfying conditions \eqref{T5.7-1} with $N$ standing for the Clarke normal cone ($N:=N^C$) and $\omega'_i$ in place of $\omega_i$ $(i=1,\ldots,n)$, and
\begin{gather*}
\sum_{i=1}^{n} \ang{x_{i}^*,x+a_i-\omega_{i}'+\omega_{i}} =\max_{1\le{i}\le{n}} \|x+a_i-\omega_{i}'+\omega_{i}\|;
\end{gather*}
\item
if $X$ is Asplund,
then, for any $\la>0$ and $\tau\in]0,1[$,
there exist points
$\omega'_i\in\Omega_i\cap{B}_\la(\omega_i)$ $(i=1,\ldots,n)$ and $x\in\rho\B$, and vectors
$x_i^*\in X^*$ $(i=1,\ldots,n)$  satisfying conditions \eqref{T5.7-1} with $N$ standing for the \Fr\ normal cone ($N:=N^F$) and $\omega'_i$ in place of $\omega_i$ $(i=1,\ldots,n)$, and
\begin{gather}\label{C5.8-4}
\sum_{i=1}^{n} \ang{x_{i}^*,x+a_i-\omega_{i}'+\omega_{i}} >\tau\max_{1\le{i}\le{n}} \|x+a_i-\omega_{i}'+\omega_{i}\|.
\end{gather}
\end{enumerate}
\end{corollary}

\begin{corollary}\label{C5.13}
Suppose $\Omega_1,\ldots,\Omega_n$ are closed subsets of a Banach space $X$,
$\omega_i\in\Omega_i$, $a_i\in{X}$ $(i=1,\ldots,n)$, $\eps>0$, $\la>0$, $\eps+\la\le\rho$, and conditions \eqref{C5.8-1} and \eqref{C5.8-2} are satisfied.
(The latter condition can be replaced by the simpler and stronger one: $\max_{1\le i\le n}\|a_i\|<\eps$.)
Then,
\begin{enumerate}
\item
there exist points
$\omega'_i\in\Omega_i\cap{B}_\la(\omega_i)$ $(i=1,\ldots,n)$ and $x\in\rho\B$, and vectors
$x_i^*\in X^*$ $(i=1,\ldots,n)$ such that conditions \eqref{D4} are satisfied with $N$ standing for the Clarke normal cone ($N:=N^C$),  $\al:=\frac{\eps}{\la}$ and $\omega'_i$ in place of $\omega_i$ $(i=1,\ldots,n)$, as well as condition \eqref{C5.8-4} with $\tau:=\frac{\rho-\eps}{\rho+\eps}$;
\item
if $X$ is Asplund,
then, for any $\tau\in]0,\frac{\rho-\eps}{\rho+\eps}[$,
there exist points
$\omega'_i\in\Omega_i\cap{B}_\la(\omega_i)$ $(i=1,\ldots,n)$ and $x\in\rho\B$, and vectors
$x_i^*\in X^*$ $(i=1,\ldots,n)$ such that conditions \eqref{D4} are satisfied with $N$ standing for the \Fr\ normal cone ($N:=N^F$), $\al:=\frac{\eps}{\la}$ and $\omega'_i$ in place of $\omega_i$ $(i=1,\ldots,n)$, as well as condition \eqref{C5.8-4}.
\end{enumerate}
\end{corollary}

\begin{corollary}\label{C5.14+}
Suppose $\Omega_1,\ldots,\Omega_n$ are closed subsets of a Banach space $X$,
$\omega_i\in\Omega_i$, $a_i\in{X}$ $(i=1,\ldots,n)$, $\eps>0$, $\rho>0$, $\eps+\rho\le\la$, and conditions \eqref{C5.8-1} and \eqref{C5.8-2} are satisfied.
(The latter condition can be replaced by the simpler and stronger one: $\max_{1\le i\le n}\|a_i\|<\eps$.)
Then,
\begin{enumerate}
\item
there exist points
$\omega'_i\in\Omega_i\cap{B}_\la(\omega_i)$ $(i=1,\ldots,n)$ and $x\in\rho\B$, and vectors
$x_i^*\in X^*$ $(i=1,\ldots,n)$ such that conditions \eqref{D3} are satisfied with $N$ standing for the Clarke normal cone ($N:=N^C$), $\al:=\frac{\eps}{\rho}$ and $\omega'_i$ in place of $\omega_i$ $(i=1,\ldots,n)$, as well as condition \eqref{C5.8-4} with $\tau:=\frac{\la-\eps}{\la+\eps}$;
\item
if $X$ is Asplund,
then, for any $\tau\in]0,\frac{\la-\eps}{\la+\eps}[$,
there exist points
$\omega'_i\in\Omega_i\cap{B}_\la(\omega_i)$ $(i=1,\ldots,n)$ and $x\in\rho\B$, and vectors
$x_i^*\in X^*$ $(i=1,\ldots,n)$ such that conditions \eqref{D3} are satisfied with $N$ standing for the \Fr\ normal cone ($N:=N^F$), $\al:=\frac{\eps}{\rho}$ and $\omega'_i$ in place of $\omega_i$ $(i=1,\ldots,n)$, as well as condition \eqref{C5.8-4}.
\end{enumerate}
\end{corollary}

\begin{Remark}\label{R5.15}
In the above three corollaries, the assumption of the existence of a common point
$\bx\in\cap_{i=1}^n\Omega_i$ of a collection of sets, used in Theorem~\ref{T5.7} and Corollaries~\ref{C5.11+} and \ref{C5.12+}, is relaxed to that of a collection of individual points $\omega_i\in\Omega_i$ $(i=1,\ldots,n)$ (which always exist as long as all the sets are nonempty).
On the other hand, if such a point
$\bx\in\cap_{i=1}^n\Omega_i$ exists, it can be used along with the collection $\omega_i\in\Omega_i$ $(i=1,\ldots,n)$ to provide additional useful estimates. Indeed, if $\|\omega_i-\bx\|\le\xi$ $(i=1,\ldots,n)$ for some $\xi\ge0$ (one can take, e.g., $\xi:=\max_{1\le i\le n} \|\omega_i-\bx\|$), then each of the above three corollaries immediately gives $\|\omega'_i-\bx\|<\xi+\la$ $(i=1,\ldots,n)$.
This simple observation plays an important role in the proof of the extended extremal principle.
It is used also in the next statement, which is a consequence of Corollaries~\ref{C5.13} and \ref{C5.14+}.
\end{Remark}

\begin{proposition}\label{P5.16}
Suppose $\Omega_1,\ldots,\Omega_n$ are closed subsets of a Banach space $X$,
$\bx\in\cap_{i=1}^n\Omega_i$,
$\xi\ge0$,
${\omega_i\in\Omega_i\cap B_{\xi}(\bx)}$,
$a_i\in{X}$ $(i=1,\ldots,n)$, $\al\in]0,1[$, $\rho>0$,
$\xi+\rho(1-\al)<\de$,
and conditions \eqref{P7} are satisfied.
Then,
\begin{enumerate}
\item
there exist points $\omega'_i\in B_{\de}(\bx)$ and vectors $x^*_i\in X^*$ $(i=1,\ldots,n)$
such that conditions \eqref{D3} are satisfied with $N$ standing for the Clarke normal cone ($N:=N^C$) and $\omega'_i$ in place of $\omega_i$, as well as condition \eqref{C5.8-4} with $\tau:=\frac{1}{1+2\al}$;
\item
there exist points $\omega'_i\in B_{\de}(\bx)$ and vectors $x^*_i\in X^*$ $(i=1,\ldots,n)$
such that conditions \eqref{D4} are satisfied with $N$ standing for the Clarke normal cone ($N:=N^C$), and $\al':=\frac{\al}{1-\al}$ and $\omega'_i$ in place of $\al$ and $\omega_i$, respectively, as well as condition \eqref{C5.8-4} with $\tau:=\frac{1-\al}{1+\al}$.
\cnta
\end{enumerate}
Suppose $X$ is Asplund and $N$ stands for the \Fr\ normal cone ($N:=N^F$).
Then,
\begin{enumerate}
\cntb
\item
for any
$\tau\in]0,\frac{1}{1+2\al}[$,
there exist points $\omega'_i\in B_{\de}(\bx)$ and vectors $x^*_i\in X^*$ $(i=1,\ldots,n)$
such that conditions \eqref{D3} are satisfied with
$\omega'_i$ in place of $\omega_i$, as well as condition \eqref{C5.8-4};
\item
for any
$\tau\in]0,\frac{1-\al}{1+\al}[$,
there exist points $\omega'_i\in B_{\de}(\bx)$ and vectors $x^*_i\in X^*$ $(i=1,\ldots,n)$
such that conditions \eqref{D4} are satisfied with $\al':=\frac{\al}{1-\al}$ and $\omega'_i$ in place of $\al$ and $\omega_i$, respectively, as well as condition \eqref{C5.8-4}.
\end{enumerate}
\end{proposition}

\begin{Proof}
In view of Remark~\ref{R5.15}, assertions (i) and (iii) are consequences of Corollary~\ref{C5.13} with $\eps:=\al\rho$ and ${\la:=(1-\al)\rho}$, while assertions (ii) and (iv) are consequences of Corollary~\ref{C5.14+} with $\eps:=\al\rho$ and $\la:=(1+\al)\rho$.
\qed\end{Proof}

\begin{Remark}
The assertion in Proposition~\ref{P5.16}(iii) improves \cite[Theorem~3.1]{KruLop12.1}, which was used in \cite{KruLop12.1} as the main tool when extending the extremal principle to infinite collections of sets.
\end{Remark}

The above proposition yields dual characterizations of approximate $\al$-stationarity (Definition~\ref{D2.17}).

\begin{corollary}\label{C5.9}
Suppose $\Omega_1,\ldots,\Omega_n$ are closed subsets of a Banach space $X$,
$\bx\in\cap_{i=1}^n\Omega_i$ and $\al\in]0,1[$.
Suppose also that the collection $\{\Omega_1,\ldots,\Omega_n\}$ is approximately $\al$-stationary at $\bx$.
Then,
\begin{enumerate}
\item
for any $\eps>0$,
there exist points $\omega_i\in\Omega_i\cap B_{\eps}(\bx)$ and vectors $x^*_i\in X^*$ $(i=1,\ldots,n)$
such that conditions \eqref{D3} are satisfied with $N$ standing for the Clarke normal cone ($N:=N^C$)ж
\item
for any $\eps>0$,
there exist points $\omega_i\in\Omega_i\cap B_{\eps}(\bx)$ and vectors $x^*_i\in X^*$ $(i=1,\ldots,n)$
such that conditions \eqref{D4} are satisfied with $N$ standing for the Clarke normal cone ($N:=N^C$), and $\al':=\frac{\al}{1-\al}$ in place of $\al$.
\cnta
\end{enumerate}
Suppose $X$ is Asplund and $N$ stands for the \Fr\ normal cone ($N:=N^F$).
Then,
\begin{enumerate}
\cntb
\item
for any $\eps>0$,
there exist points $\omega_i\in\Omega_i\cap B_{\eps}(\bx)$ and vectors $x^*_i\in X^*$ $(i=1,\ldots,n)$
such that conditions \eqref{D3} are satisfied;
\item
for any $\eps>0$,
there exist points $\omega_i\in\Omega_i\cap B_{\eps}(\bx)$ and vectors $x^*_i\in X^*$ $(i=1,\ldots,n)$
such that conditions \eqref{D4} are satisfied with $\al':=\frac{\al}{1-\al}$ in place of $\al$.
\end{enumerate}
\end{corollary}

\begin{Proof}
In each of the assertions,
let a number $\eps>0$ be given.
For parts (i) and (iii), choose a number $\xi\in]0,\frac{\eps}{2+\al}[$, while for parts (ii) and (iv), choose a number $\xi\in]0,\frac{\eps}{2-\al}[$.
By Definition~\ref{D2.17}, there exist a number $\rho\in]0,\xi[$, points $\omega_i\in \Omega_i\cap B_\xi(\bx)$ and vectors $a_i\in X$ $(i=1,\ldots,n)$ such that conditions \eqref{P7} are satisfied.
In parts (i) and (iii), we have $\xi+\rho(1+\al)<\xi(2+\al)<\eps$, and in parts (ii) and (iv), we have $\xi+\rho(1-\al)<\xi(2-\al)<\eps$.
The conclusions follow from Proposition~\ref{P5.16}.
\qed\end{Proof}

\if{
\AK{1/07/18.
Corollary~\ref{C5.9} yields the implication (i) \folgt (iii) in Theorem~\ref{T2.9}.
We do not seem to have a statement which would immediately give the implication (i) \folgt (ii).}
\HB{We can formulate corollaries after Theorem 5.4}
}\fi

\begin{Remark}
The infinitesimal statements in Corollary~\ref{C5.9} are crucial for the extended extremal principle and its extensions to infinite collections of sets.
For instance, in view of Proposition~\ref{P2.18}, Corollary~\ref{C5.9}(iii) and (iv)
immediately yield the implications, respectively, (i) \folgt (ii) and (i) \folgt (iii) in Theorem~\ref{T2.9}.
\end{Remark}

The dual characterizations of the `extremal' arrangements of collections of sets given in the statements in the first part of this section are themselves in a sense extremal properties of collections of sets.
They can be partially reversed in the setting of a general normed vector space and with \Fr\ normal cones.
We start with an `asymmetric' statement where the last set in the list plays a special role.

\begin{proposition}\label{P5.9}
Suppose $\Omega_1,\ldots,\Omega_n$ are subsets of a normed vector space $X$,
$\omega_i\in\Omega_i$
$(i=1,\ldots,n)$, and $\eps>0$.
If vectors
$x_i^*\in{X}^*$ $(i=1,\ldots,n)$ satisfy
\begin{gather}\label{P5.9-1}
x_i^*\in{N}_{\Omega_i}^F(\omega_i)\;\; (i=1,\ldots,n-1),\quad d(x_n^*,{N}_{\Omega_n}^F(\omega_n))<\eps
\end{gather}
and conditions \eqref{T5.1-3},
then there exists a $\de>0$ such that, for any $\rho\in]0,\de[$ and $\tau\in]0,1[$, there exist vectors
$a_i\in{X}$ $(i=1,\ldots,n-1)$ satisfying conditions \eqref{P6} and
\begin{gather}\label{P5.9-3}
\sum_{i=1}^{n-1}\ang{x_{i}^*,a_i} >\tau\eps\rho.
\end{gather}
\end{proposition}

\begin{Proof}
By \eqref{P5.9-1}, we can choose a vector $x^*\in{N}_{\Omega_n}^F(\omega_n)$ and a positive number $\eps'<\eps$ such that
\begin{gather}\label{P5.9P1}
\norm{x_n^*-x^*}<\eps'.
\end{gather}
Choose also numbers $\eps_1>0$ and $\eps_2>0$ such that \begin{gather}\label{P5.9P0}
n\eps_1+(n-1)\eps_2<\eps-\eps'
\qdtx{and}
(n-1)\eps_2<\eps(1-\tau).
\end{gather}
By the definition of the Fr\'echet normal cone, there exists a number $\de>0$ such that
\begin{gather}\label{P5.9P2}
\langle
x_i^*,\omega-\omega_i\rangle\le
\frac{\varepsilon_1}{\eps+1}\norm{\omega-\omega_i}
\qdtx{for all}
\omega\in\Omega_i\cap
B_{(\eps+1)\de}(\omega_i)\;(i=1,\ldots,n-1),
\\\label{P5.9P2a}
\langle
x^*,\omega-\omega_n\rangle\le
{\varepsilon_1}\norm{\omega-\omega_n}
\qdtx{for all}
\omega\in\Omega_n\cap
B_{\de}(\omega_n).
\end{gather}
Let $\rho\in]0,\de[$.
Choose vectors $a_i\in{X}$ $(i=1,\ldots,n-1)$ such that
\begin{gather}\label{P5.9P3}
\|a_i\|<\eps\rho
\quad\mbox{and}\quad
\langle x_i^*,a_i\rangle> \eps\rho\norm{x_i^*}-\varepsilon_2\rho\quad (i=1,\ldots,n-1).
\end{gather}
Then, by \eqref{T5.1-3} and \eqref{P5.9P0},
\begin{gather*}
\sum_{i=1}^{n-1}\ang{x_{i}^*,a_i} >\eps\rho-(n-1)\varepsilon_2\rho >\tau\eps\rho.
\end{gather*}
Hence, the inequality in (\ref{P6}) and condition \eqref{P5.9-3} are satisfied.
Suppose that the equality in (\ref{P6}) does not hold.
Then there exist $\omega_i'\in\Omega_i$ $(i=1,\ldots,n)$ and an $x\in\rho\B$ such that
\begin{gather*}
\omega_1'-\omega_1-a_1=\ldots =\omega_{n-1}'-\omega_{n-1}-a_{n-1}=
\omega_n'-\omega_n=x.
\end{gather*}
By \eqref{P5.9P3},
$\norm{\omega_i'-\omega_i} =\norm{x+a_i}\le\norm{x}+\norm{a_i}<(\eps+1)\rho <(\eps+1)\de$ $(i=1,\ldots,n-1)$ and $\norm{\omega_n'-\omega_n} =\norm{x}\le\rho <\de$.
Hence, by \eqref{P5.9P2}, \eqref{P5.9P2a} and \eqref{P5.9P3},
\begin{gather*}
\langle x_i^*,x\rangle= \langle x_i^*,\omega_i'-\omega_i\rangle-\langle x_i^*,a_i\rangle< -\eps\rho\norm{x_i^*}+(\varepsilon_1+\eps_2)\rho,
\\
\langle x^*,x\rangle= \langle x^*,\omega_n'-\omega_n\rangle< \varepsilon_1\rho,
\end{gather*}
and consequently, using \eqref{T5.1-3} and \eqref{P5.9P0},
$$
\langle x^*-x_n^*,x\rangle =\sum_{i=1}^{n-1}\langle x_i^*,x\rangle+\langle
x^*,x\rangle< -\eps\rho+n\varepsilon_1\rho+(n-1)\eps_2\rho<-\eps'\rho.
$$
On the other hand, by \eqref{P5.9P1},
$
\langle x^*-x_n^*,x\rangle> -\eps'\rho.
$
A contradiction.
\qed\end{Proof}

The corresponding `symmetric' statement follows immediately after applying Proposition~\ref{P5.9} to the collection of $n+1$ sets $\Omega_1,\ldots,\Omega_n$ and $X$.

\begin{corollary}\label{C5.11}
Suppose $\Omega_1,\ldots,\Omega_n$ are subsets of a normed vector space $X$,
$\omega_i\in\Omega_i$
$(i=1,\ldots,n)$, and $\eps>0$.
If vectors
$x_i^*\in{X}^*$ $(i=1,\ldots,n)$ satisfy conditions \eqref{D1}, where $N$ stands for the \Fr\ normal cone ($N:=N^F$),
then there is a $\de>0$ such that, for any $\rho\in]0,\de[$ and $\tau\in]0,1[$, there exist vectors
$a_i\in{X}$ $(i=1,\ldots,n)$ satisfying conditions \eqref{P3} and
\begin{gather}\label{P5.11-1}
\sum_{i=1}^{n}\ang{x_{i}^*,a_i} >\tau\eps\rho.
\end{gather}
\end{corollary}

In view of Proposition~\ref{P2.7}(iv), the above two statements produce dual sufficient characterizations for approximate stationarity.

\begin{corollary}
Suppose $\Omega_1,\ldots,\Omega_n$ are subsets of a normed vector space $X$ and $\bx\in\cap_{i=1}^n\Omega_i$.
If for any $\eps>0$
there exist points $\omega_i\in\Omega_i\cap B_\eps(\bx)$ $(i=1,\ldots,n)$
and vectors
$x_i^*\in{X}^*$ $(i=1,\ldots,n)$ satisfying either conditions \eqref{P5.9-1} or conditions \eqref{D1} with $N$ standing for the \Fr\ normal cone ($N:=N^F$), then
the collection $\{\Omega_1,\ldots,\Omega_n\}$ is approximately stationary at $\bx$.
\end{corollary}

\begin{Remark}
Similar to the dual necessary characterizations of extremality/stationarity properties discussed in the first part of this section, the sufficient conditions in Proposition~\ref{P5.9} and Corollary~\ref{C5.11} contain conditions \eqref{P5.9-3} and \eqref{P5.11-1}, respectively, relating the given dual vectors $x_i^*$ $(i=1,\ldots,n)$ and the primal space translation vectors $a_i$ $(i=1,\ldots,n)$ guaranteed by the statements.
In view of Remark~\ref{R5.2}(ii), such conditions seem to be an intrinsic feature of the extremality/stationarity properties, independently on whether one goes from primal space conditions to dual space ones or the other way round.
\end{Remark}

\if{
\AK{1/07/18.
Corollary~\ref{C5.11} yields the implication (iii) \folgt (i) in Theorem~\ref{T2.9}.
We do not seem to have a statement which would immediately give the implication (ii) \folgt (i).}
\HB{We do have that by Corollary 5.18}
Combining the statements of Corollaries~\ref{C5.9} and \ref{C5.11}, we arrive at the following equivalent dual space characterization of the condition {\rm (i)} in Corollary~\ref{C4.9}.
}\fi

Combining Corollaries~\ref{C5.9} and \ref{C5.11}, we can formulate a full dual characterization of approximate $\al$-stationa\-rity when either the space is Asplund or the sets are convex.

\begin{corollary}\label{C5.12}
Suppose $\Omega_1,\ldots,\Omega_n$ are closed subsets of a Banach space $X$,
$\bx\in\cap_{i=1}^n\Omega_i$, and $\al>0$.
Suppose also that
either $X$ is Asplund or $\Omega_1,\ldots,\Omega_n$ are convex.
The collection $\{\Omega_1,\ldots,\Omega_n\}$ is approximately $\al$-stationary at $\bx$
if and only if, for any $\eps>0$,
there exist points $\omega_i\in B_{\eps}(\bx)$ and vectors $x^*_i\in X^*$ $(i=1,\ldots,n)$
such that
conditions \eqref{D3} are satisfied with $N$ standing for the \Fr\ normal cone ($N:=N^F$).

Moreover, under the above conditions, if $X$ is Asplund and $\tau\in]0,\frac{1}{1+2\al}[$, or $\Omega_1,\ldots,\Omega_n$ are convex and $\tau:=\frac{1}{1+2\al}$, then, for any $\eps>0$,
points $\omega_i\in B_{\eps}(\bx)$ and vectors $x^*_i\in X^*$ $(i=1,\ldots,n)$ can be chosen to satisfy also condition \eqref{T5.7-0}, while, for any $\hat\tau\in]0,1[$, vectors $a_i\in X$ $(i=1,\ldots,n)$ in Definition~\ref{D2.17} of the approximate $\al$-sta\-tionarity can be chosen to satisfy additionally $\sum_{i=1}^{n}\ang{x_{i}^*,a_i} >\hat\tau\al\rho$.
\end{corollary}

\begin{Proof}
The `only if' part together with the condition \eqref{T5.7-0} in the `moreover' part follow from parts (i) and (iii) of Corollary~\ref{C5.9}, taking into account that for convex sets the Clarke and \Fr\ normal cones coincide.
Conversely, given any $\eps>0$, $\hat\tau\in]0,1[$, points $\omega_i\in B_{\eps}(\bx)$ and vectors $x^*_i\in X^*$ $(i=1,\ldots,n)$
satisfying
conditions \eqref{D3}, Corollary~\ref{C5.11} with $\al$ and $\hat\tau$ in place of $\eps$ and $\tau$, respectively, yields the approximate $\al$-stationarity and condition $\sum_{i=1}^{n}\ang{x_{i}^*,a_i} >\hat\tau\al\rho$.
\qed\end{Proof}

\begin{Remark}
In view of Corollary~\ref{C2.20} and Remark~\ref{R3.9}, the first part of Corollary~\ref{C5.12} yields the statements of Theorems~\ref{T2.9} and \ref{T2.14}.
\end{Remark}

\if{
\AK{1/07/18.
It would be good to avoid resorting to Remark~\ref{R2.11}(i).}
}\fi

The assumption
$x_i^*\in{N}_{\Omega_i}^F(\omega_i)$ $(i=1,\ldots,n-1)$ in Proposition~\ref{P5.9} can be relaxed (at the expense of weakening the estimates in \eqref{P6} and \eqref{P5.9-3}).
The next statement is a consequence of Propositions~\ref{P5.9} and \ref{P2.15}.

\begin{corollary}\label{C5.14}
Suppose $\Omega_1,\ldots,\Omega_n$ are subsets of a normed vector space $X$,
$\omega_i\in\Omega_i$
$(i=1,\ldots,n)$ and $\eps\in]0,1[$.
If vectors
$x_i^*\in{X}^*$ $(i=1,\ldots,n)$ satisfy conditions \eqref{D6} with $N$ standing for the \Fr\ normal cone ($N:=N^F$),
then there is a $\de>0$ such that, for any $\rho\in]0,\de[$ and $\tau\in]0,1[$, there exist vectors
$a_i\in{X}$ $(i=1,\ldots,n-1)$ satisfying conditions \eqref{P6} and \eqref{P5.9-3} with $\eps':=\eps/(1-\eps)$ in place of $\eps$.
\sloppy
\end{corollary}

Since $\eps/(1-\eps)<\eps$, the conclusions of Corollary~\ref{C5.14} are weaker than those of Proposition~\ref{P5.9}.
Observe that conditions \eqref{P6} involve a localization of the $n$th set (near $\omega_n\in\Omega_n$).
The estimates can be improved by considering localizations of all the sets.

\begin{proposition}\label{P5.13}
Suppose $\Omega_1,\ldots,\Omega_n$ are subsets of a normed vector space $X$,
$\omega_{i}\in\Omega_i$ $(i=1,\ldots,n)$, and $\eps>0$.
If vectors
$x_i^*\in{X}^*$ $(i=1,\ldots,n)$ satisfy conditions \eqref{D6} with $N$ standing for the \Fr\ normal cone ($N:=N^F$),
then there is a $\de>0$ such that, for any $\rho\in]0,\de[$ and $\tau\in]0,1[$, there are vectors
$a_i\in{X}$ $(i=1,\ldots,n-1)$ satisfying
\sloppy
\begin{gather}\label{P5.13-1}
\bigcap_{i=1}^{n-1}\Big((\Omega_i-\omega_i)\cap(\rho\B) -a_i\Big)\cap(\Omega_n-\omega_n)\cap(\rho\B) =\emptyset,\quad
\max_{1\le i\le n-1}\norm{a_i}<\eps\rho,
\end{gather}
and condition \eqref{P5.9-3}.
\end{proposition}

The proof below is a modification of that of Proposition~\ref{P5.9}.

\begin{Proof}
Choose vectors $y_i^*\in N_{\Omega_i}^F(\omega_i)$ ($i=1,\ldots,n$) and a positive number $\eps'<\eps$ such that
\begin{equation}\label{P5.13P1}
\sum_{i=1}^{n}\norm{x_i^*-y^*_i}<\eps'.
\end{equation}
Then, choose numbers $\eps_1>0$ and $\eps_2>0$ such that $\eps_1+\eps_2<\eps-\eps'$ and $\varepsilon_2<(1-\tau)\eps$.
By the definition of the Fr\'echet normal cone, there is a $\de>0$ such that
\begin{equation}\label{P5.13P2}
\langle y_i^*,\omega-\omega_i\rangle \le\frac{\eps_1}{n}\norm{\omega-\omega_i}
\qdtx{for all}
\omega\in\Omega_i\cap B_{\de}(\omega_i)\;\; (i=1,\ldots,n).
\end{equation}
Let $\rho\in]0,\de[$.
Choose vectors $a_i\in{X}$ $(i=1,\ldots,n)$ satisfying  \begin{gather*}
\|a_i\|<\eps\rho
\quad\mbox{and}\quad
\langle x_i^*,a_i\rangle> \eps\rho\norm{x_i^*}-\frac{\varepsilon_2\rho}{n-1}\quad (i=1,\ldots,n-1).
\end{gather*}
By \eqref{D6}, we have
\begin{gather}\label{P5.10P3}
\sum_{i=1}^{n-1}\langle x_i^*,a_i\rangle >\eps\rho-\varepsilon_2\rho.
\end{gather}
The inequality in (\ref{P5.13-1}) and
condition (\ref{P5.9-3}) follow.
Suppose that the equality in (\ref{P5.13-1}) is not satisfied.
Then, there exist points
$\omega_i'\in\Omega_i\cap B_\rho(\omega_i)$
$(i=1,\ldots,n)$ and an $x\in\rho\B$ such that
\begin{gather}\label{P5.10P4}
\omega_1'-\omega_1-a_1=\ldots =\omega_{n-1}'-\omega_{n-1}-a_{n-1}=
\omega_n'-\omega_n=x.
\end{gather}
Hence, making use of \eqref{P5.10P4}, \eqref{D6}, \eqref{P5.10P3}, \eqref{P5.13P2} and \eqref{P5.13P1}, we have
\begin{align*}
0&=\sum_{i=1}^{n-1}\langle x_i^*,(\omega_n'-\omega_n) -(\omega'_i-\omega_i-a_i)\rangle
\\
&=-\sum_{i=1}^{n-1}\langle x_i^*,\omega'_i-\omega_i-a_i\rangle-\langle x_n^*,\omega_n'-\omega_n\rangle
\\
&=\sum_{i=1}^{n-1}\langle x_i^*,a_i\rangle
-\sum_{i=1}^{n}\langle x_i^*,\omega_i'-\omega_i\rangle
\\
&=\sum_{i=1}^{n-1}\langle x_i^*,a_i\rangle
-\sum_{i=1}^{n}\langle y_i^*,\omega_i'-\omega_i\rangle
+\sum_{i=1}^{n}\langle y_i^*-x_i^*,\omega_i'-\omega_i\rangle
\\
&>\eps\rho-\varepsilon_2\rho
-\frac{\eps_1}{n} \sum_{i=1}^{n}\norm{\omega_i'-\omega_i}
-\eps'\max_{1\le i\le n}\norm{\omega_i'-\omega_i}
\\
&\ge(\eps-\eps'-\eps_1-\varepsilon_2)\rho>0.
\end{align*}
This contradiction proves the proposition.
\qed\end{Proof}

\if{
\AK{1/07/18.
Where the last two statements could be used?}
\HB{This one is a dual approach for Proposition 4.6 by the sense that implications \eqref{P2} implies \eqref{P4.3-2}, we also have \eqref{P2} implies \eqref{T5.1-3},\eqref{T5.1-4},\eqref{T5.1-5} with $\Omega_n=B_\rho(\bx)$, and those conditions imply \eqref{P4.3-2}.}
}\fi

\section{Conclusions}
We have exposed, analysed and refined
the core arguments used in various proofs of the {extremal principle} and its extensions as well as in primal and dual characterizations of the {approximate stationarity} and {transversality} of collections of sets, presenting a unifying theory, encompassing all existing approaches to obtaining `extremal' statements.
For that, we have examined and clarified quantitative relationships between the parameters involved in the respective definitions and statements.
Some new characterizations of extremality properties have been obtained.

\begin{acknowledgements}
The research was supported by the Australian Research Council, project DP160100854.
Hoa T. Bui is supported by an Australian Government Research Training Program (RTP) Stipend and RTP Fee-Offset Scholarship through Federation University Australia.
Alexander Y. Kruger benefited from the support of the FMJH Program PGMO and from the support of EDF.

We wish to thank PhD student Nguyen Duy Cuong from Federation University Australia for careful reading of the manuscript and helping us eliminate numerous typos, and the anonymous referees for their constructive comments and suggestions.
\end{acknowledgements}

\section*{Conflict of Interest}
The authors declare that they have no conflict of interest.

\end{document}